\documentclass[16pt]{article}
\usepackage{graphicx}
\usepackage{amssymb}

\newtheorem{theorem}{Theorem}[subsection]
\newtheorem{proposition}{Proposition}[subsection]
\newtheorem{corollary}{Corollary}[subsection]
\newtheorem{lemma}{Lemma}[subsection]

\newtheorem{definition}{Definition}[subsection]

\begin{document}

\title{Michel  Kervaire  work on  knots in higher dimensions}

\author{Fran\c coise Michel  and Claude Weber }
\date {GENEVE, 28 aout 2014}

%Si on ne tape rien dans l'une des trois accolades ci-dessus, rien n'apparait dans le texte.%

\maketitle

\begin{abstract} 
The aim of this paper is to present Michel Kervaire's work on differential  knots  in  higher dimensions  in codimension $q = 2$. In order to appreciate the importance of Kervaire's contribution, we  describe in  Sections 2 to 4 what was, at the time,  the situation in differential topology and in knot theory in codimension $q \geq 3$. In Section 5, we expose  Michel Kervaire's characterization of  the fundamental group of a knot complement. In Section 6,  we explain  Kervaire  and Levine's  work  on knot modules. In Section 7,  we detail Kervaire's construction of the ``simple knots"  classified by Jerome Levine.  Section 8 summarizes   Kervaire and Levine's results on knot cobordism. In Section 9,  we apply higher dimensional  knot theory  to singularities of complex hypersurfaces. In the Appendix, Chapters 10 to 13 are devoted to a discussion of some basic concepts, known to the experts: Signs, Seifert Hypersurfaces, Open Book Decompositions and Handlebodies. In Chapter 14, we conclude this paper with an exposition of the results of Hill-Hopkins-Ravenel on the Kervaire Invariant Problem and its consequences to knot theory in codimension two.
\\
The point of view adopted in this paper is kind of pseudo-historical. When we make explicit Kervaire's work we try to follow him closely, in order to retain some of the flavor of the original texts. When necessary, we add further contributions often due to Levine. We also propose developments which occurred later.  

\vspace{6pt}

COMMENTS, REMARKS, COUNSELS, CRITICISMS WILL BE GREATLY APPRECIATED. SEND THEM TO:
Claude.Weber@unige.ch OR TO: fmichel@picard.ups-tlse.fr

\end{abstract}

\vskip1in

\newpage

\begin{center}
\LARGE{Table of contents}
\end{center}

{\bf 1 Introduction}

1.1 Michel Kervaire's six papers on knot theory
\\
1.2 A brief description of the content of this paper
\\
1.3 What is a knot?
\\
1.4 Final remarks
\\
1.5 Conventions and notations

{\bf 2 Some tools of differential topology}

2.1 Surgery from an elementary point of view
\\
2.2 Vector bundles and parallelisability
\\
2.3 Pontrjagin method and the J-homomorphism

{\bf 3 Kervaire-Milnor study of homotopy spheres }

3.1 Homotopy spheres
\\
3.2 The groups $\Theta^n$ and $bP^{n+1}$
\\
3.3 Kervaire invariant
\\
3.4 The groups $P^{n+1}$
\\
3.5 Kervaire manifold

{\bf 4 Differential knots in codimension $\geq 3$}

4.1 On the isotopy of knots and links in any codimension
\\
4.2 Embeddings and isotopies in the stable and metastable ranges
\\
4.3 Below the metastable range

{\bf 5 The fundamental group of a knot complement}

5.1 Homotopy n-spheres in $S^{n+2}$
\\
5.2 Necessary conditions for a group to be the fundamental group of a knot complement
\\
5.3 Sufficiency of the conditions if $n \geq 3$
\\
5.4 Kervaire conjecture
\\
5.5 Groups which satisfy Kervaire conditions

{\bf 6 Knot modules}

6.1 The knot exterior
\\
6.2 Some algebraic properties of knot modules
\\
6.3 The q-th knot module when $q < n / 2$
\\
6.4 Seifert hypersurfaces
\\
6.5  Odd dimensional  knots and the Seifert form
\\
6.6 Even dimensional  knots and the torsion Seifert form  

\newpage

{\bf 7 Odd dimensional  simple links}

7.1  The q-handlebodies
\\
7.2 The realization theorem for Seifert matrices
\\
7.3 Levine's classification of embeddings of handlebodies in codimension one

{\bf 8 Knot cobordism}

8.1 Definitions
\\
8.2 The even dimensional  case
\\
8.3 The odd dimensional  case

{\bf 9 Singularities of complex hypersurfaces}

9.1 The theory of Milnor
\\
9.2 Algebraic links and Seifert forms
\\
9.3 Cobordism of algebraic links
\\
9.4 Examples

{\bf 10 Appendix I: Linking numbers and signs}

10.1 The boundary of an oriented manifold
\\
10.2 Linking numbers

{\bf 11 Appendix II:  Existence of Seifert hypersurfaces}

{\bf 12 Appendix III: Open book decompositions}

12.1 Open books
\\
12.2 Browder's Lemma 2

{\bf 13 Appendix IV: Handlebodies}

13.1 Bouquets of spheres and handlebodies
\\
13.2 Parallelisable handlebodies
\\
13.3 m-dimensional spherical links in $S^{2m+1}$

{\bf  14 Appendix V: Homotopy spheres embedded in codimension two and 
\\
the Kervaire-Arf-Robertello-Levine invariant}

14.1 Which homotopy spheres can be embedded in codimension two?
\\
14.2 The Kervaire-Arf-Robertello-Levine invariant
\\
14.3 The Hill-Hopkins-Ravenel result and its influence on the KARL invariant

{\bf Bibliography}

\newpage

\section{Introduction}

\subsection{Michel Kervaire's six papers on knot theory}

Michel Kervaire wrote six papers on knots: \cite{kerv63}, \cite{kerv65}, \cite{kerv70}, \cite{hake78-1}, \cite{hake78-2} and \cite{kewe77}.

The first one is both an account of Michel Kervaire's talk at the Symposium held at Princeton in Spring 1963 in the honor of Marston Morse and a report of discussions between several participants of the Symposium about Kervaire's results. Reading this paper is really fascinating, since one sees higher dimensional  knot theory emerging. The subject is the fundamental group of knots in higher dimensions. 

The second one is the written thesis that Michel Kervaire presented in Paris in June 1964.  In fact Kervaire had already obtained a PhD in Zurich under Heinz Hopf in 1955, but  he applied in 1964 for a position in France and, at that time, a French thesis was compulsory.  Finally the appointment did not materialize. But  the thesis text remains as an article published in `` Bulletin de la S.M.F." (\cite{kerv65}).   As this article  is  the main reason for writing   this text, we name it  Kervaire's Paris paper. It is the more important that Michel Kervaire wrote on  knot theory. It can be considered to be the foundational text on knots in higher dimensions together with contemporary papers by Jerry Levine \cite{levi65-1} , \cite{levi65-2} and \cite{levi66}. One should also add to the list the Hirsch-Neuwirth paper \cite{hine64}, which seems to be a development of discussions held during the Morse Symposium. 

To briefly present the subject of   Kervaire's  Paris paper,  we need a few definitions. A knot $K^n \subset S^{n+2}$ is the image of a differentiable embedding of a n-dimensional  homotopy sphere in $S^{n+2}$. Its exterior $E(K)$ is the complement of an open tubular neighborhood. The exterior has the homology of the circle $S^1$ by Alexander duality. In short, the subject is the determination of the first homotopy group  $\pi_q(E(K))$ which is different from $\pi_q(S^1)$. 

Later, Michel Kervaire complained that he had to rush for completing this  Paris paper in due time and he had doubts about the quality of its redaction. In fact we find this article well-written. The exposition is concise and clear, typically in Kervaire's style. The pace is slow in parts which present a difficulty and fast when things are obvious. Now, what was obvious to Kervaire in Spring 1964? Among other things, clearly Pontrjagin's construction and surgery. As these techniques are possibly not so well known to a reader fifty years later, we devote Sections 2 and 3 of this text to a presentation of these matters. 

The last chapter of Kervaire Paris paper is a first attempt at the understanding of the cobordism of knots in higher dimensions.  It contains a complete proof that an  even dimensional knot is  always cobordant to a trivial knot. For the odd dimensional knots, the subject  has a strong algebraic flavor, much related to quadratic forms, their isometry group and algebraic number theory. Kervaire liked this algebraic aspect and devoted his third paper (Amsterdam meeting in 1970) to it.

Jerry Levine spent the first months of 1977 in Geneva and he gave wonderful lectures on many aspects of knot theory. His presence had a deep influence on several members of the audience, including the two authors of this paper! Under his initiative a meeting was organized by Kervaire in Les-Plans-sur-Bex in March 1977. The proceedings are recorded in the Springer Lecture Notes volume 685, edited by Jean-Claude Hausmann. The irony of history is that it is precisely at this time that William Thurston made his first announcements, which soon completely shattered classical knot theory and 3-manifolds. See  p.44 of Cameron Gordon's paper in the Proceedings. The meeting renewed Kervaire's interest in knot theory. In the following months he wrote his last three papers on the subject.

\subsection{A brief description  of the content of this paper}

As promised, we devote Sections 2 and 3 to some background in differential topology. Mainly: vector bundles, Pontrjagin Construction, surgery,  culminating with Kervaire-Milnor. 

We felt necessary to devote Section 4 to knots in codimension $\geq 3$. Its reading is optional. One reason to do so is that the subject was flourishing at the time, thanks to the efforts of Andr\'e Haefliger and Jerry Levine. It is remarkable that Levine was present in both fields. Another reason is that the two subjects are in sharp contrast. Very roughly speaking one could say that it is a matter of fundamental group. In codimension $\geq 3$ the fundamental group of the exterior is always trivial while it is never so in codimension 2.  But more must be said. In codimension $\geq 3$ there are no PL knots, as proved by Christopher Zeeman  \cite{zeem63}. Hence everything is a matter of comparison between PL and DIFF. This is the essence of Haefliger's theory of smoothing, written a bit later. On the contrary, in codimension 2 the theories of PL knots and of DIFF knots do not much differ. 

In Section 5 we present Michel Kervaire's determination of the fundamental group of knots in higher dimensions, together with some of the results of his two papers written with Jean-Claude Hausmann about the commutator subgroup and the center of these groups. Some later developments are also presented. 

In Section 6 we expose Michel Kervaire's results on the first homotopy group $\pi_q(E(K))$ which is different from $\pi_q(S^1)$. For $q \geq 2$ these groups are in fact ${\bf Z}\lbrack t , t^{-1} \rbrack$-modules. Indeed,  Kervaire undertakes a first study of such modules, later to be called knot modules by Levine. In our presentation, we include several developments 
due to Levine. 

Up to  Kervaire's  Paris paper, most of the efforts in knot theory went to the construction of knot invariants. They produce necessary conditions for two knots to be equivalent. In Levine's paper \cite{levi65-2} a change took place. From Dale Trotter's work it was known that for classical knots,  Seifert matrices of equivalent knots are S-equivalent. Levine introduced a class of odd-dimension knots (called by him simple knots) for which the S-equivalence of the Seifert matrices is both necessary and sufficient for two knots to be isotopic. This is a significant classification result. Indeed simple knots are already present in Kervaire Paris paper, but he did not pursue their study that far. Technically,  the success of Levine's study is largely due to the fact that these knots bound a very special kind of Seifert hypersurface: a (parallelisable) handlebody. In Section 7 we present Levine's work on odd dimension simple knots.  

Section 8 is devoted to higher dimension knot cobordism. Levine reduced the determination of these groups to an algebraic problem. A key step in the argument rests on the fact that each knot is cobordant to a simple knot. 

In knot theory, the handlebodies one deals with are parallelisable and their boundary is a homotopy sphere. If we keep the  parallelisability condition  but admit any boundary, the Kervaire-Levine's  arguments are still valid. This immediately applies to the Milnor fiber of isolated singularities of complex hypersurfaces, as was first noticed by Milnor himself and developed by Alan Durfee. Section 9 is devoted to that matter. It can be considered as a posterity of Kervaire.

Sections 10 to 14 are a kind of appendix which  can be read independently. They provide basics, comments, variations on subjects treated elsewhere in our paper. In Section 10 we expose our conventions on signs which agree  with Kauffman-Neumann \cite{kane77}. This allows us to justify   the signs for  invariants of some basic  algebraic  links (examples given at  the end of Section 9).  Often, authors do not mention their signs conventions and  hence one can find other  signs in the literature. In Section 11 we prove the existence of  Seifert hypersurfaces in a more general context. In Section  12 and Section 13,  we present basics about open book decompositions and  parallelisable handlebodies which are useful in knot theory. Our aim in Section 14 is to present the beautiful result of  Mike Hill, Mike Hopkins and Doug Ravenel about the Kervaire Invariant and to tell how this affects the  theory of knots in higher dimensions. It is a spectacular way to conclude this paper.

\subsection{What is a knot ?}

\begin{definition}
A n-link  in $S^{n+q}$ is a compact oriented differential submanifold without boundary $L^n \subset S^{n+q}$. The integer $q \geq 1$ is the codimension of the link. The n-links $L^n_1$ and $L^n_2$ are equivalent if there exists a diffeomorphism $f : S^{n+q} \rightarrow S^{n+q}$ such that $f(L^n_1) = L^n_2$, respecting the orientation of $S^{n+q}$ and of the links.
\\
When $L^{n}$ is a homotopy sphere, we say that  $L^n \subset S^{n+q}$ is a n-knot in codimension $q$. 
\end{definition}

In the beginning of Section 4, we present the well-known proof that if two links are equivalent there always exists a diffeomorphism $f : S^{n+q} \rightarrow S^{n+q}$ which is isotopic to the identity and moves one link  to the other. Hence n-links are equivalent if and only if they are isotopic. 

In general the boundary of a Milnor fiber is not a homotopy sphere. It is a motivation to 
  explain, in Section 7,  how   Kervaire and Levine's works on simple odd-dimensional  knots can be generalized to simple links.

  A  link $ L^n \subset S^{n+2}$  in codimension two is always the boundary of a $n+1-$dimensional oriented smooth submanifold  $F^{n+1}$ in  $S^{n+2}$. We say that  $F^{n+1}$ is a {\it Seifert hypersurface } for  $ L^n \subset S^{n+2}$ ( some authors name it Seifert surface even when $n\geq 2$).

\subsection{Final remarks}
The aim of this paper is to pay tribute to Michel Kervaire and to make    his  work on knots of higher dimensions easier to read by  younger generations of mathematicians. Basically it is a mathematical   exposition paper. Our purpose is not to write an history of knots in higher dimensions. We apologize for not making  a list of all papers in the subject.
 When we present Kervaire's work we try to follow him closely, in order to retain some of the flavor of the original texts. When necessary, we add further contributions often due to Levine. We also propose developments which occurred later.
 With the passing of time, we  find important to present in  details results on the fundamental group of the knot complement  and on simple odd dimensional knots(and links). Indeed:
 \\
  1) The  determination of the fundamental group of the knot complement played a  key role in the beginning of higher dimensional knots theory.
  \\
  2) The higher odd  dimensional simple knots  can be classified via their relations with handelbodies. In one hand this  classification induces a  classification up to cobordism. On the other hand, it can be easily generalized  to links associated to isolated singular point of complex hypersurfaces.

We have wondered  whether to include Jerry Levine's name in the title of the paper. We have decided not to, although he certainly is the cofounder of higher dimension knot theory. But Levine pursued his work much beyond these first years, while Kervaire stopped publishing in the subject (too) early. Hence it would have been difficult to find an equilibrium between them. In fact a study of Levine's work in knot theory should be much longer than this  paper. 

\subsection{Conventions and notations}

 Manifolds and embeddings are $\mathcal{C}^{\infty}$. Usually, manifolds are compact and oriented. A manifold is {\bf closed} if compact without boundary. The boundary of $M$ is written $bM$, its interior is   $\mathring{M} $ and its closure is $\bar{M}.$
Let $L$ be a closed oriented submanifold of a closed  oriented manifold $M$.  We denote by $N(L)$ a closed tubular neighbourhood of $L$ in $M$ and by $E(L)$ the closure of $M\setminus N(L),$ i.e., $E(L)=M\setminus \mathring {N}(L).$ By definition $E(L)$ is  the {\bf  exterior } of $L$ in $M.$

Fibres of vector bundles are vector spaces  over the field of real numbers ${\bf R}$. The {\bf rank } of a vector bundle over a connected base is the dimension of its fibres.

\vskip.5in

\section{Some tools of  differential topology}

An extremely useful reference for this  section and the next one is provided by Andrei Kosinski's book~\cite{kosi93}. Basics about differential manifolds are beautifully presented by Morris Hirsch in \cite{hirs76}.

\subsection{Surgery from an elementary point of view}

Originally, surgery was ``just" a way to transform a differential manifold into another one as in \cite{miln61}. 

\begin{definition}
The {\bf (standard) sphere} $S^n$ of dimension $n$ is the set of points $x = (x_0 , x_1 , \dots , x_n) \in {\bf R}^{n+1}$ such that $\sum x_i^2 = 1$. The {\bf ball} $B^{n+1}$ of dimension $(n+1)$ is the set of points $x = (x_0 , x_1 , \dots , x_n) \in {\bf R}^{n+1}$ such that $\sum x_i^2 \leq 1$. The {\bf open ball} $\mathring{B}^{n+1}$ of dimension $(n+1)$ is defined by $\sum x_i^2 < 1$.
\end{definition}

Consider the product of spheres $S^a \times S^b ~~a \geq 0 ~, ~b \geq 0$. This manifold is the boundary of $S^a \times B^{b+1}$ and of $B^{a+1} \times S^b$. Suppose now that we have $(S^a \times B^{b+1}) \subset M'$ where $M'$ is a manifold of dimension $m = a + b + 1$. If $M'$ has a boundary we suppose that the embedding is far from $bM'$. We consider then $M = M' \setminus (S^a \times \mathring{B}^{b+1})$ and we construct $M'' = M \cup (B^{a+1} \times S^b)$ with $M \cap (B^{a+1} \times S^b) = S^a \times S^b$.

\begin{definition}
One says that $M''$ is obtained from $M'$ by a {\bf surgery} along $S^a \times \lbrace 0 \rbrace$. The manifold $M$ is the {\bf common part} of $M'$ and $M''$. 
\end{definition}

{\bf Remarks.} 1) The process is reversible: $M'$ is obtained from $M''$ by surgery along $\lbrace 0 \rbrace \times S^b$.
\\
2) $M'$ is obtained from the common part $M$ by attaching the  cell $e^{b+1} = \lbrace u \rbrace \times B^{b+1}$ of dimension $b+1$ and then a cell $e^m$ of dimension $m$.
\\
3) Analogously $M''$ is obtained from $M$ by attaching the cell $e^{a+1} = B^{a+1} \times \lbrace v \rbrace$ of dimension $a+1$ and then a cell  of dimension $m$.  

\begin{definition}
We call $S^a \times \lbrace 0 \rbrace \subset M'$ and $\lbrace 0 \rbrace \times S^b \subset M''$ the {\bf scars} of the surgeries.
\end{definition}

\subsection{Vector bundles and parallelisability}

In order to apply surgery successfully, we need to know when a sphere differentiably embedded in a differential manifold has a trivial normal bundle. Here are some answers.

\begin{definition}
A vector bundle   over a complex $X^k$ is {\bf trivial} if it is isomorphic to a product bundle. A {\bf trivialisation} is such an isomorphism.
\end{definition}

\begin{definition}
A manifold is {\bf parallelisable} if its tangent bundle $\tau M$ is trivial.
\end{definition}

{\bf Comment.} Parallelisable manifolds were introduced and studied by Ernst Stiefel in his thesis \cite{stie35} written under the direction of Heinz Hopf. Both were looking for conditions satisfied by a manifold when it is the underlying space of a Lie group. A trivialisation is called by Stiefel a parallelism. Stiefel proved that,  among spheres, $S^1 , S^3 , S^7$ are parallelisable. The following result was proved by Michel Kervaire \cite{kerv58} and also by Bott-Milnor  \cite{bomi58}. A key step is provided by Bott periodicity theorem \cite{bott59}.

\begin{theorem}
$S^1 , S^3 , S^7$ are the only spheres which are parallelisable.
\end{theorem}

{\bf Notations.}  Let $\epsilon^r$ denote the trivial bundle of rank $r$ over an unspecified basis. 

\begin{definition}
A vector bundle $\eta$ is {\bf stably trivial} if $\eta \oplus \epsilon^r$ is trivial for some integer $r \geq 0$. 
\end{definition}

\begin{proposition}
Let $X^k$ be a complex of dimension $k$. Let $\eta$ be a vector bundle over $X^k$ of rank $s$ with $k < s$. Suppose that $\eta$ is stably trivial. Then it is trivial. 
\end{proposition}

{\bf Proof.} See \cite{miln61} Lemma 4.

To be parallelisable is a very strong condition on a manifold, quite difficult to handle. Next condition is easier to work with.  

\begin{definition}
A manifold $M$ is {\bf stably parallelisable} if its tangent bundle is stably trivial. 
\end{definition}

{\bf Comments.} 1) If $M$ is stably parallelisable then $\tau M \oplus \epsilon^1$ is trivial, by the proposition above.
\\
2) The same proposition implies that a  compact connected manifold  $M^n$ with non-empty boundary is stably parallelisable if and only if it is parallelisable, since it has the homotopy type of a complex of dimension $n-1$.  
\\
3) The standard embedding of $S^n$ in ${\bf R}^{n+1}$ shows that $S^n$ is stably parallelisable for all $n$. 
\\
4) More generally,   suppose that $M^n$ is embedded in a stably parallelisable manifold with trivial normal bundle. Then $M^n$ is stably parallelisable.

\begin{proposition}
Let $M^n$ be a stably parallelisable manifold of dimension $n$, embedded in a stably parallelisable manifold $W^N$ of dimension $N$ with $N \geq 2n+1$. Then the normal bundle $\nu$ of $M$ in $W$ is trivial. 
\end{proposition}

{\bf Proof.} Since $W$ is stably parallelisable,  we have $\tau W \oplus \epsilon^1 = \epsilon^{N+1}$. Hence $\tau M \oplus \nu \oplus \epsilon^1 = \epsilon^{N+1}$. Since $M$ is stably parallelisable we have $\nu \oplus \epsilon^{n+1} = \epsilon^{N+1}$. Since the rank of $\nu$ is strictly larger than the dimension of $M$ proposition 1.2.1 applies. 

{\bf Comments.} 1)   In practice,  last proposition is applied to  perform surgery on embedded spheres in stably parallelisable manifolds.
\\
2) The proposition is NOT TRUE in general if $N = 2n$. See below the subsections about Kervaire invariant and about the groups $P^{n+1}$.

\subsection{Pontrjagin method and the J-homomorphism}

Pontrjagin method is also called Thom-Pontrjagin method or construction. See \cite{pont55} and \cite{thom54}. Its importance lies in the fact that it ties together (stably) parallelisable manifolds and (stable) homotopy groups of spheres. Originally, Pontrjagin wished to compute homotopy groups of spheres via differential topology. In the hands of Kervaire and Milnor it provided the crucial link between homotopy groups of spheres and groups of homotopy spheres (the novice reader should breathe deeply and read the last sentence a second time). 

We can see in retrospect that Michel Kervaire's early work as a graduate student was an ideal preparation  for a young topologist in the fifties.   Heinz Hopf asked Michel Kervaire to read Pontrjagin's announcement note \cite{pont50} and to provide proofs where needed. See \cite{kerv56} bottom p.220. Michel Kervaire's description of the starting point of his thesis is quite premonitory: ``Une vari\'et\'e ferm\'ee $M_k$ \'etant plong\'ee dans un espace euclidien ${\bf R}_{n+k}$ avec un champ de rep\`eres normaux..." (A closed manifold $M_k$ being embedded in a euclidean space ${\bf R}_{n+k}$ with a field of normal frames...). See \cite{kerv56} middle p.219. The first chapter of Michel Kervaire thesis  in Zurich is devoted to a detailed presentation of Pontrjagin construction.

Here is a short reminder about Pontrjagin method. We shall only need the stable version of it. 

Let $M^n$ be a closed stably parallelisable manifold. We do not assume that $M^n$ is connected. Choose an integer $m$ such that $m \geq n+2$. By Whitney's theorems \cite{whit36} we can embed $M^n$ in $S^{n+m}$ and any two embeddings are isotopic. By subsection 2.2 its normal bundle is trivial. {\bf We choose} a trivialisation $F$ of it. More precisely, we choose $m$  sections of the normal bundle $e_1 , e_2 , \dots , e_m$ such that $e_1(x) , e_2(x) , \dots , e_m(x)$ are linearly independent for each $x \in M^n$. These sections provide a map $\psi : (U; bU)  \rightarrow (B^m; bB^m)$ where $U$ denotes a closed tubular neighbourhood of $M^n$ in $S^{n+m}$. The sphere $S^m$ is identified with $B^m / bB^m$. Let $\bar{\psi} : U / bU \rightarrow B^m / bB^m = S^m$ be the induced map. Extend $\bar{\psi}$ to a map $\Psi : S^{m+n} \rightarrow S^m$ by sending the complement  of $U$ in $S^{n+m}$ to the smashed boundary $bB^m$. It is easily verified that the homotopy class of $\Psi$ depends only on $(M , F)$. This is the essence of  {\bf Pontrjagin construction}. 

For short we call the couple $(M^n \subset S^{n+m} , F)$ a {\bf framed n-manifold} (in fact a manifold in $S^{n+m}$ with a trivialisation (framing) of its normal bundle). Two framed n-manifolds $M^n_i$ for $i = 0 , 1$ are {\bf framed cobordant} if there exists a framed $W^{n+1} \subset S^{n+m} \times I$ such that $W \cap S^{n+m} \times \lbrace i \rbrace = M_i$ consistently with the framings. Framed cobordism classes constitute an abelian group $\Omega_{n,m}^{fr}$ under disjoint union.

Using transversality techniques one proves that the group $\Omega_{n,m}^{fr}$  is isomorphic to $\pi_{n+m}({S^m})$. Since we have assumed that $m \geq n+2$ the groups  $\pi_{n+m}({S^m})$ and $\Omega_{n,m}^{fr}$ are  independent of $m$. We denote them by $\pi^S_n$ and $\Omega_n^{fr}$. The group $\pi^S_n$ is   the {\bf n-th~stable homotopy group of spheres} (also called the n-th stem).  By Jean-Pierre Serre thesis \cite{serr51}, $\pi_n^S$ is a finite abelian group.

\newpage

In \cite{kerv59} Michel Kervaire gave interpretations of several constructions in the theory of  homotopy groups of spheres via Pontrjagin method. For instance, the  method  yields an easy description of  the Hopf-Whitehead J-homomorphism as follows. First, let us observe that homotopic trivialisations give rise to framed cobordant manifolds and hence produce the same element of $\pi^S_n$. Now take for manifold $M^n$ the n-th sphere $S^n$ standardly embedded in $S^{n+m}$. Its normal bundle  has a standard trivialisation  and hence (homotopy classes of)  different trivialisations are in natural bijection with $\pi_n(SO_m)$.   Pontrjagin construction restricted to $M^n = S^n$ produces a homomorphism $J_n : \pi_n(SO_m) \rightarrow \pi^S_n$. Since we assume that $m \geq n+2$ the group $\pi_n(SO_m)$ does not  depend on $m$ and is usually denoted by $\pi_n(SO)$. The homomorphism $J_n : \pi_n(SO) \rightarrow \pi^S_n$ is the {\bf stable J-homomorphism} in dimension $n$. 

Just at the right time Raoul Bott \cite{bott59} computed the groups $\pi_n(SO)$. His celebrated results {\bf (Bott periodicity)} are as follows:

\begin{theorem}
1) $\pi_n(SO)$ depends only on the residue class of $n ~~ mod~8$.
\\
2) $\pi_n(SO)$ is isomorphic to ${\bf Z} / 2$ if $n \equiv 0$ or $1 ~~mod~8$ (of course we assume $n > 0$).  
\\
3) $\pi_n(SO)$ is isomorphic to ${\bf Z}$ if $n \equiv 3$ or $7 ~~mod~8$.
\\
4) $\pi_n(SO)$ is equal to 0 in the four  other cases $n \equiv  2 , 4 , 5 , 6 ~~mod ~8$. 
\end{theorem}

Frank Adams' important results  \cite{adam65-2} about the stable J-homomorphism are stated in the next two theorems. The first one is indispensable to prove that a homotopy sphere of dimension $n$ is stably parallelisable when $(n-1)$ is congruent to $0$ or $1 \ mod \ 8$. See Subsection 3.2.

\begin{theorem}
The homomorphism $J_n$ is injective when  $\pi_n(SO)$ is isomorphic to ${\bf Z} / 2$.  
\end{theorem}

When $\pi_n(SO)$ is isomorphic to ${\bf Z}$ the image $ImJ_n \subset \pi^S_n$ is a finite  cyclic subgroup whose elements are by their very construction easy to describe. Hence it was important to determine what this subgroup is. Milnor-Kervaire in \cite{mike58} (with the help of Atiyah-Hirzebruch \cite{athi59} to get rid of a possible factor 2) gave a ``lower bound" (i.e. a factor) of its order. Then Adams proved that this ``expected value" is the right one. Another possible factor 2 was eliminated by the proof of the Adams conjecture. See \cite{adam70} p.529-532. The proof of Adams conjecture  also established  that $ImJ_n$ is a direct factor of $\pi^S_n$. The final statement is: 

\begin{theorem}
Let us write $n = 4k - 1$. Then the order of the image of  $J_{4k-1} : \pi_{4k-1}(SO) \rightarrow \pi^S_{4k-1}$ is equal to the denominator den$(B_k / 4k)$ where $B_k$ is the k-th Bernoulli number, indexed as by Friedrich Hirzebruch in~\cite{hirz66}.  
\end{theorem}

It is fortunate for topologists  that this denominator is computable (it is the ``easy part" of Bernoulli numbers, much related to von Staudt theorems).  See \cite{adam65-1} for explicit formulas.

\vskip.5in

\section{Kervaire-Milnor study of homotopy spheres}

\subsection{Homotopy spheres}

\begin{definition}
A differential  closed manifold $\Sigma^n$ of dimension $n$ is a { \bf homotopy sphere} if it has the homotopy type of the standard sphere $S^n$. We shall always assume that homotopy spheres are oriented.
\end{definition}

{\bf Comments.} 1) Classical results of algebraic topology imply that $\Sigma^n$ is a homotopy sphere if and only if $\pi_1(\Sigma^n) = \pi_1(S^n)$ and $H_i(\Sigma^n ; {\bf Z}) = H_i(S^n ; {\bf Z})$ for all $i = 0 , 1 , \dots , n$. 
\\
2) When \cite{kemi63} was written it was just known (thanks to Stephen Smale  and John Stallings) that a homotopy sphere of dimension $n \geq 5$ is in fact homeomorphic to the standard sphere $S^n$. See \cite{smal60} and \cite{stal60}.
\\
3) It is known today that  a homotopy sphere of dimension $n$ is homeomorphic to $S^n$ for any  value of $n$. In  fact a stronger result is known. Consider a compact and connected topological n-manifold which has the homotopy type of the n-sphere. Then this manifold is homeomorphic to $S^n$. We could say that: ``the topological Poincar\'e conjecture is true in all dimensions". Stated in this form, the result is due to Newman \cite{newm66} for $n \geq 5$ and to Freedman \cite{free82} for $n = 4$.  The proof is classical for $n \leq 2$, but it needs the triangulation of surfaces. For $n = 3$ the proof is tortuous. Apply first Moise \cite{mois52} to triangulate the homotopy sphere. Then smooth it by Munkres \cite{munk60}. Finally apply Perelman!

\begin{proposition}
A homotopy sphere of dimension $n \geq 5$ is diffeomorphic to the standard sphere if and only if it bounds a contractible manifold.
\end{proposition}

{\bf Proof.} This is an immediate consequence of Smale h-cobordism theorem. 

\subsection{The groups $\Theta^n$ and $bP^{n+1}$}

 We consider the set $\Lambda^n$ of oriented homotopy spheres up to orientation preserving diffeomorphism. On $\Lambda^n$  a composition law is defined, called {\bf connected sum} and written $\Sigma_1 \sharp \Sigma_2$. See \cite{miln59-2} for a comprehensive study of this operation. It is commutative, associative and has the standard sphere as zero element. In other words, it is a commutative monoid. By the validity of the topological Poincar\'e conjecture, $\Lambda^n$ classifies the oriented differential structures on the n-sphere up to orientation preserving diffeomorphism. To avoid confusing readers who go through papers of the fifties and sixties, we keep the terminology ``homotopy spheres". 

\begin{definition}
Let $W^{n+1}$ be an oriented compact (n+1)-manifold and let $M^n_1$ and $M^n_2$ be two oriented n-manifolds. Then $W^{n+1}$ is an oriented h-cobordism between $M^n_1$ and $M^n_2$ if:
\\
1) the oriented boundary $bW$ is diffeomeorphic to $M^n_1 \coprod -M^n_2$, where the $-$sign denotes the opposite orientation, 
\\
2) both inclusions $M^n_i \hookrightarrow W^{n+1}$ are homotopy equivalences.
\end{definition}

Kervaire and Milnor prove that the h-cobordism relation is compatible with the connected sum operation. The quotient of $\Lambda^n$ by the the h-cobordism equivalence relation is an abelian group written $\Theta^n$. The inverse of $\Sigma$ is $-\Sigma$.

{\bf Comments.} 1) If $n \geq 5$ by Smale's h-cobordism theorem $\Theta^n$ is in fact isomorphic to $\Lambda^n$. Note that Kervaire-Milnor say explicitly that their paper does not depend on Smale's results. But of course Smale's results are needed to interpret $\Theta^n$ in terms of differential structures  on $S^n$. 
\\
2) Since homotopy spheres are oriented a chirality question is present here. A homotopy sphere $\Sigma^n$ is achiral (i.e. possesses an orientation reversing diffeomorphism) if and only if it represents an element of order $\leq 2$ in $\Theta^n$ (assume $n \neq 4$).  

\newpage

The starting point of \cite{kemi63} is the following result.

\begin{theorem}\label{T:hompar}
Any homotopy sphere  $\Sigma^n$ is  stably parallelisable. 
\end{theorem}

{\bf Rough idea of the proof.}  We first need a definition.

\begin{definition}
A closed connected differential manifold $M^n$ is {\bf almost parallelisable} if $M^n \setminus \lbrace x \rbrace$  or equivalently $M^n \setminus \mathring {B}^n$ is parallelisable.
\end{definition}

Clearly a homotopy sphere is almost parallelisable. For an almost parallelisable manifold $M^n$ there is one obstruction to stable parallelisability  which  is an element of $\pi_{n-1}(SO)$.  By Pontrjagin construction this obstruction lies in the kernel  of the J-homomorphism. If $\pi_{n-1}(SO)$ is isomorphic to ${\bf Z} / 2$ by Adams theorem the J-homomorphism is injective in these dimensions and hence the obstruction vanishes.  If $\pi_{n-1}(SO)$ is isomorphic to ${\bf Z}$ a non-zero multiple of  the obstruction is equal to the signature of $M^n$ and hence equal to 0 if $M^n$ is a homotopy sphere. 

{\bf Comment.} In fact the argument proves the following result: An almost parallelisable n-manifold $M^n$ is stably parallelisable if $n \equiv 1 , 2 , 3 $ mod 4. If $n \equiv 0$ mod 4 then $M^n$ is stably parallelisable if and only if its signature vanishes.

Kervaire and Milnor can then apply Pontrjagin method to elements of $\Theta^n$ as follows. Let $\Sigma^n$ represent an element $x$ of $\Theta^n$. Embed $\Sigma^n$ in $S^{n+m}$ for $m$ large $(m \geq n+2)$. By Proposition 2.2.2 its normal bundle $\nu$ is trivial. {\bf We choose a trivialisation of} $\nu$. Pontrjagin method produces then an element $f_n(x)$ of $\pi_n^S$ which does not depend on the choice of $\Sigma^n$ to represent $x$ nor on the embedding of $\Sigma^n$ in $S^{n+m}$. It depends however on the choice of the trivialisation of $\nu$. The indeterminacy lies in the subgroup $ImJ_n$ of $\pi_n^S$. We obtain therefore a homomorphism 

$$f_n : \Theta^n \rightarrow \pi_n^S/ImJ_n = Coker J_n$$

The main objective of \cite{kemi63} is to determine the kernel and cokernel of this homomorphism by surgery. 

Kervaire-Milnor results are the following. First the homomorphism $f_n$ is almost always surjective.  More precisely:

\begin{theorem}(Kervaire-Milnor)
1) For  $n \equiv ~0 , 1 , 3~~mod ~4$ the homomorphism $f_n : \Theta^n \rightarrow Coker J_n$ is surjective.
\\
2) For $n \equiv ~2 ~~mod ~4$ there is a homomorphism $KI_n : \pi^S_n \rightarrow {\bf Z} / 2$ called the {\bf Kervaire invariant} such that  $f_n : \Theta^n \rightarrow Coker J_n$ is surjective if and only if $KI_n$ is the trivial homomorphism. 
\end{theorem}

The definition and properties of Kervaire invariant are postponed to next subsection.

{\bf Comment.} Usually the J-homomorphism $J_n$ is not surjective. But, since $f_n$ is almost surjective, we can often represent an element of $\pi_n^S$ by an embedded homotopy sphere, instead of the standard sphere.

The kernel of $f_n$ is very interesting. It is a subgroup of $\Theta^n$ denoted by $bP^{n+1}$ in \cite{kemi63}. The authors prove the following facts. The proofs are quite non-trivial and make a heavy use of surgery.

\begin{theorem}
The group $bP^{n+1}$ is a finite cyclic group. More precisely: 
\\
1) It is the trivial group if $n+1$ is odd. 
\\
2) If $n+1 \equiv 2~mod~4$ it is either $0$ or ${\bf Z}/2$. More precisely $bP^{4k+2} = {\bf Z} / 2$ if and only if $KI_{4k+2}$ is trivial. 
\\
3) If $n+1 = 4k$ then $bP^{n+1}$ is non-trivial. Its order  is equal to 

$$   2^{2k-2}(2^{2k-1} - 1) num(4B_k / k) $$
\end{theorem}

{\bf Comment.}  The numerator of $4B_k / k$ is the ``hard part" of Bernoulli numbers. It is a product of irregular primes (Kummer results about Fermat conjecture) and tends to infinity at a vertiginous speed. See \cite{bern??}.

Summing up we obtain Kervaire-Milnor short (indeed not quite short) exact sequence (where we extend the definition of $KI_n$ to other values of $n$ to be the 0-homomorphism):

$$0 \rightarrow bP^{n+1} \rightarrow \Theta^n \rightarrow \pi_n^S / ImJ_n \rightarrow ImKI_n \rightarrow 0$$

This exact sequence is valid for all values of $n \geq 1$. In fact for $n \leq 6$ one has $\Theta^n = 0 = bP^{n+1}$. This is due to Kervaire-Milnor for $n \geq 4$ and Perelman for $n = 3$. Moreover one has $ \pi_n^S / ImJ_n = 0 = ImKI_n$ for $n = 1 , 3 , 4, 5$ and $ \pi_n^S / ImJ_n = {\bf Z} / 2 = ImKI_n$ for $n = 2 , 6$. Typically, in these last two dimensions, the non-trivial element in $\pi_n^S / ImJ_n$ is represented by a framing on $S^1 \times S^1$ and on $S^3 \times S^3$. 
\\
 Serious affairs begin with $n = 7$.

For us an important fact is expressed in the following remark. It is an immediate consequence of Pontrjagin method.

{\bf Remark.} A homotopy sphere $\Sigma^n$ represents an element of $bP^{n+1}$ if and only if it bounds a parallelisable manifold (since such a manifold has necessarily non-empty boundary ``parallelisable" is equivalent to ``stably parallelisable"). 

{\bf Comments.} 1) Kervaire-Milnor short exact sequence  relates $\Theta^n$ with  two objects which are in the heart of mathematics. Both are quite difficult to compute explicitly: 
\\
i) $Coker J_n = \pi_n^S / Im J_n$ is the ``hard part" of $\pi_n^S$. It is known (Adams) that it can be identified with a direct summand of $\pi_n^S$ and is essentially accessible via spectral sequences (Adams, Novikov). 
\\
ii) The hard part of Bernoulli numbers is a mysterious subject. Kummer results imply that it is a product of irregular primes and that every irregular prime appears at least once in the hard part of some $B_k$.  A look at existing tabulations is impressive. See \cite{bern??}.

2) It is known today \cite{brum68} that Kervaire-Milnor short exact sequence splits at $\Theta^n$. Michel Kervaire tried to prove this in 1961-62 (letter to Andr\'e Haefliger dated Jan. 13 1962). 

3) From existing tables of $\pi_n^S$ and $B_k$ (for instance on the Web) one can determine $\Theta^n$ without pain for roughly $n \leq 60$. Very likely, specialists of the 2-primary component of $Coker J_n$ can improve the computations to $n \leq 100$ (and maybe more?). See \cite{rave04} for $n \leq 60$.

{\bf The rival group $\Gamma^n$.}  Historically (and conceptually)   two groups are  in competition: $\Theta^n$ and $\Gamma^n$. The second one was first defined by Thom (1958) in his program to smooth PL manifolds. Its definition is lucidly given by Milnor in \cite{miln59-2}. See also \cite{miln64}. Here it is.  Milnor says that an oriented differential structure  $\Sigma^n$ on the n-sphere is a twisted n-sphere if $\Sigma^n$ admits a Morse function with exactly two critical points. From Morse theory, this is equivalent to say that $\Sigma^n$ is obtained by gluing two closed n-balls along their boundary with an orientation preserving diffeomorphism. Then we consider the subset $\Gamma^n$ of $\Lambda^n$ which consists  in elements represented by twisted spheres. It is a group (abelian) since there is an isomorphism:

$$ \Gamma^n =  \frac {Diff^+(S^{n-1})}{rDiff^+(B^n)}$$

Here $Diff^+(S^{n-1})$ denotes the group of orientation preserving diffeomorphisms of $S^{n-1}$ and $Diff^+(B^n)$ the group of orientation preserving diffeomorphisms of the closed ball $B^n$. The letter $r$ denotes the restriction homomorphism. 
\\
Today, it is known that: $\Gamma^n = 0$ for $n \leq 3$ by Smale \cite{smal59} and Munkres \cite{munk60} and $\Gamma^4 = 0$ by Cerf \cite{cerf68}. Of course $\Theta^n = 0$ for $n = 1 , 2$ and by Perelman $\Theta^3 = 0$. Kervaire-Milnor proved  that $\Theta^4 = 0$. Now Smale's results  imply that $\Lambda^n = \Gamma^n = \Theta^n$ for $n \geq 5$. Hence, finally,  $\Gamma^n = \Theta^n$ for all values of $n$.

{\bf Remarks.} 1) For  $n\neq 4$, every homotopy $n-$sphere can be obtain by gluing two discs along their boundary.
\\
2) For $n=4,$ $\Gamma^4 = 0$ means that every differential structure on  the 4-sphere constructed by gluing two closed 4-balls by a diffeomorphism of their boundary is diffeomorphic to the standard differential structure. But it is unknown whether every differential structure on the 4-sphere can be obtained by such a gluing. On the other hand,  $\Theta^4 = 0$  means that  every homotopy $4-$sphere is h-cobordant to $S^4.$
But we cannot apply Smale's h-cobordism theorem, which requires that the dimension of the h-cobordism is $\geq 6$. 
\\
Moreover it is known today that the {\bf differential} simply connected h-cobordism theorem is not valid in dimension 5. This is an immediate consequence of two results:

1) Simply connected   differential 4-manifold which are homeomorphic are h-cobordant. This follows from Wall in  \cite{wall64}.
\\
2) Many examples of 4-dimension differential manifolds which are homeomorphic but not diffeomorphic are known. The first such examples were constructed by Simon Donaldson.
\\
Note that $n = 4$ is the only dimension for which we do not know if $\Lambda^n$ is a group.

\subsection{Kervaire invariant}

We  define Kervaire invariant $KI_{4k+2} : \pi^S_{4k+2} \rightarrow {\bf Z} / 2$  with the help of  Pontrjagin's  construction. In this setting $KI_{4k+2}$ is defined as an obstruction to framed surgery. The fact that Arf invariant is an obstruction to surgery in dimensions $n \equiv~ 2~mod~4$ was first announced by Milnor in \cite{miln59-1} where, we think, Arf is quoted for the first time in differential topology. Surprisingly, Kervaire does not cite Arf in \cite{kerv60}, but Kervaire-Milnor \cite{kemi63} do. 

Let $x \in \Omega_{4k+2}^{fr}  = \pi_{4k+2}^S$. The element $KI_{4k+2} (x) \in {\bf Z} / 2$ is the obstruction to represent $x$ by a homotopy sphere. Here are some details. For simplicity assume $k \geq 1$.

We represent $x$ by a differential manifold $M^{4k+2}$ differentiably embedded in $S^{(4k+2)+m}$ $(m \geq (4k+2)+2)$ with a framing of its normal bundle. By easy surgery (as explained by Milnor in \cite{miln61}) we can assume that $M^{4k+2}$ is $2k-$connected. The group $H_{2k+1}(M^{4k+2} ; {\bf Z})$  is a free abelian group of even rank $2r$ since the intersection form $I$ is alternate and unimodular.  By Hurewicz theorem $H_{2k+1}(M^{4k+2} ; {\bf Z}) = \pi_{2k+1}(M^{4k+2})$.  Since $M^{4k+2}$ is simply connected, by Whitney's elimination of double points \cite{whit44} every element  of $H_{2k+1}(M^{4k+2} ; {\bf Z})$ can be represented by an embedded sphere $S^{2k+1} \subset M^{4k+2}$. The normal bundle $\nu$ of $S^{2k+1}$ is stably trivial of rank $(2k+1)$. We summarise the prerequisites  in a lemma. See \cite{kosi93} Appendix 1.5 for proofs.

\newpage

\begin{lemma}
(1) Isomorphism classes of oriented vector bundles of rank $r$ over the sphere $S^{d+1}$ are in a natural bijection with $\pi_d(SO_r)$. This is a special case of Feldbau classification theorem.
\\
(2) Stably trivial vector bundles of rank $(2k+1)$ over $S^{2k+1}$ are in bijection with the kernel $\mathcal {K}$ of $\pi_{2k}(SO_{2k+1})  \rightarrow \pi_{2k}(SO)$.
\\
(3) This kernel is cyclic, generated by the  tangent bundle $\tau$ of $S^{2k+1}$.
\\
(4) The order of $\tau$ in $\mathcal{K}$ is equal to 1 if and only if $k = 0 , 1 , 3$ (this corresponds to the dimension of  the spheres which are parallelisable). For other values of $k$ the order of $\tau$ is equal to 2. 
\end{lemma}

There are now two cases.

1) Assume that $2k+1 \neq 1 , 3 , 7$. In this case $\nu \in \mathcal{K}$ is non-necessarily trivial.  Let $o(\nu) \in {\bf Z} / 2$ be equal to $0$ if $\nu$ is trivial  and to $1$ if $\nu =  \tau$. The correspondence  $(S^{2k+1} \subset M^{4k+2}) \mapsto \nu \mapsto o(\nu)$ gives rise to a quadratic form 

$$q : H_{2k+1}(M^{4k+2} ; {\bf Z} / 2) \rightarrow {\bf Z} / 2$$

More precisely we have the equality in $ {\bf Z} / 2$ : $q(y_1 + y_2) = q(y_1) + q(y_2) + I_2(y_1 , y_2)$. Here $I_2$ denotes the intersection form $I$ reduced mod 2. Since the equality takes place in ${\bf Z} / 2$  the quadratic form $q$ is not determined by the bilinear form $I_2$.    

{\bf Definition.} Let $(e_1 , e_2 , \dots , e_r , f_1 , f_2, \dots , f_r)$ be a symplectic basis of $H_{2k+1} (M^{4k+2} ; {\bf Z}) $. The Arf invariant Arf$(q) \in {\bf Z} / 2$ is defined as 

\vspace{.2in}

\centerline{Arf$(q) = \sum_{i = 1}^{r} q(e_i)q(f_i)$}

{\bf Claim.} We have:
\\
i) Arf$(q)$  is the only obstruction to transform $M^{4k+2}$ by framed surgery to a homotopy sphere. 
\\
ii) Arf$(q)$ depends only on the element $x \in \Omega_{4k+2}^{fr} =  \pi_{4k+2}^S$ and not on the manifolds $M^{4k+2}$ chosen to represent it. 

 {\bf By definition} the correspondence $x \in \Omega_{4k+2}^{fr} =  \pi_{4k+2}^S \mapsto M^{4k+2} \mapsto $Arf$(q)$ gives rise to {\bf Kervaire invariant}

$$KI_{4k+2} : \pi_{4k+2}^S \rightarrow {\bf Z} / 2$$

2) Assume that $4k+2 = 2 , 6 , 14$. In this case  the normal bundle $\nu$ is trivial. However there is still an obstruction to  transform $M^{4k+2}$ by a {\bf framed} surgery to a homotopy sphere. It is also given by the Arf invariant of a quadratic form on $H_{2k+1} (M^{4k+2} , {\bf Z} / 2)$. The statements of the claim are also valid in this case. For $4k+2 = 2$ this way of reasoning goes back to Pontrjagin \cite{pont55} who actually proved that $KI_2 : \pi_2^S \rightarrow {\bf Z} / 2$ is an isomorphism.

The ingredients of the proofs are scattered in \cite{kemi63} section 8 and in fact already in \cite{kerv60} for $4k+2 = 10$. A comprehensive presentation is in \cite{levi83} p.81-84 and in Kosinski's book \cite{kosi93} Section X. 

{\bf Note.} Another approach of Kervaire invariant as an obstruction to surgery uses immersion theory. See \cite{wall70}.  There are several other definitions of Kervaire invariant,  more relevant to stable homotopy groups of spheres.  William Browder's paper \cite{brow69} relates Kervaire invariant to Adams spectral sequence. See also several papers by Edgar Brown.  We have merely presented the beginning of a long story.

The {\bf Kervaire Invariant Problem} is to determine for each $n = 4k+2$ whether the homomorphism $KI_n$ is trivial or not. An important step in Kervaire's paper \cite{kerv60} is the proof that $KI_{10}$ is trivial.    In April 2009, Mike Hill, Mike Hopkins and Doug Ravenel  \cite{hihora09} announced that they have solved the Kervaire invariant problem. Soon afterwards, they published a complete proof. Prior to Hill-Hopkins-Ravenel it was known that:

(i) $KI_{8k+2}$ is  trivial for $k \geq 1$ (Brown-Peterson (1966) in \cite{brpe66});
\\
(ii) $KI_{4k+2}$ is trivial if $(4k+2) \neq 2^u - 2$ (Browder (1969) in \cite{brow69});

{\bf However:}
\\
(1) $KI_2$ is non-trivial (Pontrjagin (1955!) in \cite{pont55});
\\
(2) $KI_6$ and $KI_{14}$ are non-trivial (Kervaire-Milnor (1963) in \cite{kemi63});
\\
(3) $KI_{30}$ is non-trivial (Barratt (1969) in \cite{barr69});
\\
(4) $KI_{62}$ is non-trivial (Barratt-Jones-Mahowald (1982) in \cite{bajoma82}).

Right after the publication of \cite{kemi63} it was generally conjectured that $KI_{4k+2}$ is trivial if $4k+2 \neq 2 , 6 , 14 $.  After the publication of Browder's paper  it was conjectured by several algebraic topologists that $KI_{4k+2}$ is non-trivial in the dimensions left open by Browder, i.e. $4k+2 = 2^u - 2$. This was confirmed for $4k+2 = 30 , 62$.  Mahowald derived important consequences of the second conjecture. The answer provided by Hill-Hopkins-Ravenel is close to the first one.

\begin{theorem}
(Hill-Hopkins-Ravenel) The homomorphism $KI_{4k+2}$ is trivial for  all (4k+2)  except 2 , 6 , 14 , 30 , 62 and perhaps 126.
\end{theorem}

\subsection{The groups $P^{n+1}$}

Kervaire and Milnor did not actually define the group $P^{n+1}$ in their paper. But they knew everything about it, since, in fact, they compute it.  Only the ``boundary quotient group" $bP^{n+1}$ is present.  As we shall meet the manifolds involved in the definition of the group $P^{n+1}$ (as Seifert hypersurfaces for knots in higher dimensions) we say a few words about it. Roughly speaking elements of $P^{n+1}$ are represented by parallelisable $(n+1)-$manifolds which have a homotopy sphere as boundary. 

More precisely,  let $F^{n+1}$ be a $(n+1)$-dimension manifold embedded in $\bf {R}^{(n+1) + m}$ with a trivialisation of its normal bundle. We suppose that $m \geq (n+1) + 2$ and that the boundary $bF$ is a homotopy sphere. An equivalence relation ({\bf framed cobordism}) between these manifolds is defined as follows. $F_0$ and $F_1$ are framed cobordant if there exists a framed submanifold $W^{n+2}$ embedded in $\bf {R}^{(n+1) + m} \times \lbrack 0 , 1 \rbrack$ such that:
\\
$W \cap \bf {R}^{(n+1) + m} \times \lbrace i \rbrace = F_i$ for $i = 0 , 1$
\\
$bW = F_0 \cup ( -F_1) \cup U$ where $U$ is an h-cobordism between $bF_0$ and $bF_1$. 

One of Kervaire-Milnor important result  is   the computation of the groups $P^{n+1}$. It goes as follows:

\begin{theorem}

$P^{n+1} = 0$ if $(n+1)$ is odd.
\\
$P^{n+1} = \bf {Z}$ if $(n+1)$ is divisible by 4.
\\
$P^{n+1} = \bf {Z} / 2$ if $(n+1) \equiv 2~~ mod~4$.

\end{theorem}

We offer some comments on this theorem.
\\
The proof that $P^{n+1} = 0$ if $(n+1)$ is odd is highly non trivial. It occupies $20$ pages in the paper.

\newpage

Here is a  description of how the isomorphisms with ${\bf Z}$ and ${\bf Z} / 2$ are obtained. 

For $n+1$ even, let us write $n+1 = 2q$. Assume $q \geq 3$. By easy surgery arguments (essentially \cite{miln61}), Kervaire-Milnor prove that every element of $P^{2q}$ can be represented by a $(q-1)-$connected manifold, say $W^{2q}$. Now $H_q(W ; {\bf Z}) = \pi_q(W)$ is a finitely generated free abelian group. By Whitney's elimination of double points \cite{whit44}, since $q \geq 3$ every element of $H_q(W ; {\bf Z})$ can be represented by an embedded sphere $S^q \hookrightarrow W^{2q}$. The problem is to determine its normal bundle $\nu$. Since $W^{2q}$ is stably parallelisable $\nu$ is stably trivial, but non necessarily trivial. The problem now splits, depending on the parity of $q$. 

If $q$ is even, a stably trivial vector bundle $\eta$ of rank $q$ on the sphere of dimension $q$ is classified by an element $o(\eta) \in {\bf Z}$. For the embedded  $S^q \hookrightarrow W^{2q}$ the element $o(\nu)$ which classifies the normal bundle $\nu$ is equal to the self-intersection of $S^q$ in $W^{2q}$. It follows that  the intersection bilinear form $I_W$ on $H_q(W ; {\bf Z})$ classifies the normal bundle of embedded spheres representing  elements of $H_q(W ; {\bf Z})$.  By  Lemma 3.3.1 and Theorem 13.2.1, the form $I_W$ is even. It turns out that the obstruction to obtain a contractible manifold by framed surgery from $W^{2q}$ is the signature $\sigma (W)$ of $I_W$. Since $I_W$ is even, $\sigma (W)$ is divisible by 8 and the sought homomorphism $P^{2q} \rightarrow {\bf Z}$ is given by the signature divided by 8. By construction this homomorphism is injective. 

If $q$ is odd, we are in the same situation we met in the definition of Kervaire invariant. Now the manifolds we consider have a boundary, but this does not matter. The homomorphism
$P^{2q} \rightarrow {\bf Z} / 2$ is the Arf invariant defined as before. By construction this homomorphism is injective. 

There remains to prove that the homomorphisms $P^{2q} \rightarrow {\bf Z}$ or ${\bf Z} / 2$ are onto. 

We know that  elements of $P^{2q}$ can be represented by framed manifolds $W^{2q}$ which are $(q-1)$-connected. Since the boundary $bW$ is a homotopy sphere, by Poincar\'e duality and Hurewicz theorem $W^{2q}$ has the homotopy type of a wedge (``bouquet") of spheres of dimension $q$. If $q \geq 3$, it results from Smale's h-cobordism theorem that $W^{2q}$ is a q-handlebody. See Section 13 for more on handlebodies.

\begin{definition}
A {\bf q-handlebody} is a manifold $W^{2q}$ of dimension $2q$ which is obtained by attaching handles of index $q$ to a disc of dimension $2q$. 
\end{definition}

{\bf In summary} elements of $P^{2q}$ ($q \geq 3$) can be represented by parallelisable  q-handlebodies such that the intersection form on $H_q (W , {\bf Z})$ is {\bf unimodular}, i.e. of determinant $\pm 1$.

The unimodularity condition insures that the boundary $bW$ is a homology sphere. If $q \neq 2$ the boundary $bW$ of the handlebody is indeed a homotopy sphere.

To prove that the homomorphism $P^{2q} \rightarrow {\bf Z}$ or ${\bf Z} / 2$ is onto the approach is the following. The technical ingredient is the plumbing construction. 
\\
1) If $q$ is even the $E_8$-plumbing produces a q-handlebody which is parallelisable with signature 8 and boundary a homotopy sphere.
\\
2) If $q$ is odd, say $q = 2k+1$,   the argument splits in two sub-cases. 
\\
$2_1$) If $2q = 4k+2$ is equal to 2, 6 or 14 we consider the product $S^q \times S^q$ with a suitable normal framing. If we dig a hole in it, i.e. if we consider  $S^q \times S^q \setminus \mathring {D}^{2q}$ we obtain a framed q-handlebody with Arf invariant 1. 
\\
$2_2$) If $2q = 4k+2$ is not equal to 2 , 6 or 14 we consider the $q$-disc bundle associated to the tangent bundle of $S^q$. Then we plumb two copies of it. We obtain a parallelisable q-handlebody $\widetilde {W}^{2q}$  with boundary a homotopy sphere. Its Arf invariant is equal to 1. Note that in this case we do not have a problem analogous to the Kervaire invariant problem, since our manifolds are allowed to have a boundary. 

The continuation of the story is in next subsection.

\subsection{Kervaire manifold}

In  \cite{kerv60},  Kervaire constructed a closed topological manifold of dimension $10$  which does not have the homotopy type of a differential manifold. We now tell the story.

{\bf Remark} If  $k$ is such that Kervaire invariant $KI_{4k+2}$ is trivial, then the Arf invariant is equal to $0$ for all closed differential stably parallelisable manifolds of dimension $(4k+2)$, essentially by definition of Kervaire invariant.

Consider the manifold $\widetilde {W}^{4k+2}$ we have just constructed. Its boundary is homeomorphic to the sphere $S^{4k+1}$. We glue topologically a $(4k+2)$-disc on its boundary to obtain a manifold $Q^{4k+2}$.  By construction,  it is a closed topological manifold, equipped with a differential structure outside  an open disc. If   the differential structure  of $\widetilde {W}^{4k+2}$  extends over the disc,  the so obtained differentiable manifold  $Q_{diff}^{4k+2}$   would be stably parallelisable  (it could be framed) and by construction would have Arf invariant 1. It is impossible when Kervaire invariant $KI_{4k+2}$ is trivial.
\\
Note that the differential structure on $\widetilde {W}^{4k+2}$ extends if and only if the boundary $\partial \widetilde {W}^{4k+2}$ is diffeomorphic to the standard sphere.

In his paper \cite{kerv60} Kervaire wanted more (in dimension 10)! So  let us forget  the differential structure on $Q^{4k+2}$ outside  a disc and let us consider  $Q^{4k+2}$ only as a topological manifold. This topological manifold is  the {\bf Kervaire manifold} in dimension $4k+2$.  

In \cite{kerv60},  Kervaire proceeds as follows. First he shows, non trivially,   that $KI_{10}$ is trivial.
Then he proves that $Q^{4k+2}$ does not have the homotopy type of a differential manifold in any dimension $(4k+2)$ such that $KI_{4k+2}$ is trivial.   Here are the main steps of the proof which proceeds by contradiction.

Step 1. Construct an invariant $\Phi (V) \in {\bf Z} / 2$ for every triangulable  $(2k)$-connected manifold $V$ of dimension $(4k+2)$ in such a way that $\Phi (V) = KI_{4k+2}(V)$ if $V$ is a framed differential manifold (it is not required that the differential structure is compatible with the triangulation). The definition of $\Phi$ uses cohomology operations. It is a homotopy type invariant.
\\
Step 2. Prove that $\Phi (Q^{4k+2}) = 1$.
\\
Step 3. Prove that if  there exists a differential manifold $Q_{diff}^{4k+2}$ homotopy equivalent to  $Q^{4k+2}$, this manifold $Q_{diff}^{4k+2}$ can be framed. If $2k$ is not congruent to $0 ~mod~8$ the obstruction to framing vanishes by Bott. But Milnor gave a general proof of Step 3 in an appendix to Brown-Peterson paper \cite{brpe66}.
\\
Step 4. Now,  $KI_{4k+2}(Q_{diff}) = 1$ from steps $1$ and $2$. It  gives the wanted  contradiction  each time that  $KI_{4k+2}$ is trivial.

{\bf Consequence of the Hill-Hopkins-Ravenel result.} The Kervaire  manifold $Q^{4k+2}$ does not have the homotopy type of  a differential manifold  for $4k+2 \neq 30 , 62 , 126$.  It does admit one for $4k+2 = 30$ and $62$.  The case $4k+2 = 126$ is still open.

{\bf Historical note.} When Kervaire wrote his paper Smale's (or Stallings') result on the higher dimension Poincar\'e conjecture was not known.  To know that the boundary  $b\widetilde{W}^{10}$ of his plumbing is homeomorphic to the sphere $S^9$ Kervaire relied on a lemma due to Milnor which uses crucially that the plumbing has only two factors. At the same time (even before) Milnor knew about his $E_8$-plumbing. See \cite{miln59-1}. Once Smale's result was known it was easy to construct (by attaching topologically a disc on the boundary of the $E_8$-plumbing) PL manifolds of dimension congruent to 12 mod 16 which do not have the homotopy type of a differential manifold. The contradiction is obtained from Hirzebruch's signature formula instead of Kervaire invariant. The argument is much simpler. Differential topology was moving fast around the year 1960 !

\vskip.5in

\section{ Differential knots  in  codimension  $\geq 3$ }

\subsection { On the isotopy of  knots and links in any codimension }

\begin{definition}  
{\bf A (spherical) n-knot $K^n$ in codimension $q$} is a n-homotopy sphere differentiably embedded in the sphere $S^{n+q}$. Two n-knots $K_0^n$ and $K_1^n$ are equivalent if there exists an orientation preserving diffeomorphism $f : S^{n+q} \rightarrow S^{n+q}$ such that $f(K_0^n) = K_1^n$ (respecting orientations; recall that homotopy spheres are oriented).  
\end{definition}

{\bf Comments.} 1) As Andr\'e Haefliger told us, once homotopy spheres had been identified it seemed natural to investigate how they embed in Euclidean space (or equivalently in a sphere). Observe that the theory of knots in high codimension (a little more than the dimension of the sphere) is isomorphic to  the ``abstract" theory.
\\
2) For people acquainted with classical knot theory, it might seem better to  consider only (images of)  embeddings of the standard sphere $S^n$. But in fact this restriction complicates matters. Images of   embeddings of the standard sphere constitute a subset (in fact a subgroup if $q \geq 3$) of n-knots. 
\\
3) Beware that an orientation preserving diffeomorphism of $S^{n+q}$ is usually not isotopic to the identity. However the following proposition is true.

\begin{proposition}
Let  $Q_0$ and $Q_1$ be two compact subsets of $S^N$. Suppose that there exists an orientation preserving  diffeomorphism $\Phi : S^N \rightarrow S^N$ such that $\Phi(Q_0) = Q_1$. Then there exists a diffeomorphism $\Psi : S^N \rightarrow S^N$ {\bf isotopic to the identity} such that $\Psi (Q_0) = Q_1$.
\end{proposition} 

{\bf Comments.} 1) As a consequence knots $K_0$ and $K_1$ which are equivalent by a diffeomorphism preserving the orientation of the ambient sphere are isotopic. This is also  true for links whatever the codimension may be. Indeed if $\Psi_t$ is the ambient isotopy which connects $\Psi_0$ to $\Psi_1$ then $\Psi_t (K_0)$ is an ambient isotopy which connects $K_0$ to $K_1$. 
\\
2) If $N \leq 3$ it is not necessary to replace $\Phi$ by $\Psi$ since in this case any orientation preserving diffeomorphism of $S^N$ is isotopic to the identity. This is non-trivial and was proved by Smale and also by Munkres   for $N = 1 , 2$ and by Cerf for $N = 3$. 

The proof of the proposition relies on next lemma. For a proof see \cite{munk60} Lemma 1.3 p.523.

\begin{lemma}
Let $f : S^N \rightarrow S^N$ be a diffeomorphism preserving the orientation. Let $B^N \subset S^N$ be a differential ball. Then $f$ is isotopic to a diffeomorphism $f^*$ which is the identity on $B^N$. 
\end{lemma}

To prove the lemma one uses the following result, due   to Jean Cerf \cite{cerf61}.

\begin{theorem}
Let $M^N$ be a connected and oriented differential N-dimension manifold. It does not need  to be compact and may have a boundary. Let $f_0$ and $f_1$ be two orientation preserving embeddings of the ball $B^N$ in $Int(M)$. Then there exists a diffeomorphism $\varphi : M^N \rightarrow M^N$, isotopic to the identity, such that $\varphi \circ f_0 = f_1$.
\end{theorem}

This result is also the key tool to prove that the connected sum of  oriented differential manifolds is well defined up to orientation preserving diffeomorphism. For a proof see for instance \cite{hirs76} near p.117.
  
{\bf Proof of the proposition.}    

If $Q_1 = S^N$ there nothing to prove. Hence we may suppose that $Q_1 \neq S^N$ and let $B^N$ be a differentiable ball such that $Q_1 \subset B^N \subset S^N$. We apply the lemma with $f = \Phi^{-1}$ and $B^N$ as above. Let $f^*$ be a diffeomorphism as in the lemma. 
\\
If $x \in Q_0$ then $\Phi (x) \in Q_1$ and hence $f^*(\Phi(x)) = \Phi (x)$. Hence $f^* \circ \Phi$ sends also $Q_0$ to $Q_1$. On the other hand let $f_t$ be an isotopy which connects $f_0 = f = \Phi^{-1}$ to $f_1 = f^*$. Let $h_t = f_t \circ \Phi$. We have $h_0 = f_0 \circ \Phi = \Phi^{-1} \circ \Phi = id$ and $h_1 = f_1 \circ \Phi = f^* \circ \Phi$. Hence $\Psi = f^*\circ  \Phi$ is the diffeomorphism we are looking for.

{\bf End of proof of the proposition.}

\subsection { Embeddings and isotopies in the stable and metastable ranges}

The paper which initiated the theory of differential knots in codimension $\geq 3$ is Andr\'e Haefliger's paper \cite{haef62-2}. Here is a short history of the events.

In 1936 and 1944 Hassler Whitney published (see \cite{whit36} and \cite{whit44}) two foundational papers about differential manifolds. His results  can be summarised as follows.

\begin{theorem}(Embeddings and isotopies  in the stable range)
Let $M^n$ be a  differential n-dimensional  manifold. Then:
\\
1) $M^n$ can be differentiably embedded in ${\bf R}^{2n}$;
\\
2) any two differentiable embeddings of $M^n$ in ${\bf R}^{2n+2}$ are isotopic;
\\
3) if $n \geq 2$ and if $M^n$ is connected, any two differentiable embeddings of $M^n$ in ${\bf R}^{2n+1}$ are isotopic.
\end{theorem}

Around 1960, Haefliger presented  a vast theory of embeddings of differential manifolds in differential manifolds,  announced in \cite{haef61}. As far as homotopy spheres are concerned his results are as follows.

\begin{theorem}(Embeddings and isotopies  in the metastable range) Let $\Sigma^n$ be a  homotopy n-sphere. Then:
\\
1) $\Sigma^n$ can be differentiably embedded in ${\bf R}^N$ if $2N \geq 3(n+1)$;
\\
2) Any two differentiable embeddings of $\Sigma^n$ in ${\bf R}^N$ are isotopic if $2N > 3(n+1)$. 
\end{theorem}

\subsection{Below the metastable range}

At about the same time (1960) Christopher Zeeman announced in  \cite{zeem61} the following result, proved in \cite{zeem63} (from Smale PL theory it results that PL homotopy spheres of dimension $\geq 5$ are PL isomorphic to $S^n$).

\begin{theorem}
Any two piece-wise linear (PL) embeddings of the sphere $S^n$ in ${\bf R}^N$ are PL isotopic if $N \geq n+3$. 
\end{theorem}

Naturally the question arose: what happens below the metastable range for differentiable  embeddings of homotopy spheres? By general position the complement is simply connected if $q \geq 3$.  From Zeeman's result no topological invariant  can be obtained from the knot complement, since it is homeomorphic to $S^{q-1} \times B^{n+1}$ . For small values of $n$ the first critical values are $n = 3$ and $N = 6$. It was a surprise when Haefliger announced that there exist infinitely many isotopy classes of differentiable embeddings of $S^3$ in $S^6$. The full result is the following~\cite{haef62-2} which indicates that there is a important difference between PL and DIFF.

\begin{theorem}
Let $k$ be an integer $\geq 1$. Then there exist infinitely many  isotopy classes of differentiable embeddings of $S^{4k-1}$ in $S^{6k}$. 
\end{theorem}

To make later a parallel  with knots in codimension 2, we give a brief idea of Haefliger's approach. The following key definition is the starting point. 

\begin{definition}
Let $M_i^n$ for $i = 0, 1$ be two oriented manifolds differentiably embedded in $S^{n+q}$. They are said to be {\bf h-cobordant} it there exists an oriented manifold $W^{n+1}$ differentiably and properly embedded in $S^{n+q} \times \lbrack 0 , 1 \rbrack$ such that:
\\
1) $W^{n+1} \cap S^{n+q} \times \lbrace i \rbrace = M_i^n$ for $i = 0 ,1$;
\\
2) $W^{n+1}$ is an oriented h-cobordism between $M_0^n$ and $M_1^n$. In particular $bW = M_1 - M_0$.
\end{definition}

One of the many  Smale's results of \cite{smal62} reads as follows.

\begin{theorem} (h-cobordism implies isotopy)
Suppose that the differential manifolds $M_i^n$ in $S^{n+q}$  are h-cobordant. Suppose that the h-cobordism $W^{n+1}$ is simply connected and  that $n \geq 5$ and  $q \geq 3$. Then $M_0^n$ and $M_1^n$ are isotopic.
\end{theorem}

There are two steps in Smale's proof. In the first step, simple connectivity and $n \geq 5$ are used to deduce (h-cobordism theorem) that $W^{n+1}$ is a product. In the second step, the hypothesis $q \geq 3$ is used to straighten the product $W^{n+1}$ inside $S^{n+q} \times \lbrack 0 , 1 \rbrack$. The hypothesis $q \geq 3$ is crucial here as we shall see later in the section  dedicated to ``knot cobordism in codimension 2". 

 For $q\geq 3$, Haefliger proves that isotopy classes of embedded n-dimensional  homotopy spheres in $S^{n+q}$ form an abelian group $\Theta^{n+q,n}$ under the connected sum operation.   This is noteworthy, since knots in codimension 2 do not form a group but only a monoid. The inverse is lacking. Haefliger denotes by $\Sigma^{n+q,q}$ the subgroup which consists in embeddings of the standard sphere $S^n$. The main result of \cite{haef62-2} is that $\Sigma^{6k,4k-1}$ is isomorphic to the integers ${\bf Z}$. The isomorphism is constructed as follows.

Suppose that $K$ is an embedded $S^{4k-1}$ in $S^{6k}$. Consider $S^{6k}$ as the boundary of the disc $B^{6k+1}$. Using results of Kervaire, Haefliger proves that $K$ bounds in $B^{6k+1}$ a framed submanifold $V^{4k}$. Very roughly speaking, there is an obstruction to transform $V^{4k}$ by surgeries in order to obtain an embedded disc bounding $K$ in $B^{6k+1}$. This obstruction is an integer which provides the isomorphism between $\Sigma^{6k,4k-1}$ and ${\bf Z}$. The generator is Haefliger's famous construction based on a high dimensional version of the Borromean rings. 

Once \cite{haef62-2} was published the question arose to determine the groups $\Theta^{n+q,n}$ and $\Sigma^{n+q,n}$ in general. The answer was provided by Jerome Levine in \cite{levi65-2}. To understand Levine's results it is necessary to go back a littler bit in time. Kervaire and Milnor realised that ``their" short exact sequence is in fact a consequence of  a large diagram of groups which consists in four braided long exact sequences. This was very likely the intended content of Part II of ``Groups of Homotopy Spheres", which never appeared. Levine proved that  the groups $\Theta^{n+q,n}$  are also parts of a diagram  of four braided long exact sequences. His diagram is the unstable version of Kervaire-Milnor's braided diagram. See also \cite{haef64}. Levine says also how to recapture $\Sigma^{n+q,n}$ as a subgroup of $\Theta^{n+q,n}$.

A consequence of Levine's result is that the abelian groups $\Theta^{n+q,n}$ are finitely generated and of rank $\leq 1$ over ${\bf Q }$ if $q \geq 3$. Thus we can say that below the metastable range there exist  differentiably knotted homotopy spheres but that, if $q \geq 3$, there are rather few of them. This contrasts sharply with the situation $q = 2$ as we shall see. 

{\bf Comment.}  Andr\'e Haefliger has shown us letters that Jerome Levine wrote to him in 1962 and 1963. They indicate that, since 1962,  Levine was  very actively  working on spherical knots in codimensions $q \geq 2$. Here is an example  he described in a letter dated May 7, 1963 and was apparently never published. We present this example since it ties together several concepts which are put forward in our paper.

Let $2k \geq 4$. Let $V^{4k}$ be a parallelisable ${2k}-$handlebody with boundary $bV$ diffeomorphic to $S^{4k-1}$. Suppose that the signature of $V$ is non-zero. Remark  that this signature is a multiple of  $8\times order( bP^{4k})$ and hence quite large. 
\\
Since $V^{4k}$  is a parallelisable handlebody   it is not hard to embed it in $S^{4k+1}$. In fact it embeds in many ways, the key notion here being the Seifert matrix (see below Section 7). Let $\phi : V^{4k}\rightarrow S^{4k+1}$ be such an embedding. The spherical knot $K$ that Levine considers is $\phi(bV) \subset S^{4k+1}$. Let $j_N : S^{4k+1} \hookrightarrow S^N$ be the standard embedding for any $N \geq 4k+1$. 

\begin{proposition}
$j_N\phi(bV) \subset S^N$ is differentiably knotted for any $N$ such that $4k+1 \leq N \leq 6k-1$.
\end{proposition}

{\bf Proof of the proposition.} It is enough to prove the statement for $N = 6k-1$. 
\\
Consider $S^{6k-1}$ as the equator of $S^{6k}$. It bounds two balls $B_+^{6k}$ and $B_{-}^{6k}$. Keeping $K = bV$ fixed in $S^{6k-1}$ we can push the interior of $V^{4k}$  in the interior of $B_+^{6k}$ (we neglect $j_N\phi$).
\\
From now on we argue by contradiction. Suppose that $K$ is differentiably unknotted in $S^{6k-1}$. Hence it bounds a differentiable disc $B^{4k}$ in  $S^{6k-1}$. We push the interior of  $B^{4k}$ in the interior of $B_{-}^{6k}$. The union $V^{4k} \cup B^{4k}$ is a differentiable manifold $\hat {V}$ of dimension $4k$ differentiably embedded in $S^{6k}$. It is almost parallelisable with non-zero signature. This is impossible. Here is why. 
\\
Let $p_k \in H^{4k}(\hat {V} ; {\bf Z})$ be the k-th Pontrjagin class of $\hat {V}$. Since $\hat {V}$ is almost parallelisable, by Hirzebruch's signature formula the signature of $\hat {V}$ is detected by $p_k$. Hence $p_k$ is non-zero. 
\\
By Whitney duality $\bar {p}_k= -p_k$ where $\bar {p}_k$ is the k-th Pontrjagin class of the normal bundle $\nu$ of $\hat {V}$ in $S^{6k}$. Next lemma provides the contradiction.

\begin{lemma}
Let $M^{4k} \subset S^{6k}$ be a closed, connected and oriented differential submanifold. Then $\bar {p}_k \in H^{4k}(M^{4k} ; {\bf Z})$ vanishes. 
\end{lemma}

{\bf Proof of the lemma.} The normal bundle $\nu$ of $M^{4k}$ in $S^{6k}$ is of rank 2k. Let $\bar {e} \in H^{2k}(M^{4k} ; {\bf Z})$ be the Euler class of the normal bundle $\nu$.     From Theorem 31 of Milnor's Lecture Notes  \cite{miln74} applied to the normal bundle $\nu$ we have $\bar {p}_k = \bar {e}^2$. But $\bar {e}$ vanishes since it is the normal Euler class of $M^{4k}$ embedded in $S^{6k}$, by Thom's definition of the Euler class. See \cite{miln74} Theorem 14.

\newpage

\section{The fundamental group of the knot complement}

From now on, a n-knot will mean a spherical n-knot in codimension $q = 2$. A 1-knot is often called a {\bf classical knot}. A {\bf standard n-knot} is the image of an embedding of the standard n-sphere.  

Question for mathematicians (topologists): how can we represent (or construct) n-knots, for $n \geq 2$? There is a difficulty there and very likely it prevented the advancement of knots in higher dimensions before the sixties. At that time two constructions  became available: surgery and embeddings (in $S^{n+1} $) of parallelisable $(n+1)-$manifolds  which have a homotopy sphere as boundary (Seifert hypersurfaces).

The paper \cite{kewe77} contains a  survey of the theory  of  n-knots   before 1960.

\subsection{Homotopy n-spheres embedded in $S^{n+2}$}

\begin{proposition} 
Let $L^n$ be a closed oriented  manifold without boundary embedded in $S^{n+2}$. Then its normal bundle is trivial.
\end{proposition}

{\bf Proof of proposition  5.1.1} Let us observe first that the normal bundle is orientable. Since it has rank 2,  it is enough to prove that it has a nowhere vanishing section. The only obstruction to construct  such a section is the Euler class $e \in H^2(L^n ; {\bf Z})$. Since $H^2(S^{n+2} ; {\bf Z}) = 0$ Thom's formula for the Euler class implies that $e = 0$. Observe that the same proof works if $S^{n+2}$ is replaced by an orientable manifold $W^{n+2}$ such that $H^2(W^{n+2} ; {\bf Z}) = 0$.

{\bf Construction of knots by surgery.}

Michel Kervaire uses this result in the following way to construct several n-knots by surgery. Suppose that  the n-knot $K$ is the standard sphere  $S^n$  embedded in $S^{n+2}$. Thanks to the proposition above we can perform surgery on $K$. Moreover, if $n \geq 2$ there is a unique trivialisation of the normal bundle. Therefore, there is a well-defined manifold $M^{n+2}$ obtained by doing surgery on the  knot  $K \subset S^{n+2}$.  Michel Kervaire's idea is to construct  first this manifold  $M^{n+2}$ and then perform the surgery backwards to get the knot.  Observe that this procedure works only for knots represented by an embedding of the standard sphere.  

Let then $K \times B^2 \subset S^{n+2}$. Let $E$ be the {\bf knot exterior} defined by $S^{n+2} \setminus (K \times \mathring{B}^2)$.  The advantage of the knot exterior over the knot complement $S^{n+2} \setminus K$ is that it is a compact manifold with boundary, having the same homotopy type of the complement. 

{\bf Lovely observation.} The exterior $E$ is the common part of the surgery.

Denote by $M^{n+2}$ the manifold obtained by doing surgery on $K$ in $S^{n+2}$. Remark that, by construction, $M^{n+2}$ contains the scar (an embedded circle) $\gamma \subset M^{n+2}$ such that if we perform a surgery on $\gamma$ we get the manifold $S^{n+2}$ with its knot $K$ inside as the scar. Next proposition is an easy consequence of Alexander and Poincar\'e dualities.

\begin{proposition}
Let $n \geq 2$ and let $M^{n+2}$ be the manifold obtained by surgery on a standard n-knot in $S^{n+2}$. Then $M^{n+2}$ has the homology of $S^1 \times S^{n+1}$. Explicitly:
\\
$H_i(M^{n+2} ; {\bf Z}) = {\bf Z}$ if $i = 0 , 1 , n+1 , n+2$ and $H_i(M^{n+2}  ; {\bf Z}) = 0$ otherwise.  
\end{proposition}

We shall see later several applications of Kervaire's construction of a manifold like $M^{n+2}$. More  informations on the knot exterior  are given in Subsection 6.1.

{\bf Seifert hypersurfaces.}

Next theorem was proved by Michel Kervaire in \cite{kerv63}. But it was certainly known to most differential topologists at the time of the meeting in the honor of Marston Morse.

\begin{theorem}
(Kervaire) Let $\Sigma^n$ be a homotopy n-sphere differentiably embedded in $S^{n+2}$. Then there exists an orientable (n+1)-differentiable manifold $F^{n+1} \subset S^{n+2}$ such that $bF^{n+1} = \Sigma^n$.
\end{theorem}

\begin{definition}
A submanifold as $F^{n+1}$ is called a {\bf Seifert hypersurface} of the knot $\Sigma^n$.
\end{definition}

\begin{corollary}
A homotopy n-sphere can be embedded in $S^{n+2}$ if and only if it bounds a parallelisable manifold. In other words, the elements of $\Theta^n$ which can be differentiably embedded in $S^{n+2}$ are exactly the elements of $bP^{n+1}$.
\end{corollary}

The existence of Seifert hypersurfaces  is proved in a  more general context  in Section 11. The corollary is discuss in Proposition 14.1.1.

{\bf Remark on the groups $\Theta^{n+q,n}$ for $q \geq 3$.} 
\\
We have just stated  that every homotopy n-sphere embedded in $S^{n+2}$ has a trivial normal bundle and bounds a framed  $(n+1)-$dimension submanifold in $S^{n+2}$. Both these statements do not hold for homotopy spheres embedded in codimension $\geq 3$. More precisely Levine in \cite{levi65-2} constructs two homomorphisms:
\\
1) a homomorphism $\Theta^{n+q,n} \rightarrow \pi_{n-1} (SO_q)$ which expresses the obstruction of the normal bundle of $\Sigma^n \hookrightarrow S^{n+q}$ to be trivial. 
\\
2) a homomorphism $\Theta^{n+q,n} \rightarrow \pi_n (G_q , SO_q)$ which expresses the obstruction of the embedded $\Sigma^n \hookrightarrow S^{n+q}$ to bound a framed submanifold in $S^{n+q}$. 
See Levine's paper  \cite{levi65-2} for a definition of the H-space $G_q.$
\\
These two homomorphisms take place in two long exact sequences which are the essential constituent of Levine's diagram of four braided long exact sequences  (the unstable version of Kervaire-Milnor's).

\subsection{Necessary conditions for a group to be the fundamental group of a knot complement}  

Michel Kervaire's paper which got knots in higher dimensions really started was his determination of the fundamental group.

\begin{theorem}
(Kervaire \cite{kerv63}) Let $K$ be a n-knot ($n \geq 1$). Let $\Gamma$ be the fundamental group of the complement $S^{n+2} \setminus K$. Then $\Gamma$ satisfies the following conditions.
\\
1) $H_1(\Gamma) = {\bf Z}$
\\
2) $H_2(\Gamma) = 0$
\\
3) $\Gamma$ is finitely presented
\\
4) $\Gamma$ is the normal closure of a single element
\end{theorem}

{\bf Comments.} 1) The reason why we state Michel Kervaire's conditions in this order will be clear in what follows.
\\
2) Usually the $\bf weight$ of a group is defined to be the minimal number of elements of the group such that the normal closure of this set of  elements is the whole group. Condition 4) can hence be stated: $\Gamma$ is of weight 1.  Another way to state this condition is: there exists an element $z \in \Gamma$ such that every element of $\Gamma$ is a product of conjugates of $z$ and $z^{-1}$.
\\
3) $H_i(\Gamma)$ denotes the i-th homology group of $\Gamma$ with integer coefficients {\bf Z} and  with the trivial action of $\Gamma$ on {\bf Z}.
\\
4) We recall that $H_1(\Gamma) = \Gamma / \Gamma'$ is the abelianisation of $\Gamma$.
\\
5) To understand $H_2(\Gamma)$ the key result is Hopf theorem on $H_2(\Gamma)$.

\begin{theorem}
(Hopf theorem) Let $X$ be a connected ``good" topological space  (for instance a C.W. complex) with fundamental group $\Gamma$. Let $h_2 : \pi_2(X) \rightarrow H_2(X ; {\bf Z})$ be the Hurewicz homomorphism. Then the quotient $H_2(X ; {\bf Z}) / Im(h_2)$ depends only on $\Gamma$ and is by definition  $H_2(\Gamma)$.
\end{theorem}

\begin{corollary}
The following statements are equivalent.
\\
1) $H_2(\Gamma) = 0$.
\\
2) For all spaces $X$ as in the theorem, the  Hurewicz homomorphism $h_2 : \pi_2(X) \rightarrow H_2(X ; {\bf Z})$ is onto.
\\
3) There exists a space $X$ such that $H_2(X ; {\bf Z}) = 0$.
\end{corollary}

{\bf Proof of Theorem 5.2.1}

Conditions 1) and 2) are an immediate consequence of Alexander duality. Indeed let $Z \subset S^{n+2}$ be a compact subset of $S^{n+2}$ which has the \u{C}ech cohomology of $S^n$. Then $S^{n+2} \setminus Z$ has the singular homology of $S^1$.  In particular Conditions 1) and 2) are also true for wild embeddings of $S^n$ in $S^{n+2}$.

For Conditions 3) and 4) to be satisfied requires more on the embedding. To simplify matters let us assume that $K$ is an homotopy  n-sphere topologically embedded in $S^{n+2}$ with a neighbourhood homeomorphic to $K \times B^2$. 

Let $E$ be the exterior $S^{n+2} \setminus (K \times \mathring {B}^2)$. It is clear that $E$ has the same homotopy type as the knot complement $S^{n+2} \setminus K$. Since $E$ is a compact topological manifold, its fundamental group is finitely presented. This proves that Condition 3) is satisfied. 

We can reconstruct $S^{n+2}$ from the exterior by attaching to $E$ a 2-cell $e^2 = \lbrace x \rbrace \times B^2$ for some $x \in K$ and then a (n+2)-cell $e^{n+2}$. Let $z \in \pi_1 (E)$ be represented by the loop $\lbrace x \rbrace \times \partial B^2$. Since $n \geq 1$ and since $S^{n+2}$ is simply connected, $E \cup e^2$ is simply connected. Hence the normal subgroup of $\pi_1 (E)$ generated by $z$ is equal to the whole group. This proves Condition 4).  The loop $\lbrace x \rbrace \times \partial B^2$ and  its class $z \in \pi_1 (E)$  are often called a {\bf meridian}. 

{\bf End of proof of Theorem 5.2.1}

{\bf A comment about meridians and killers.} 
\\
The following terminology seems to get established. A {\bf killer} $\gamma$ is an element of a group $\Gamma$ such that its normal closure is equal to the whole group. Define two killers $\gamma$ and $\gamma^*$ to be equivalent if there exists an automorphism $\Psi: \Gamma \rightarrow \Gamma$ such that $\Psi (\gamma) = \gamma^*$. It is quite natural to ask whether there exist knot groups with killers non equivalent to a meridian. The answer is yes. In \cite{siwhwi10}  Silver, Whitten and Williams prove that there exist  classical knot groups with infinitely many non equivalent killers. For instance, this is the case for rational knots and for torus knots. The authors put forward the exciting (or risky?) conjecture that every classical knot group (of course non abelian) has infinitely many non equivalent killers. 
\\
Jonathan Hillman in \cite{hill89} calls a killer a {\bf weight element}.  He introduces the following convenient terminology. A {\bf weight orbit} is the equivalence class of a weight element through group automorphisms; a {\bf weight class} is the equivalence class through conjugation.  He proves that knot groups which are metabelian have only one   weight class. Classical knot groups are not metabelian by Neuwirth \cite{neuw65}.

\subsection{Sufficiency of the conditions if $n \geq 3$}

We begin this subsection by presenting Kervaire's construction of knots by surgery. 

\begin{definition}
Let $n \geq 2$. A {\bf knot manifold $M^{n+2}$} is a closed, connected and oriented differential manifold of dimension $(n+2)$ such that:
\\
1) $H_*(M^{n+2} ; {\bf Z}) = H_*(S^1 \times S^{n+1} ; {\bf Z})$
\\
2) $\pi_1(M)$ has weight 1.
\end{definition}

Surgery on a spherical n-knot in $S^{n+2}$ produces a manifold satisfying Condition 1) (see Proposition 5.1.2) together with a killer (see the proof of Theorem 5.2.1). 

Conversely let $M^{n+2}$ be a knot manifold and $z \in \pi_1(M)$ a weight element. Following Kervaire, we show now how to obtain a spherical n-knot from this data by surgery if $n \geq 3$.

Let $z \in \pi_1 (M)$ be the chosen killer.  We represent this element by an embedded circle $\zeta : S^1 \hookrightarrow M$. Note that the isotopy class of this embedding is well defined, once $z$ is chosen (only its weight class matters).

{\bf Claim.} The normal bundle of  $\zeta : S^1 \hookrightarrow M$ is trivial and there are two homotopy classes of trivialisations of this bundle.

The first assertion is true since $M$ is orientable and the second one is true since homotopy classes of trivialisations are in bijection with $\pi_1 (SO_{n+1})$. 

We choose a trivialisation $S^1 \times B^{n+1} \hookrightarrow M$ and perform a surgery on $S^1 \times B^{n+1}$. Let $M'$ be the manifold obtained by the surgery.

\begin{proposition}
The manifold $M'$ is a homotopy sphere.
\end{proposition}

{\bf Proof of Proposition 5.3.1}

Let $V = M \setminus (S^1 \times \mathring {B}^{n+1})$ be the common part of $M$ and $M'$. 

The manifold $M$ is obtained from $V$ by attaching a (n+1)-cell and a (n+2)-cell. Hence $V$ and $M$ have the same n-skeleton. Since $n \geq 3$ we have $\pi_1 (V) = \pi_1 (M) = \Gamma$ and $H_i (V  ; {\bf Z}) = H_i (M  ; {\bf Z}) = 0$ for $2 \leq i \leq n-1$. 
\\
On the other hand $M'$ is obtained from $V$ by attaching a 2-cell and a (n+2)-cell. By construction the 2-cell ``kills" the element $z$. Since $z$ is of infinite order in $H_1 (V  ; {\bf Z})$ the 2-cell does not create torsion in $H_2$. Since $z$ is a killer, the fundamental group $\pi_1 (M')$ is trivial. Moreover $H_i (M'  ; {\bf Z}) = 0$ for $1 \leq i \leq n-1$. Hence, by Poincar\'e duality, $M'$ is a homology sphere and, since its  fundamental group is trivial, $M'$ is a homotopy sphere. 

{\bf End of proof of Proposition 5.3.1}

The homotopy sphere $M'$ has some differential structure but we cannot guarantee that this structure is the standard one. Since its dimension $(n+2)$ is $\geq 5$ we can  perform a connected sum with its inverse in the group $\Theta^{n+2}$ to obtain the standard differential structure. This connected sum can take place on a small (n+2)-ball. We  write $M'_{mod}$ for this modified differential structure such that $M'_{mod} = S^{n+2}$.

Where is the knot $K$ in $M' _{mod}= S^{n+2}$ ? Well $M'_{mod}$ is constructed by surgery. The common part is the knot exterior, and the knot is the scar of the surgery in $M'_{mod}$. Let us remark for future reference that the knot exterior has the same n-skeleton as $M$. 

{\bf Remark.} One defect of the construction of knots by surgery is that it produces only standard knots. This is the price we pay to work in the differential category. It is possible to mimic this construction in the PL or in the TOP category, once adequate hypothesis are formulated to ensure the existence of trivial normal bundles. See Kervaire-Vasquez \cite{keva67} for the PL case and Hillman \cite{hill89} for the TOP case. 

{\bf Gluck reconstruction.} We have noticed in the claim above that the normal bundle of the embedding $\zeta : S^1 \rightarrow M^{n+2}$ has two trivializations. Hence we can construct a priori two n-knots by surgery. By definition, these two knots differ by  Gluck reconstruction. Clearly they have the same exterior. We shall see in Section 6 that here are at most two knots with the same exterior. 

{\bf The case $n = 2$.} Let $M^4$ be a knot manifold and $z \in \pi_1 (M)$ be a weight element. We perform a surgery on $z$. We get a homotopy  4-sphere $\Sigma^4$ together with a 2-knot in it. Today it is not known if we can modify the differential structure on $\Sigma^4$ to get a 2-knot in $S^4$. Fifty years ago this was presumably seen as a shortcoming. Today it is quite the contrary. Maybe this is the way to construct counterexamples to the Differential Poincar\'e Conjecture in dimension 4! One difficulty (among several) is to construct 4-dimension knot manifolds. Examples are due to Cappell and Shaneson. See \cite{cash76-2}.

Kervaire's realization theorem is the following.

\begin{theorem}(Kervaire)
Let $\Gamma$ be a group satisfying Conditions 1) to 4) of theorem 5.2.1. Suppose that $n \geq 3$. Then there exists an embedding of the standard n-sphere $K \hookrightarrow S^{n+2}$ such that the fundamental group of the exterior is isomorphic to $\Gamma$. 
\end{theorem}

The first step in the proof of the realization theorem is next proposition.

\begin{proposition}
Let $\Gamma$ be a group satisfying Conditions 1) , 2) and 3) of the theorem. Let $n \geq 3$. Then there exists a differential manifold $M^{n+2}$, orientable, compact, connected, without boundary such that:
\\
i) $\pi_1 (M) = \Gamma$
\\
ii) $H_i (M  ; {\bf Z}) = {\bf Z}$ if $i = 0 , 1 , n+1 , n+2$ and $= 0$ otherwise. 
In other words, additively $H_* (M  ; {\bf Z}) = H_* (S^1 \times S^{n+1}  ; {\bf Z})$.
\end{proposition}

Clearly, the theorem follows from a  surgery on the  knot manifold provided by the proposition applied to a group satisfying Conditions 1) to 4).

Before  the proof of Proposition 5.3.2 we state and prove another proposition. 

\begin{proposition}
Let $\Pi$ be group satisfying Conditions  2) and 3). Then there exists a compact and connected C.W. complex $Y$ of dimension 3 such that:
\\
i) $\pi_1(Y) = \Pi$.
\\
ii) $H_i (Y  ; {\bf Z}) = 0$ for $i \geq 2$.  
\end{proposition}

{\bf Proof of Proposition 5.3.3}

We present the proof with some details in order to show how Condition 2) and Hopf theorem are used.  (Conditions 1) and 4) are not necessary).
\\
Choose arbitrarily a finite presentation of the group $\Pi$. Then construct a finite C.W. complex $Y'$ of dimension 2 with $\pi_1 (Y') = \Pi$, following the instructions contained in the presentation: one 0-cell, the 1-cells correspond to  the generators and  the 2-cells correspond to the  relations. We then consider the chain complex of $Y'$, with integer coefficients.

$$C_2(Y') \stackrel {\partial_2}{\longrightarrow} C_1 (Y') \stackrel {\partial_1}{\longrightarrow} C_0 (Y')$$

We have $H_2 (Y'  ; {\bf Z}) = Ker (\partial_2)$ since $C_3 (Y') = 0$. Hence $H_2 (Y'  ; {\bf Z})$  is a finitely generated free abelian group, of rank say $s$. Let $\lbrace g_1 , \dots , g_s \rbrace$ be a basis of $H_2 (Y'  ; {\bf Z})$. By construction $\pi_1 (Y') = \Pi$ and since $H_2 (\Pi) = 0$, by Hopf theorem, the  Hurewicz homomorphism $h_2 : \pi_2 (Y') \rightarrow H_2(Y'  ; {\bf Z})$ is onto.  Represent then each $g_i$ by a map $\varphi_i : S^2 \rightarrow Y'$.  We attach $s$ 3-cells $e_1^3 , \dots , e_s^3$ on $Y'$, where  $e_i^3$ is  attached by the map $\varphi_i$ for $i = 1 , \dots , s$. We denote by $Y$ the 3-dimension C.W. complex thus obtained. We claim that this is the complex we are looking for. 
\\
i) Since $Y$ is obtained by attaching cells of dimension 3 to $Y'$ we have $\pi_1 (Y) = \pi_1 (Y') = \Pi$.
\\
ii) Consider the following portion of the chain complex of $Y$:

$$C_3(Y) \stackrel {\partial_3}{\longrightarrow} C_2 (Y) \stackrel {\partial_2}{\longrightarrow} C_1 (Y)$$

By construction we have $Ker {\partial_2} = Im {\partial_3}$ and hence $H_2(Y  ; {\bf Z}) = 0$. Now $H_3 (Y  ; {\bf Z}) = Ker (\partial_3)$ since there are no 4-cells in $Y$. But $Ker (\partial_3) = 0$ since the basis elements of $C_3 (Y) $ are sent by $\partial_3$ to linearly independent elements.

{\bf End of proof of Proposition 5.3.3}

{\bf Proof of Proposition 5.3.2}

Let $\Gamma$ be a group satisfying Kervaire conditions 1) to 4). Let $Y$ be a C.W. complex supplied by Proposition 5.3.3 such that $\pi_1 (Y)= \Gamma$. The idea to construct the manifold $M^{n+2}$ is to first construct a manifold $N^{n+3}$ which has the complex $Y$ for a spine. The manifold $M^{n+2}$ is then the boundary $\partial N^{n+3}$. The fact that $\partial N^{n+3}$ has the required properties is easy. Since the codimension of $Y$ in $N^{n+3}$ is $n\geq 3$, we have $\pi_1 (M) = \pi_1 (N) = \pi_1 (Y)$. Poincar\'e duality implies that the homology of $M$ is as desired.

There are at least two (closely related indeed) ways to construct $N^{n+3}$. 

The first method is to attach handles to a (n+3)-dimension disc. The handles are of index 1 , 2 and 3 and attached step by step in a way similar to the way the complex $Y$ is constructed. This is the path followed by Michel Kervaire in \cite{kerv63} p. 110-111. In fact Michel Kervaire uses the surgery technique (certainly it is not a surprise). This is equivalent to construct a thickening of $Y$. Michel Kervaire's argument works for $n \geq 3$.

A more expeditious method is to construct a finite 3-dimension simplicial complex $Y^*$ homotopy equivalent to $Y$. We then embed $Y^*$ piece-wise linearly in ${\bf R}^{n+3}$. The manifold $N^{n+3}$ is a regular neighbourhood of $Y^* \hookrightarrow {\bf R}^{n+3}$. Since the simplicial complex $Y^*$ is of dimension 3, by general position it can be embedded in ${\bf R}^{n+3}$ for $n \geq 4$. But, since $H_3 (Y^*) = 0$ it can also be embedded in ${\bf R}^6$. Hence both methods produce knots in $S^{n+2}$ for $n \geq 3$.  
\\
{\bf End of proof of Proposition 5.3.2. }

{\bf Remark.} The proof given above for Theorem 5.3.1 shows that any group satisfying Kervaire conditions can be realised as the fundamental group of the exterior of an embedding of the standard n-sphere. We have seen above that a homotopy n-sphere can be embedded in $S^{n+2}$ if and only if it bounds a parallelisable manifold. In \cite{kerv63} Kervaire proves, following an argument of Milnor, that  any group satisfying Kervaire's conditions can be realised as the fundamental group of the exterior  of an embedding of any given element of $bP^{n+1}$.

\subsection{Kervaire conjecture}

{\bf Question:} Can a free product $G \ast {\bf Z}$ be a group satisfying Kervaire conditions?

More precisely we suppose that $G$ is finitely presented and satisfies $H_1(G) = 0 = H_2(G)$. Then clearly $G \ast {\bf Z}$ satisfies Kervaire conditions 1) , 2) and 3). Is it possible that this free product also satisfies condition 4)? In \cite{kerv63} Kervaire studied the free product $I_{120} \ast {\bf Z}$ where $I_{120}$ is the binary icosahedral group of order 120. He proved that $I_{120} \ast {\bf Z}$ is not of weight 1. See p.117. Incidentally this proves that Condition 4) is not a consequence of Conditions 1) , 2) and 3). Kervaire asked whether it is possible that a non-trivial free product $G \ast {\bf Z}$ can ever be of weight 1. The negative answer to this question is known as {\bf Kervaire Conjecture}, although Michel Kervaire told us that he asked the question (in particular to Gilbert Baumslag) but did not formulate the conjecture.

It is known that Kervaire Conjecture is true for torsion-free groups. See Klyachko \cite{klya93} and Fenn-Rourke \cite{fero96}.

{\bf Remark.} Let $G$ be a group satisfying the following three conditions:
\\
i) $H_1(G) = 0$.
\\
ii) $H_2(G) = 0$.
\\
iii) $G$ is finitely presented.
\\
In \cite{kerv69} Kervaire proved that these are exactly the groups which are the fundamental group of a homology sphere of dimension $\geq 5$. His arguments are reminiscent of those about knot groups.

\subsection{Groups which satisfy Kervaire conditions}

Theorem 5.2.1 says that Kervaire conditions are necessary for all $n \geq 1$ and Theorem 5.3.1 says that they are sufficient if $n \geq 3$. Question: can we improve on Theorem 5.3.1?
We shall see that the answer is no. It is easy to deduce from Artin's spinning construction (1925) that there are inclusions of the following sets of groups up to isomorphism:

$\lbrace$Groups of 1-knots$\rbrace$ $\subset$ $\lbrace$Groups of 2-knots$\rbrace$ $\subset$ $\lbrace$Groups of 3-knots$\rbrace$ $=$ $\lbrace$Groups of n-knots for $n \geq 3$$\rbrace$

It is known that each inclusion is strict (see below). Among the set of groups (up to isomorphism) which satisfy Kervaire conditions the subset of groups of classical knots is very small.  2-knot groups constitute an interesting intermediate class. Hillman's book \cite{hill89} contains a lot of information about them. 

Let $\Gamma$ be a n-knot group (for any $n$) and let $\Gamma'$ be its commutator subgroup. Of course Condition 1) says that $\Gamma / \Gamma' = {\bf Z}$. Implicitly we often assume that $\Gamma \neq {\bf Z}$. 

\newpage

{\bf Questions:} Can $\Gamma'$ be abelian? If the answer is positive, what abelian groups can be isomorphic to $\Gamma'$?

For $n = 1$, the answer is no. Indeed from Neuwirth's analysis of the commutator subgroup of a 1-knot group \cite{neuw65},  the group $\Gamma'$ contains a free subgroup of rank $\geq 2$ (the fundamental group of an incompressible Seifert surface). 

The case $n = 1$ is in sharp contrast with the case $n \geq 3$, since Hausmann-Kervaire prove in \cite{hake78-1} the following result.

\begin{theorem}
Let $G$ be a finitely generated abelian group. Then there exists a n-knot group $\Gamma$ for $n \geq 3$ such that $\Gamma'$ is isomorphic to $G$ if and only if $G$ satisfies the following conditions:
\\
1. $r_G \neq 1 , 2$;
\\
2. $r_G(2^k) \neq 1 , 2$ for every  $k \geq 1$;
\\
3. $r_G(3^k)$ is equal to 1 for at most one value of $k \geq 1$.
\end{theorem}

\begin{definition}
For any abelian group $G$ (non-necessarily finitely generated)  the rank $r_G$ is the dimension of the ${\bf Q}$-vector space $G \otimes_{\bf Z} {\bf Q}$. For a prime $p$ and an integer $k \geq 1$ the number of factors isomorphic to ${\bf Z} / p^k$ in the torsion subgroup of $G$ is denoted by $r_G(p^k)$. 
\end{definition}

Hausmann and Kervaire also prove that  an abelian $\Gamma'$ can be non finitely generated and give some examples. It seems that a complete classification of non finitely generated abelian  $\Gamma'$ is not known.

The case $n = 2$ is somewhat in between. See \cite{levi77-2} and \cite{hill89}. 

\begin{theorem}
If the commutator subgroup $\Gamma'$ of a 2-knot group $\Gamma$ is abelian, then $\Gamma '$ is isomorphic to one of the following  groups:  ${\bf Z}^3$ , ${\bf Z} \lbrack {1 \over 2} \rbrack$ or a finite cyclic group of odd order. 
\end{theorem}

{\bf Question:} What is the centre $Z({\Gamma}$) of a knot group $\Gamma$? 

The answer for the classical case $n = 1$ is the following:
\\
1. If $Z(\Gamma)$ is non-trivial it is isomorphic to the integers ${\bf Z}$. This result is due to Neuwirth \cite{neuw65}.
\\
2. Burde-Zieschang \cite{buzi66} proved that $Z(\Gamma) = {\bf Z}$ if and only if the knot is a torus knot.

Again the situation is very different if $n \geq 3$. Next result is also due to Hausmann-Kervaire. See  \cite{hake78-2}.

\begin{theorem}
Let $G$ be a finitely  generated abelian group. Then there exists a knot group $\Gamma$ for $n \geq 3$ such that $Z(\Gamma)$ is isomorphic to $G$. 
\end{theorem}

It seems that it is unknown if the centre can be non finitely generated. 

The case $n = 2$ is again in between. For more details see \cite{hill89}.

\begin{theorem}
The centre $Z(\Gamma)$ of a 2-knot group is of rank $\leq 2$. If the rank is equal to 2, then $Z(\Gamma)$ is torsion free.
\end{theorem}

Hillman  gives examples of 2-knot groups $\Gamma$ with $Z(\Gamma)$ isomorphic to ${\bf Z}^2~ ,~ {\bf Z}~ ,~ {\bf Z} \oplus {\bf Z} / 2~ ,~ {\bf Z} / 2$.

\vskip.5in

\section{Knot modules}

\subsection{The knot exterior}

Let us recall J.H.C. Whitehead theorem.

\begin{definition} Let $X$ and $Y$ be two connected C.W. complexes with base-point. A base-point preserving map $f : X \rightarrow Y$   is a homotopy equivalence if there exists a base-point preserving map $g : Y \rightarrow X$ such that the two composition maps $f \circ g$ and $g \circ f$ are homotopic to the identity, keeping the base-points fixed.

\end{definition}

\begin{theorem}
Let $X$ and $Y$ be two connected C.W. complexes with base-point. Let $f : X \rightarrow Y$ be a base-point preserving map. The following three statements are equivalent.
\\
1. $f$ is a homotopy equivalence. 
\\
2. $f$ induces isomorphisms $f_j : \pi_j (X) \rightarrow \pi_j (Y)$ for all integers $j \geq 1$.
\\
3. $f_1 : \pi_1 (X) \rightarrow \pi_1 (Y)$ is an isomorphism and $\widetilde {f}_j : H_j (\widetilde {X} ; {\bf Z}) \rightarrow H_j (\widetilde {Y} ; {\bf Z})$ is an isomorphism for every $j \geq 2$ where $\widetilde {f} : \widetilde {X} \rightarrow \widetilde {Y}$ is a lifting of $f$ between  universal coverings.
\end{theorem}

Knot invariants are provided mostly by the knot exterior $E(K)$ (and its coverings). Since $E(K)$ has the homology of a circle, the homology of $E(K)$ is useless.  Well, let us look at the homotopy groups $\pi_j (E(K))$ for $j \geq 2$. Let us recall that a connected C.W. complex $X$ is said to be {\bf aspherical} if $\pi_j (X) = 0$ for $j \geq 2$. Note that nothing is required on the fundamental group (of course!). By Whitehead's theorem, $X$ is aspherical if and only if its universal covering is contractible. 

{\bf Question:} Can a knot exterior be aspherical?

The answer is very contrasted. By Papakyriakopoulos \cite{papa57} the knot exterior of a 1-knot is always aspherical. But, for $n \geq 2$ if the exterior of a n-knot is aspherical then the knot group is isomorphic to $\bf Z$. This result is due to Dyer-Vasquez \cite{dyva73}. Hence for $n \geq 2$,  a n-knot exterior which is aspherical has the homotopy type of a circle. Next theorem was known when Michel Kervaire wrote his Paris thesis. See \cite{levi65-1}.

\begin{theorem}
(Levine's Unknotting Theorem) Let $K^n \subset S^{n+2}$ be a n-knot such that $E(K)$ has the homotopy type of a circle. Assume that  $n \geq 3$. Then $K$ is the trivial knot. 
\end{theorem}

{\bf Comments.} 1) Levine's result is true in the differential and also in the PL (locally flat) categories. It was already known that the result is true in the topological (locally flat)
 category, thanks to Stallings \cite{stal63}. 
\\
2) For $n = 1$ Levine's result is true by Dehn's Lemma. For $n = 2$ the result is true in the TOP category by Freedman. See Hillman's book \cite{hill89} p.4-5. Indeed the hypothesis $\pi_1(E(K)) = {\bf Z}$ is enough.

Let us now come back to the homotopy groups $\pi_j (E(K))$ for $j \geq 2$. Since the fundamental group is never trivial, these groups are not only  abelian groups but modules over the group ring ${\bf {Z}} \pi_1 (E(K))$.  Except when the fundamental group is infinite cyclic, this group ring is never commutative. It seems that the homotopy groups $\pi_j (E(K))$ are little studied when the fundamental group is not infinite cyclic.
\\
Let us write the infinite cyclic group in a multiplicative way as $T = \lbrace t^k \rbrace$ for $k \in {\bf Z}$. Hence the group ring ${\bf Z} T$ of the infinite cyclic group can be identified with the ring of Laurent polynomials ${\bf Z} \lbrack t , t^{-1} \rbrack$.

\newpage

\begin{definition}
Let $\widehat {E}(K) \rightarrow E(K)$ be the infinite cyclic covering of $E(K)$, i.e. the covering associated with the abelianisation epimorphism of the fundamental group $\pi_1 (E(K)$. A {\bf knot module} is a homology group $H_i (\widehat {E}(K) ; {\bf Z})$ equipped with the ${\bf {Z}}T$-module structure induced by the Galois transformations.
\end{definition}

In the case of 1-knots, knot modules are often called Alexander modules.

In Chapter II of his Paris paper  \cite{kerv65}, Michel Kervaire investigates the following situation: $q \geq 2$ is an integer and $K \subset S^{n+2}$ is a n-knot such that $\pi_i (E(K)) = \pi_i (S^1)$ for $i < q$ and $\pi_q (E(K)) \neq 0$. 

{\bf Problem.} To determine $\pi_q (E(K)) = H_q (\widehat {E}(K) ; {\bf Z})$ as a ${\bf Z}T$-module. 
\\
In today's words Kervaire wishes to understand the ``first non-trivial knot module" but he states his results in terms of $\pi_q (E(K))$. Later he regretted to have adopted  this viewpoint, since many of his arguments are valid for knot modules in general, not only for the first non-trivial one. Thus he found Jerome Levine's general  study of knot modules  in \cite{levi77-1} much better. 

Next lemma contains a useful information. It is proved along the way by Levine in \cite{levi65-1}. It is also a consequence of Poincar\'e duality in the infinite cyclic covering. 

\begin{lemma}
Let $K^n \subset S^{n+2}$ be a n-knot such that $\pi_i (E(K)) = \pi_i (S^1)$ for $1 \leq i \leq q=(n+1) / 2$. Then $E(K)$ has the homotopy type of $S^1,$ and hence the knot is trivial.
\end{lemma}

Thus the problem  of the determination of $\pi_q (E(K))$   splits in three cases which shall be treated separately   in subsections below: 
\\
1. $q < n / 2$ ; 
\\
2. n is even and $q = n / 2$ ; 
\\
3. n is odd and $q = (n+1) / 2$.

{\bf Historical Note.} Clearly, a question is:
How much is detected of a n-knot from its exterior? 
Here are some answers. 

A) A n-knot is determined by its exterior if $n = 1$ by Gordon-Luecke \cite{golu89}. 

B) For $n \geq 2$ there are at most two n-knots with a diffeomorphic exterior.  This result is due to Gluck \cite{gluc61} for $n = 2$, to Browder \cite{brow67} for $n \geq 5$ and to Lashof-Shaneson \cite{lash69} for $n \geq 3$. An equivalent statement is that there are at most two admissible meridians on the boundary $bE(K)$. The proofs of these theorems are stated for embeddings of the standard sphere.
\\
In \cite{cash76-1} Cappell and Shaneson prove that there do exist pairs of inequivalent n-knots with the same exterior.

C) Let us now consider the homotopy type of the exterior.  There are two ways to proceed: either we consider the exterior $E(K)$ or the pair $(E(K) , bE(K))$.  
\\
If $n = 1$ by the asphericity of knots, the homotopy type of $E(K)$ is determined by the fundamental group $\pi_1 (E(K))$ and the homotopy type of the pair $(E(K) , bE(K))$ is determined by the pair of groups (fundamental group, peripheral subgroup). The work of Waldhausen \cite{wald68} and Johannson \cite{joha79} implies then that the homotopy type of the pair $(E(K) , bE(K))$ determines the exterior. 
\\
If $n \geq 3$, Lashof-Shaneson \cite{lash69} prove that the homotopy type of the pair $(E(K) , bE(K))$ determines the differential type of the knot exterior,  provided that the fundamental group of the exterior is infinite cyclic. Their proof assumes that the n-knot is represented by an embedding of the standard sphere.

\subsection{Some algebraic properties of knot modules}

Let us recall some easy properties of the commutative ring with unit ${\bf Z}T$.

A) ${\bf Z}T$ is an integral domain (no zero divisors) but it is not a principal ideal domain, although it might be tempting to think so since ${\bf Q}T$ is. As a consequence there is no decomposition of finitely generated ${\bf Z}T$-modules as a direct sum of cyclic modules. In fact a classification of finitely generated ${\bf Z}T$-modules seems inaccessible. 

B) ${\bf Z}T$ is a noetherian ring (Hilbert basis Theorem). As a consequence, every submodule of a finitely generated ${\bf Z}T$-module is also finitely generated and finitely generated  modules are finitely presented. A surjective endomorphism of a finitely generated ${\bf Z}T$-module  is an isomorphism.

C) ${\bf Z}T$ is a unique factorisation domain (Gauss Lemma). One consequence is that Fitting ideals of a knot module have a gcd. This  gives rise to Alexander invariants.

We now turn to knot modules and recall a definition due to Levine. Notice that the letter  "K"  is the first letter of Knot as well as Kervaire's.  

\begin{definition}
A ${\bf Z}T$-module $H$ is said to be {\bf of type K} if it is finitely generated and if the multiplication by   $(t-1)$  is an isomorphism of $H$. 
\end{definition}

\begin{theorem}
(Kervaire)
Knot modules are of type K.
\end{theorem}

{\bf Proof.} The knot exterior $E(K)$ can be triangulated as a finite simplicial complex. Hence the chain group $C_i(\widehat {E}(K) ; {\bf Z})$ is a finitely generated free ${\bf Z}T$-module. Since the ring ${\bf Z}T$ is noetherian the cycles are  finitely generated as a module and hence the quotient module $H_i(\widehat {E}(K) ; {\bf Z})$ is also a finitely generated module.

Let us now prove that the multiplication by $(t-1)$ is an isomorphism. The key ingredient used by Kervaire is the spectral sequence of a covering map. Since the Galois group is the integers (written multiplicatively) $T$,  this spectral sequence degenerates at the $E_2$ term and it is better replaced by an exact sequence due originally to Serre (Appendix to his thesis) and much used later in knot theory. 

\begin{lemma}
((t-1) lemma.)
Let $X$ be a connected complex (simplicial or C.W.) and let $\pi_1(X) \rightarrow T$ be a surjective homomorphism. Let $p: \widehat {X} \rightarrow X$ be the projection of the infinite cyclic covering associated to this homomorphism. Then $(t-1) : H_i(\widehat {X}) \rightarrow H_i(\widehat {X})$ is an isomorphism for all $i \geq 1$ if and only if $H_*(X)$ is isomorphic to $H_*(S^1)$.
\end{lemma}

{\bf Proof of the lemma.}

It is easy to see that one has a short exact sequence of chain complexes (with integer coefficients)

$$0 \rightarrow C_*(\widehat {X}) \stackrel {t-1}{\longrightarrow} C_*(\widehat {X}) \stackrel {p_*}{\longrightarrow} C_*(X) \rightarrow 0$$

Its associated homology long exact sequence is the object to be considered.

$$\dots \rightarrow H_{i+1}(X) \rightarrow H_i(\widehat {X})  \stackrel {t-1}{\longrightarrow} H_i(\widehat {X})  \stackrel {p_i}{\longrightarrow} H_i(X) \rightarrow \dots $$

This exact sequence implies that $(t-1) : H_i(\widehat {X}) \rightarrow H_i(\widehat {X})$ is an isomorphism for all $i \geq 2$ if and only if $H_*(X)$ is isomorphic to $H_*(S^1)$. The argument also works for $i = 1$, since $ (t-1) : H_0 (\widehat {X})  {\longrightarrow} H_0 (\widehat {X})$ is the zero homomorphism, the connecting homomorphism $\delta : H_1 (X)  {\longrightarrow} H_0 (\widehat {X})$ is a isomorphism and $p_1=0$. 

{\bf End of proof of the lemma.}

We apply the lemma to the knot exterior $E(K)$ which has the homology of a circle.

{\bf End of proof of the theorem.}

An easy consequence of the theorem is that a knot module is a torsion ${\bf Z}T$-module.  

\begin{definition}
Let $H$ be  a ${\bf Z}T$-module of type K. We write  $Tors(H)$ for the {\bf Z}-torsion sub-module of $H$ and $f(H)$  (following Levine) for the quotient module $H / Tors(H)$.
\end{definition}

\begin{theorem}
(Kervaire Tors(H))
Let $H$ be a ${\bf Z}T$-module of type K. Then the cardinal of $Tors(H)$ is finite.
\end{theorem}

For a proof see Kervaire's  Lemme II.8 in \cite{kerv65} or Levine's Lemma 3.1 in \cite{levi77-1}.

{\bf Caution.} $f(H)$ is finitely generated as a module. A lot of examples from classical knot theory show that it is not finitely generated as an abelian group, in general.

Here are some deeper properties of ${\bf Z}T$.

D) ${\bf Z}T$ is of (global) homological dimension 2. For details see \cite{macl63} Chap. VII.  This means that the following three equivalent properties are satisfied.
\\
(i) $Ext_{{\bf Z}T}^{k}(H,H') = 0$ for all ${\bf Z}T$-modules $H$ and $H'$ and all integers $k >2$
\\
(ii) For every ${\bf Z}T$-module $H$ and for every exact sequence with projective ${\bf Z}T$-modules $P_{i}$ for $i = 0 , 1$:
\\
$$0 \rightarrow C_2 \rightarrow P_1 \rightarrow  P_0 \rightarrow H \rightarrow 0$$
\\
the module $C_2$ is also projective.
\\
(iii) For every ${\bf Z}T$-module $H$, there exist projective resolutions  of length at most $2$.

E) ${\bf Z}T$ satisfies Serre Conjecture: every finitely generated projective ${\bf Z}T$-module is free. See \cite{swan78}.

A consequence is the following.

\begin{proposition}
Let $H$ be a knot module. Let $F_1 \stackrel {\Pi}{\longrightarrow} F_0 \rightarrow H \rightarrow 0$ be  a finite presentation of $H$. Then the kernel of $\Pi$ is a finitely generated free ${\bf Z}T$-module. Hence, there exists a resolution  $0 \rightarrow F_2 \rightarrow F_1 \rightarrow F_0 \rightarrow H \rightarrow 0$ of $H$ where the  $F_i$ are finitely generated free ${\bf Z}T$-modules, for $i = 0 , 1 , 2$.
\end{proposition}

Here is a pretty result of Levine. See Proposition 3.5 in  \cite{levi77-1}.

\begin{theorem}
Let $H$ be a ${\bf Z}T$-module of type K. Then the following statements are equivalent.
\\
(1) $Tors(H) = 0$.
\\
(2) $H$ is of homological dimension 1.
\\
(3) $H$ possesses a square presentation matrix.
\end{theorem}

\subsection{The q-th knot module when $q < n / 2$}

The theorem proved by Kervaire (\cite{kerv65}) is the following. 

\begin{theorem}
Let $1 < q < n / 2$ and let $H$ be a ${\bf Z}T$-module. Then there exists a n-knot $K^n \subset S^{n+2}$ with $H_i (\widehat {E}(K) ; {\bf Z}) = 0$ if $0 < i < q$ and $H_q (\widehat {E}(K) ; {\bf Z}) = H$ if and only if $H$ is of type K.
\end{theorem} 

{\bf Proof of the theorem.}

We already know that type K is necessary.  So let us prove that this condition is sufficient. The proof shows that the knot can be realised by an embedding of the standard differential sphere.   In the proof we shall need the following statement of Hurewicz theorem.

\begin{theorem}
(Hurewicz)
Let $X$ be a connected C.W. complex such that $\pi_i (X) = 0$ for $1 \leq  i < k$. Then $\pi_k (X) \rightarrow H_k (X , {\bf Z})$ is an isomorphism and $\pi_{k+1} (X) \rightarrow H_{k+1} (X , {\bf Z)}$ is onto.
\end{theorem}

The proof follows the same steps  as the proof of  the characterization of the knot group. 

\begin{proposition}
Let $H$ be a ${\bf Z}T$-module of type K. Then there exists a connected finite C.W. complex  $Y$ of dimension $\leq (q + 2)$ such that:
\\
i) $\pi_1 (Y) = {\bf Z}$.
\\
ii) $H_i (\widehat {Y} ; {\bf Z}) = 0$ for $1 < i < q$
\\
iii) $H_q (\widehat {Y} ; {\bf Z}) = H$
\\
iv) $H_i (\widehat {Y} ; {\bf Z}) = 0$ for $q < i$. 
\end{proposition}

If we admit the proposition, the proof of the theorem is achieved as follows.

The homology exact sequence of the infinite cyclic covering $\widehat {Y}$ reveals that $Y$ has the homology of a circle. We choose a simplicial complex of dimension $(q+2)$ which has the same homotopy type of $Y$ and denote it still by $Y$.  It can be piece-wise linearly embedded in ${\bf R}^{n+3}$ since its top dimension homology vanishes. This is an easy consequence of Arnold's Shapiro's Theorem 7.2 in \cite{shap57}.  We denote  by $M$ the boundary of a regular neighbourhood $N$ of the embedded $Y$ in ${\bf R}^{n+3}$.  The PL  $(n+2)-$manifold $M$ can be smoothed by Hirsch-Whitehead. It has the same (q+1)-type as $Y$. The proof is achieved by surgery  on $M$ as in the characterization of the knot group, since $M$ is a knot-manifold.

{\bf End of the proof of the theorem.}

{\bf Proof of the proposition.}

We choose a finite presentation of the ${\bf Z}T$-module $H$: $F_{q+1} \stackrel{d_{q+1}}{\longrightarrow} F_q \rightarrow H \rightarrow 0$. We know from the preceeding subsection that the kernel $Ker(d_{q+1})$ is free and finitely generated. We then have a free resolution

\centerline{$0 \rightarrow F_{q+2} \stackrel{d_{q+2}}{\longrightarrow} F_{q+1} \stackrel{d_{q+1}}{\longrightarrow} F_q \rightarrow H \rightarrow 0$}

We denote by $r_i$ the ${\bf Z}T$-dimension of $F_i$ for $i = q+2 , q+1 , q$. 

We shall construct the C.W. complex $Y$ skeleton by skeleton, using the chosen  resolution of $H$. We begin with a 0-cell $*$ and attach to it one 1-cell. We thus obtain a circle. We attach to the  0-cell $*$,  $r_q$ cells of dimension $q$. We thus obtain a wedge of $S^1$ with $r_q$ spheres of dimension $q$. We denote this complex by  $Y_q$. The infinite cyclic covering $p_q : \widehat {Y}_q \rightarrow Y_q$ is also its universal covering. By construction $C_q(\widehat {Y}_q ; {\bf Z})$ is isomorphic to $F_q$ as a ${\bf Z}T$-module. 
\\
Let $\lbrace \epsilon_{q+1}^1 , \epsilon_{q+1}^2 ,   \dots , \epsilon_{q+1}^{r_{q+1}}  \rbrace$ be a ${\bf Z}T$-basis of $F_{q+1}$. Consider the element $d_{q+1}(\epsilon_{q+1}^{i}) \in F_q = H_q (\widehat {Y}_q ; {\bf Z}) = \pi_q(\widehat {Y}_q) = \pi (Y_q)$. For each $1 = 1 , \dots , r_{q+1}$ we attach a $(q+1)$-cell $e_{q+1}^i$ to $Y_q$ by the element $d_{q+1}(\epsilon_{q+1}^{i})$. We thus obtain a C.W. complex $Y_{q+1}$ of dimension $(q+1)$. 
\\
Let $p_{q+1} : \widehat {Y}_{q+1} \rightarrow Y_{q+1}$ be the universal covering. Far from dimensions $0$ and $1$ the chain complex of $\widehat {Y}_{q+1}$ is isomorphic over ${\bf Z}T$ to $0 \rightarrow F_{q+1} \rightarrow F_q \rightarrow 0$. Hence we have $H_q(\widehat {Y}_{q+1} ; {\bf Z}) = H$ and $H_{q+1} (\widehat {Y}_{q+1} ; {\bf Z}) = Ker (d_{q+1}) = F_{q+2}$. 

Let us now consider the exact homology sequence of the infinite cyclic covering $p_{q+1} : \widehat {Y}_{q+1} \rightarrow Y_{q+1}$. We have
\\
$\dots \rightarrow F_{q+2} = H_{q+1} (\widehat {Y}_{q+1} ; {\bf Z}) \stackrel {t-1}{\longrightarrow} F_{q+2} = H_{q+1} (\widehat {Y}_{q+1} ; {\bf Z}) \stackrel {p_{q+1}}{\longrightarrow} H_{q+1}(Y_{q+1} ; {\bf Z}) \rightarrow  H = H_q (\widehat {Y}_{q+1} ; {\bf Z}) \stackrel {t-1}{\longrightarrow} H = H_q (\widehat {Y}_{q+1} ; {\bf Z}) \rightarrow \dots$

By hypothesis, $H$ is a ${\bf Z}T$-module of type K. Hence $(t-1) : H_q(\widehat Y_{q+1} , {\bf Z}) \rightarrow  H_q(\widehat Y_{q+1} , {\bf Z})$ is an isomorphism. We thus get an exact sequence:

$$F_{q+2} = H_{q+1}(\widehat Y_{q+1} , {\bf Z}) \stackrel {t-1}{\longrightarrow} F_{q+2} = H_{q+1}(\widehat Y_{q+1} , {\bf Z}) \stackrel {p_{q+1}}{\longrightarrow} H_{q+1}(Y_{q+1} ; {\bf Z}) \rightarrow 0  ~~~~~~  (\star) $$

Since $Coker ((t-1) : {\bf Z}T \rightarrow {\bf Z}T)$ is isomorphic to the integers ${\bf Z}$ the exact sequence   is isomorphic to:

$({\bf Z}T)^{r_{q+2}} \rightarrow ({\bf Z}T)^{r_{q+2}} \rightarrow {\bf Z}^{r_{q+2}} \rightarrow 0$

{\bf Claim.} The Hurewicz homomorphism $\pi_{q+1}(Y_{q+1}) \rightarrow H_{q+1} (Y_{q+1} ; {\bf Z})$ is onto.

{\bf Proof of the claim.} The Hurewicz homomorphism is equal to the composition of the following homomorphisms: ~
$\pi_{q+1}(Y_{q+1}) \rightarrow \pi_{q+1}(\widehat {Y}_{q+1}) \rightarrow H_{q+1} (\widehat {Y}_{q+1} ; {\bf Z}) \rightarrow H_{q+1} (Y_{q+1} ; {\bf Z})$. 

The first homomorphism is the inverse of the isomorphism on homotopy groups induced by the covering projection. The second homomorphism is onto by the second part of the Hurewicz theorem. The third homomorphism is onto by the sequence  $(\star) $ above (it is here that the hypothesis ``$H$ is of type K" is used).

{\bf End of proof of the claim.}

Let $\lbrace  \epsilon_{q+2}^1 , \epsilon_{q+2}^2 , \dots , \epsilon_{q+2}^{r_{q+2}}   \rbrace$ be a free basis of $F_{q+2}$. Consider the element $p_{q+1}(d_{q+2}(\epsilon_{q+2}^j)) \in H_{q+1} (Y_{q+1} ; {\bf Z})$. Since the Hurewicz homomorphism $\pi_{q+1}(Y_{q+1}) \rightarrow H_{q+1} (Y_{q+1} ; {\bf Z})$ is onto, we can attach to $Y_{q+1}$ a (q+2)-cell along the element $p_{q+1}(d_{q+2}(\epsilon_{q+2}^j))$. We do this for $j = 1 , 2 ,..., r_{q+2}$ and thus obtain a (q+2)-dimensional  complex $Y_{q+2}$. Note that the elements $p_{q+1}(d_{q+2}(\epsilon_{q+2}^j))$ for $j = 1 , 2 , \dots , r_{q+2}$ constitute a basis of the free abelian group $H_{q+1} (Y_{q+1} ; {\bf Z})$. The complex $Y_{q+2}$   is the complex $Y$ we are looking for. Indeed, by construction the ${\bf Z}T$-chain complex of $\widehat {Y}_{q+2}$ is isomorphic to $0 \rightarrow F_{q+2} \rightarrow F_{q+1} \rightarrow F_q \rightarrow 0$ (far from dimensions 0 and 1).  Therefore the homology modules of $\widehat {Y}_{q+2}$ are as annonced. 
{\bf End of proof of the proposition.}

{\bf Remark.} Our proof follows the same lines as Kervaire's, with some improvements brought by Levine.
\\
1) The spine $Y$ replaces Kervaire's direct construction by surgery of the manifold $M$ which is then surgerised on the generator of the fundamental group to produce the n-knot.
\\
2) Kervaire constructs a n-knot with first non-trivial knot module isomorphic to $H$. Since he does not use homological dimension nor the fact that projective ${\bf Z}T$-modules are free, he cannot claim that this knot has only one non-trivial knot module (see the Paris paper p.243 lines 4 to 9). But as shown by Levine  in Knot Modules (see Lemma 9.4), it is true that this   knot has only one non-trivial knot module .

\subsection{Seifert hypersurfaces}
Let $K^n \subset S^{n+2}$ be a n-knot. We denote by $N(K)$ a closed  tubular neighbourhood of $K^n$ and $E(K)$ its exterior. Moreover we choose an oriented meridian $m$ on the bounbary of $E(K)$.

\begin{definition} A {\bf Seifert hypersurface} of $K^n$ is an oriented  and connected  $(n+1)$-submanifold  $F^{n+1}$ of  $ S^{n+2}$  which has $K^n$ as oriented boundary.
\end{definition}

{\bf Basic facts:}

I)  A $n$-knot  has always  Seifert hypersurfaces.
This basic property of $n$-knots was well-known  to  the topologists in the $60'$. We give  in Section 11 (Appendix II),  a proof of the existence of Seifert hypersurfaces  in a more general situation. To obtain
a  Seifert hypersurface $F^{n+1} $ of $K^n$, we prove  the existence of  a differentiable map $\psi :\ E(K) \rightarrow {\bf S}^1$ such that $\psi$ is of degree $+1$ on $m$  and such that the restriction $\psi'$ of $\psi$ on $bN(K)$ extends to a fibration $\tilde \psi : N(K) \rightarrow B^2$  where  $ K^n=\tilde \psi ^{-1}(0)$.
 We choose a  regular value  $z$ of $\psi$. We  take the following notations: $F=( \psi ^ {-1} (z) )$  and $K'=bF$. By construction $K'$ is isotopic to $K^n$ inside $N(K)$.  We can choose in $N(K)$  a collar $C$ such that $bC=K^n\cup (-K')$.  By construction $ F^{n+1}= F\cup C$  is a Seifert hypersurface of $K^n$. But to study $E(K)$ it is more convenient to consider  its deformation retract $F=F^{n+1}\cap E(K)$. We  also say that $F$ is  a  Seifert hypersurface of $K^n$.

 II)   As  $F$ is oriented,   $F$  has a trivial normal bundle. Hence, we can choose an embedding $\alpha :I \times F \rightarrow E(K)$ where
   $I = \lbrack -1 , +1 \rbrack$, such that  the image $\alpha (I\times F)=N(F)$  is  an oriented compact tubular  neighbourhood of $F$ in $S^{n+2}$ and  $F= \alpha (\lbrace 0 \rbrace  \times F)$. We use the following notation: $F_+=\alpha (\lbrace +1 \rbrace \times F)$ (resp. $ F_-=\alpha (\lbrace -1  \rbrace \times F)$). If  $X \subset F$ let $i_+ (X)=\alpha (\lbrace +1 \rbrace \times X),(resp.,i_-(X)=\alpha (\lbrace -1 \rbrace \times X) )$. 
     So we obtain:$$i_+:F\rightarrow E(F) ,   (resp. i_-:F\rightarrow  E(F)).$$
   Let   the exterior  of $F$ be defined as follows: $E(F)= E(K)\setminus \alpha ( (\rbrack  -1 , +1\lbrack ) \times F)$.

    As the boundary of $F$ is a homotopy sphere,  Poincare duality for $F$, followed by  Alexander duality for $F$ and $E(F)$, implies the following proposition (very often used by Kervaire):

     \begin{proposition} Let $1\leq j \leq n$. Then $ H_{j}(F , {\bf Z})$  and $H_{j}(E(F) , {\bf Z})$ are isomorphic.

     \end{proposition}

     III) By construction,  $F$ enables  us to  obtain a presentation of  knot-modules as follows.
     
 The map  $\psi $ induces the abelianisation epimorphism of the fundamental group 
$ \rho  :\pi_1 (E(K)) \rightarrow \bf{Z}$.
Let  $ p : \widehat {E}(K) \rightarrow E(K)$ be   the covering associated to $ \rho $. By definition,   $p$ is the infinite cyclic covering associated to $K^n$.

By pull-back we  obtain a differentiable map $ \widehat \psi: \widehat{E}(K) \rightarrow   \bf {R} $ such that  $  \psi \circ  p  = exp \circ \widehat{\psi },$ where exp is the exponentiel map. Then there exist a regular fiber $\widehat F$ of $\widehat  \psi$ which is diffeomorphic
by $p$ to $F$. 

    As  $p^{-1}(E(F))$ has infinitely many connected copies, we  choose one of them and denote it by $E_0(F)$. The others are obtain from $E_0(F)$ by the action of the Galois group $T$ isomorphic to  $\bf {Z}$. The restriction  $p_0$ of $p$ on $E_0(F)$ is a diffeomorphism on $E(F)$. Then  $p_0^{-1} : E(F)\rightarrow E_0(F) $ induces an embedding $ \widehat {i} : E(F)\rightarrow   \widehat {E}(K) $.
     
    We will  also denote $i_+,i_-, \widehat {i} $  the homomorphisms   induced on the homology groups  by the maps $i_+,i_-, \widehat {i} .$

Levine  presents knot-modules   with the help of the following theorem (see  \cite{levi77-1}, 14.2):     
\begin{theorem}
 Let $0<j<n+1$, the following sequence is an exact sequence of  ${\bf Z}T$-modules :
$$(*)~~~  0 {\longrightarrow}  H_{j}(F , {\bf Z}) \otimes {\bf Z}T  \stackrel {d_j}{\longrightarrow}  H_{j}(E(F) , {\bf Z})  \otimes {\bf Z}T \stackrel { \iota_j}{\longrightarrow} H_{j}( \widehat {E}(K) ; {\bf Z}) \rightarrow 0, $$
 where,  for $x\in    H_{j}(F, {\bf Z})$,  $y\in   H_{j}(E(F),{\bf Z} $ and $a(t)\in {\bf Z}T$, we  have  $ d_j (x\otimes a(t))=(i_+(x)\otimes ta(t))-(i_-(x)\otimes a(t))$, and  $ \iota_j (y\otimes a(t))=a(t)( \widehat {i}  (y))$.

\end{theorem}

We now introduce Levine's concept of simple knots. In \cite{levi70} he  proves the following theorem:

\begin{theorem} Let $1< q \leq (n+1)/2$.  If   a n-knot    $K^n \subset S^{n+2}$ is   such that $\pi_i (E(K)) = \pi_i (S^1)$ for $1 \leq i < q$,  then   $K^n$ has a  Seifert hypersurface $F^{n+1}$ which is $(q-1)$-connected.  Moreover, we  can find  a  $(q-1)$-connected Seifert hypersurface  such that the embeddings $\widehat {i} \circ i_+, (resp.\widehat {i} \circ i_-)$ induce  injective homomorphisms  $H_q (F,{\bf Z}) \rightarrow   H_{q}( \widehat {E}(K) ; {\bf  Z}) $, such a $F$ is  said to be {\bf  minimal}.
\end{theorem}

\begin{definition}  When $n$ is even, let $q=n/2$, when $n$ is odd let $q=(n+1)/2$.
 A   n-knot,  $K^n \subset S^{n+2}$,  such that $\pi_i (E(K)) = \pi_i (S^1)$ for $ i < q$  is  a {\bf {simple n-knot}}.  
\end{definition}

{\bf Comments:} When $n=1$ or $n=2$ all $n$-knots are simple. When $2\leq q$ and $ n=2q-1$,  the above theorem implies that we can choose a $(q-1)$-connected Seifert hypersurface $F^{2q}$ for $K^{2q-1}$.   In this case,  $F^{2q}$ has the homotopy type of a bouquet of $q$-spheres. If $3\leq q$, this implies  as explained in Section 13 (Appendix IV), that such $F^{2q}$ is a $q$-handlebody.

\subsection{ Odd-dimensional   knots and  the Seifert form}

Let $K^{2q-1} \subset S^{2q+1}$ be an odd dimension knot. We denote by $N(K)$ a closed   tubular neighbourhood of $K^{2q-1}$ and  by $E(K)$ its exterior. We choose a Seifert hypersurface  $F^{2q}$ of $K^{2q-1}$ and we consider  $F=F^{2q}\cap E(K)$.  We keep the notations of Subsection 6.4. 
     
     Let us consider the free $\bf {Z}$-module   $H=(H_q (F ; {\bf Z}))/Tors$, where $Tors$ is the ${\bf Z}-$torsion
subgroup of $H_q (F ; {\bf Z})$. We also denote by $i_+, (resp. \  i_-)$ the homomorphism induces  by $i_+,(resp. \ i_-)$ on  the q-cycles of $F$.

 For   $(x,y)$ in  $H\times H$,  we choose  representative  q-cycles $x'$  and $y'$. We define the  {\bf Z}-bilinear form

$$\mathcal {A} : H \times H\rightarrow {\bf Z}$$

by  $\mathcal {A} (x , y) = \mathcal {L}_{S^{2q+1}} (x' , i_+ (y'))$ where $ \mathcal {L}_{S^{2q+1}} (x' , i_+ (y'))$ is the linking number between the two q-cycles  $x'$ and $i_+(y')$ in $S^{2q+1}$ as defined in  Section 10 (the Appendix I):
\begin{definition} The bilinear form $\mathcal {A}$ is  the {\bf Seifert form} associated to $F^{2q}$.
\end{definition}

{\bf {The Seifert form determines the intersection form }}

Let $\mathcal {I} : H \times H  \rightarrow {\bf Z}$ be the intersection form on $H=(H_q (F ; {\bf Z}))/Tors$. Let $\mathcal {A}^T$ be the transpose of $\mathcal {A}$ defined as $\mathcal {A}^T (x , y) = \mathcal {A} (y , x) = \mathcal {L}_{S^{2q+1}} (y' ; i_+ (x'))$. The following formula is  often quoted up to sign in the literature:

$$(-1)^q \mathcal {I}  = \mathcal {A} + (-1)^q \mathcal {A}^T$$

We briefly present the main steps of the proof, since we shall encounter several useful formulae along the way. We have:
\\
$ \mathcal {L}_{S^{2q+1}} (x' , i_-(y')) = \mathcal  {L}_{S^{2q+1}} (i_+ (x') , y')$ by an easy translation argument,
\\$\mathcal {L}_{S^{2q+1}} (x' , i_-(y')) 
 = (-1)^{q+1} \mathcal {L}_{S^{2q+1}} (y' , i_+ (x'))$ by the symmetry of the linking number,
\\
$(-1)^{q+1} \mathcal {L}_{S^{2q+1}} (y' , i_+ (x')) = (-1)^{q+1} \mathcal {A}^T (x , y)$ by definition of $ \mathcal {A}^T$.

Hence:

$$\mathcal {L}_{S^{2q+1}} (x' , i_+ (y')) - \mathcal {L}_{S^{2q+1}} (x' , i_- (y')) = \mathcal {A} (x , y) + (-1)^q \mathcal {A}^T (x , y)$$

{\bf Lemma.} The left hand side  $\mathcal {L}_{S^{2q+1}} (x' , i_+ (y')) - \mathcal {L}_{S^{2q+1}} (x' , i_- (y'))$ of the last equality is equal to 
$(-1)^q \mathcal {I} (x , y)$.

{\bf Proof.} We have $\mathcal {L}_{S^{2q+1}} (x' , i_+ (y')) - \mathcal {L}_{S^{2q+1}} (x' , i_- (y')) = \mathcal {L}_{S^{2q+1}} (x' , i_+ (y') - i_-(y')) = \mathcal {I}_{S^{2q+1}} (x' , C)$ where $ \mathcal {I}_{S^{2q+1}} (.,.)$ is the intersection number between the $q$-cycle $ x'$ and   a $(q+1)$-chain  $C$ such that $\partial C = i_+ (y') - i_-(y')$ .

We claim that

 $$ \mathcal {I}_{S^{2q+1}} (x' , C) = (-1)^q \mathcal{I} (x , y)$$.

 Here is why. Observe that the intersection points in $S^{2q+1}$ between $x'$ and $C$ are the intersection points in $F$ between $x'$ and $y'$. Let $P$ be one of them. We have to compare the contribution of $P$ to both intersection numbers. Let $(e_1 , \cdots , e_q)$ be a frame at $P$  representing the orientation of $x'$ and let $(\epsilon_1 , \cdots , \epsilon_q)$ be one for $y'$. The contribution of $P$ to $\mathcal{I} (x , y)$ compares the frame $(e_1 , \cdots , e_q , \epsilon_1 , \cdots , \epsilon_q)$ to a frame $(f_1 , \cdots , f_{2q})$ at $P$ which represents the orientation of $F$. Now $C$ can be chosen to be equal to $I \times y'$ since $\partial C = i_{+} (y') - i_-(y')$. Hence the contribution of $P$ to $\mathcal {I}_{S^{2q+1}} (x' , C)$ compares the frame $(e_1 , \cdots , e_q , I , \epsilon_1 , \cdots , \epsilon_q)$ to the frame $(I , f_1 , \cdots , f_{2q})$. To do that we have to move $I$ in $(e_1 , \cdots , e_q , I , \epsilon_1 , \cdots , \epsilon_q)$ into the first place. This introduces the factor $(-1)^q$. QED.

\begin{proposition}
 Let ${\bf {E}} $ be a ${\bf Z }$-basis of  $H=(H_q (F ; {\bf Z}))/Tors$. Let $A$ be the  matrix of  the bilinear form $\mathcal {A}$ in the chosen basis  ${\bf E}$ and  $A^T$ be the transpose of $A$. Then there exists a basis of  $H'= (H_q (E(F) ; {\bf Z}))/Tors$, such that, in the chosen basis,  $A$ is  the matrix of $i_+$ and $((-1)^{(q+1)}A^T) $ the matrix of $i_-$.
\end {proposition}
{\bf Proof.}

 In the homology exact sequence of the pair $ (N(F)\subset {\bf S^{(2q+1)}})$,  the connecting  homomorphism  $$   H_{q+1}({\bf S^{(2q+1)} } , N(F); {\bf Z})  \stackrel {\delta  } {\longrightarrow }  H_{q}(N(F); {\bf Z}), $$ is an isomorphism. Let $\delta ^{-1}$ be its inverse.
 Let $ex$ be the excision  isomorphism and $D$ be Poincare duality. We consider the sequence of isomorphisms:
 $$  H_{q+1}({\bf S^{(2q+1)} } , N(F); {\bf Z}) \stackrel {ex} {\longrightarrow}  H_{q+1}(E(F) ,bE(F);  {\bf Z})   \stackrel {D}{\longrightarrow}  H^q(E(F); {\bf Z}) .$$ 
 Let $H''$ be the free abelian group $( H^q(E(F); {\bf Z}))/Tors. $ The isomorphism $D\circ ex \circ \delta ^{-1}$ induces an isomorphism $\alpha:H\rightarrow H''.$ 
 On the other hand the universal coefficient homomorphism induces an isomorphism $u$ where:
 $$u: H'' \rightarrow  Hom_{{\bf Z}}(H';{\bf Z}).$$
  Let $(x,y)\in (H\times H)$ and let  us  also   write $i_+(y')$  for  the class of $i_+(y')$ in $H'$. By  construction  $\mathcal {A} (x , y) = \mathcal {L}_{S^{2q+1}} (x' , i_+ (y'))= \epsilon ((u\circ \alpha) (x))(i_+(y')),$ where  $\epsilon $ is a  sign  equal to  $1$,  $ (-1)^{(q+1)}$ or $ (-1)^q $, according to  the conventions  of sign   used in  dualities. One can choose a basis for $H''$ such that the matrix of $\alpha $ is equal to $\epsilon I$, where $I$ is the identity matrix. We take the dual basis for $H'$, and then $I$ is the matrix of $u$. Hence, with the chosen basis,   $A$ is the matrix of $i_+:H\rightarrow H'.$
 
 As  $\mathcal {L}_{S^{2q+1}} (x , i_-(y)) = (-1)^{q+1} \mathcal {L}_{S^{2q+1}} (y , i_+ (x))$, we have that   $((-1)^{(q+1)}A^T) $  is the matrix of $i_-$ in the same chosen basis.

Theorem (6.4.1) and the above proposition imply the following corollary:

\begin{corollary}
The matrix $(At+(-1^q)A^T)$ is a presentation matrix of the  ${\bf Z}T$-module  $(H_{q}( \widehat {E}(K) ; {\bf Z}) )/Tors$, where $Tors$  is the ${\bf Z }$-torsion of $ H_{q}( \widehat {E}(K) ; {\bf Z}) $.

\end{corollary}

\subsection{Even-dimensional  knots and the torsion Seifert form}

Let $X$ be a (finite) CW-complex. We denote by $H_j(X)$ the homology group $H_j(X ; {\bf Z})$, by $T_j(X)$ the ${\bf Z}$- torsion subgroup of $H_j(X)$ and by $F_j(X)$ the ${\bf Z}$-torsion-free quotient $H_j(X) / T_j(X)$. 

In his theorem II.2 Kervaire gave necessary and sufficient conditions on a ${\bf Z}T$-module to be the $q$-th knot module of a simple $2q$-knot. These conditions are expressed via a presentation matrix, derived from a minimal Seifert hypersurface. In this sense they are not intrinsic. Later, in ``Knot Modules", Levine coined intrinsic conditions in terms of Blanchfield duality, adequately reinterpreted. It follows from both Kervaire and Levine work that there are no more conditions on the torsion free part $F_q \widehat {E}(K)$ besides the usual ones on a knot module (finite type and multiplication by $(1-t)$ is an isomorphism). The novelty is about $T_q \widehat {E}(K)$. Kervaire shows that it has a presentation matrix which looks formally as the presentaion matrix for $F_q \widehat {E}(K)$ of a simple $(2q-1)$-knot. Here is Kervaire's original statement of his theorem II.2 (with minor changes to adapt notations).

\begin{theorem}
Let $H$ be a ${\bf Z}T$-module and $q$ be an integer $\geq 3$. There exists a knot $K^{2q} \subset S^{2q+2}$ with $\pi_i (S^1) \cong \pi_i (S^{2q+2} \setminus K^{2q})$ for $i < q$ and $H \cong \pi_q (S^{2q+2} \setminus K^{2q})$ if and only if the multplication by $(1-t)$ is an isomorphism of $H$ and $H$ possesses a presentation of the form
$$\left( x_1 , \cdots , x_{\alpha} ; \sum_j (ta_{ij} - b_{ij}) x_j , d_jx_j  \right)$$
where the integers  $d_{k} \in  \bf{ Z}$ with $1 \leq  k \leq \alpha$  and  the square integer matrices $A=(a_{ij})$ and $B=(b_{ij})$ satisfy the relation
$$d_ia_{ij} + (-1)^{q+1}d_jb_{ji} = 0 $$
for all couple of indices $(i , j)$ with $1 \leq i \leq \alpha$ and $1 \leq j \leq \alpha$.
\\
$\lbrack$Some $d_k$'s are allowed to be equal to $0$. They correspond to the torsion free part of the module.$\rbrack$
\end{theorem}

We now tell the story about $T_q \widehat {E}(K)$. This is the part of the module where the $d_k$'s are $\neq 0$ and where the relation between the matrices $A$ and $B$ matters. 

Let $M^m$ be a closed (connected) oriented manifold of dimension $m$. For each integer $i$ such that $0 \leq i \leq m$ there is a pairing 
$$TI_i : T_i(M) \times T_{m-i-1}(M) \rightarrow {\bf Q} / {\bf Z}$$ 
defined as follows. Let $x \in T_i(M)$ and let $z \in T_{m-i-1}$. Represent $x$ resp $z$ by disjoint cycles $\xi$ resp $\zeta$. Suppose that $x$ is of order $d$. Let $\gamma$ be a chain (with integer coefficients) of dimension (i+1) such that $\partial \gamma = d\xi$ Then 
$$TI_i(x ; z) = \frac {1}{d} (\gamma \centerdot \zeta)$$
where  $\centerdot$ denotes the integral intersection number. Main properties of the maps $TI_i$ are:
\\
(0) they are well defined with values in ${\bf Q} / {\bf Z}$;
\\
(1) they are ${\bf Z}-$bilinear;
\\
(2) they are non-degenerate, meaning that the adjoint $TI_i^{\sharp} : T_i (M) \rightarrow Hom_{\bf Z}(T_{m-i-1}(M) ; {\bf Q} / {\bf Z} ) \simeq  T_{m-i-1} (M)$ is an isomorphism;
\\
(3) $TI_i = (-1)^{mi+1} TI_{m-i-1}$

A special case takes place when $m = 2q+1$ and $i = q$. Consider 
$$TI_q : T_q (M) \times T_q (M) \rightarrow {\bf Q} / {\bf Z}$$
This is a bilinear, non-degenerate, $(-1)^{q+1}-$symmetric form on the finite group $T_q (M)$. It is usually called the ``linking form" on the torsion group $T_q (M)$. We prefer to think about it as a torsion intersection form, keeping the name ``linking" for the coupling presented below. See also Wall in \cite{wall67}.

These facts were discovered by several mathematicians in the twenties, including Veblen, Alexander and de Rham. A classical reference is Section 77 of Seifert-Threllfall's book. They made more precise Poincar\'e's duality on the torsion subgroups. Of course, all this was known to Kervaire when he addressed the question of the structure of the q-th knot module of a simple 2q-knot. He used the technique, without mentioning the underlying concepts. 
 
\begin{lemma}
Let $F$ be a  $(2q+1)-$dimensional manifold differentiably embedded in $S^{2q+2}$, with $\partial F$ a homology sphere. Let $E(F)$  be the exterior of $F$ in $S^{2q+2}$. Let $j$ be an integer with $0 \leq j \leq 2q$. Then $H_j(F ; {\bf Z})$ is isomorphic to $H_j( E(F) ; {\bf Z})$.
\end{lemma}

{\bf Proof.} ~ By Poincar\'e duality $H_j(F)$ is isomorphic to $H^{2q+1-j} (F ;  \partial F)$, which is isomorphic to $H^{2q+1-j} (F)$ since $\partial F$ is a homology sphere. By Alexander duality $H^{2q+1-j} (F)$ is isomorphic to $H_{2q+2-(2q+1-j)-1} (E(F)) = H_j(E(F))$.

Alexander duality yields the {\bf torsion linking forms}
$$TL_j : T_j (F) \times T_{2q-j} (E(F)) \rightarrow {\bf Q} / {\bf Z}$$
where $j$ is an integer such that $1 \leq j \leq 2q-1$. It is defined as follows.
Let $x \in T_j (F)$ and $z \in T_{2q-j} (E(F))$. Let $\xi ~, ~\zeta~,~ \gamma~,~d$ be as above.  Let $\eta$ be an integral $(2q-j+1)$-chain in $S^{2q+2}$ such that $\partial \eta = \zeta$. Then
$$TL_j (x , z) = \frac{1}{d} (\gamma \centerdot \eta)$$
Properties of the torsion linking forms are:
\\
(0) they are well defined with values in ${\bf Q} / {\bf Z}$;
\\
(1) they are {\bf Z}-bilinear;
\\
(2) they are non degenerate, meaning that the adjoint provides an isomorphism from $T_j (F)$ to $T_{2q-j} (E(F))$. 

We should remember that cohomology was not available before 1935. Naturally there is a definition of these pairings making use of cohomology and the universal coefficient theorem. 

Now we suppose that $K^{2q}$ is a simple $2q$-knot in $S^{2q+2}$ and that $F$ is a minimal Seifert hypersurface for $K$, as Kervaire did. Following Gutierrez and Kojima (who formalised what Kervaire wrote) we define the {\bf torsion Seifert forms} 
$$T\mathcal {A}_{\pm} : T_q (F) \times T_q (F) \rightarrow {\bf Q} / {\bf Z}$$ 
as follows. Let $x$ and $y$ be in $T_q (F)$. Then 
$$T\mathcal {A}_{\pm} ( x , y ) = TL_q (x , i_{\pm} y)$$

We recall that a Seifert hypersurface $F$ for a simple $2q$-knot is {\bf minimal} if it is $(q-1)$-connected and if $i_{\pm} : H_q (F) \rightarrow H_q (Y)$ are injective. Levine proved (implicitly)  that such hypersurfaces exist (see ``Unknotting spheres ..." middle p.14 ). For a minimal hypersurface, the homomorphisms $i_{\pm} : T_q (F) \rightarrow T_q (Y)$ are isomorphisms since they are injective between finite isomorphic groups. Hence the torsion Seifert forms are non-degenerate in this case.

Kervaire proved that the torsion Seifert forms for even dimension knots have formally much in common with Seifert forms for odd dimension knots. Indeed:

\begin{proposition}
We have the equality $T\mathcal {A}_+ = (-1)^q T\mathcal {A}_{-}^T$ where $T\mathcal {A}_{-}^T$ denotes the transposition of $T\mathcal {A}_-$.
\end{proposition}

{\bf Indications on the proof.} Kervaire chooses a ``basis" $\{ \xi _i\}$ for $T_q (F)$ (with $ \xi_i $ of order $d_i$) and a dual basis $\{x_j\}$ for $T_q (E(F))$. The correspondent Seifert matrices are denoted $(a_{ij})$ for $T \mathcal {A}_+$ and $(b_{ij})$ for $T \mathcal {A}_-$. By definition $i_+(\xi_i) = \sum a_{ij}x_j$ and $i_-(\xi_i) = \sum b_{ij} x_j$. One has $a_{ij} = d_j TL_q (\xi _j , i_+ (\xi _i))$ and $b_{ij} = d_j TL_q( \xi_j , i_-(\xi_i))$.
\\
At the bottom of p.248, Kervaire proves that $d_i a_{ij} = (-1)^q d_j b_{ji}$. Hence
\\
 $d_id_j TL_q (\xi_j , i_+ (\xi_i)) = (-1)^q d_j d_i TL_q (\xi_i , i_- (\xi _j))$.
\\
Both sides are integers. Dividing by $d_i d_j$ we get the equality, with values in ${\bf Q} / {\bf Z}$.

Not surprisingly we get next corollary.

\begin{corollary}
We have the equality $T \mathcal {A}_+ + (-1)^{q+1}T \mathcal {A}_+^T = -TI_q$
\end{corollary}

Gutierrez (in \cite{guti72}) proved that several of these results are true without assuming that the even dimension knot is simple.

Without stating these results explicitly, Kervaire  uses them to get a presentation of the q-th knot module $H_q(\widehat {E} (K) ; {\bf Z})$. Classical arguments provide the exact sequence (For more details see Levine's  ``Knot Modules" p.43) 
$$ 0 \rightarrow H_q (F) \otimes {\bf Z}T \stackrel {\psi}{\longrightarrow} H_q (E(F)) \otimes {\bf Z}T \longrightarrow H_q(\widehat {E} (K) ; {\bf Z}) \rightarrow 0 $$
The homomorphism $\psi$ is classically expressed as $\psi (x \otimes 1) = i_+ (x) \otimes t - i_-(x) \otimes 1$ and hence the two Seifert matrices are the tool to get a presentation of the q-th knot module. It happens that the previous short exact sequence restricts to the short exact sequence (see ``Knot Modules" p.44)
$$ 0 \rightarrow T_q (F) \otimes {\bf Z}T \stackrel {\psi}{\longrightarrow} T_q (E(F)) \otimes {\bf Z}T \longrightarrow T_q(\widehat {E} (K)) \rightarrow 0 $$
Recall that Kervaire proved that $T_q(\widehat {E} (K))$ is a finite group. We have proved that $T_q (F) = T_q (E(F))$ and an easy Mayer-Vietoris argument proves that $T_q (F) = T_q(\widehat {E} (K))$. 
\\
We also have a short exact sequence for the torsion free quotients:
$$ 0 \rightarrow F_q (F) \otimes {\bf Z}T \stackrel {\psi}{\longrightarrow} F_q (E(F)) \otimes {\bf Z}T \longrightarrow F_q(\widehat {E} (K)) \rightarrow 0 $$

The last short exact sequence provides a presentation for the  $\bf {Z}-$torsion free quotient $F_q(\widehat {E} (K))$. The penultimate short exact sequence provides a presentation for the $\bf {Z}-$torsion subgroup $T_q(\widehat {E} (K))$ with the relation $d_i a_{ij} = (-1)^q d_j b_{ji}$ proved by Kervaire. Of course we must take into account that the ``basis" elements are of finite order and add the corresponding relations. Thus the presentation of $T_q(\widehat {E} (K))$ is not square. But the presentation for $F_q(\widehat {E} (K))$ is.

After having established that the $q$-th knot module has presentations which have the form we have just described, Kervaire proved the opposite realisation result. The proof extends over more than three pages (p.249-252). It relies much on Kervaire explicit knowledge of $(q-1)$-connected, parallelisable manifolds of dimension $(2q+1)$ with boundary a homotopy sphere. 

{\bf A brief look at what happened after Kervaire's work on simple $2q$-knots.}
\\
1. Maurizio Gutierrez (in \cite{guti72}) introduced ${\bf Q} / {\bf Z}$ into the picture and generalised some of Kervaire's results for dimensions different from $q$.
\\
2. Jerome Levine (in \cite{levi77-1}) in his exhaustive study on knot modules defined a bilinear, $\pm$ symmetric, non-degenerate form on $T_q(\widehat {E} (K))$ with values in ${\bf Q} / {\bf Z}$ , usually denoted by $\lbrack ~,~\rbrack$. Among many other things he proved a realisation result for torsion knot modules equipped with such a form.
\\
3. Sadayoshi Kojima (in \cite{koji79})  wrote a useful paper which clarifies many concepts.

\vskip.5in

\section{Odd-dimensional  simple links}

\subsection{$q$-Handlebodies}

Let $K^{2q-1} \subset S^{2q+1}$ be a simple  odd dimensional  knot.  As explained at the end of 6.4, we can  choose a Seifert hypersurface  $F^{2q}$ of $K^{2q-1}$ which is a $q-$handlebody.  In this section we consider a more general situation.

\begin{definition}
A $(2q-1)$ {\bf simple   link}  in  $ S^{2q+1}$ is the boundary of a $q-$handlebody $F^{2q}$ embedded in $S^{2q}$: $$bF^{2q}\subset S^{2q+1}$$
\end{definition}

\newpage

We first define what is a handle presentation  of a $q-$handlebody.

\begin{definition}

Let $\phi_j : (S^{q-1} \times B^q)_j  \rightarrow S^{2q-1}$, $1\leq j \leq k$ be $k$ disjoint embeddings.  Let $W^{2q}$ be the quotient of  $$(B^{2q} \ \amalg_j ( B^q \times B^q)_j) $$
 by the identification  $x=\phi_j(x) $, for all    $x\in ( S^{q-1} \times B^q)_j=( bB^q \times B^q)_j $. Let $\pi$ be  the corresponding  quotient map:
  $$\pi : (B^{2q} \ \amalg_j ( B^q \times B^q)_j) \rightarrow W^{2q}$$

The image of  $ \pi ((B^q \times \lbrace 0\rbrace)_j) $ by this identification is the core $C_j$ of the handle
$B_j=\pi (( B^q \times B^q)_j)$, and  $\pi (B^{2q} \ \amalg_j ( B^q \times \lbrace 0\rbrace )_j)$ is  the skeleton $\sigma (W)$ of $W^{2q}.$

A $(2q)$- manifold $F^{2q}$  diffeomorphic to  such  a $W^{2q}$ is a $q-$handlebody  which has $W^{2q}$ as  handle presentation. The  $(q-1)$-dimensional link  $L$ in $S^{2q-1}$  defined by :
 $$ L= \amalg_j \ ( L_j= \phi_j( (S^{q-1} \times  \lbrace 0 \rbrace )_j  ) ) \subset  S^{2q-1}$$
 is the  {\bf attaching link}  of this handle presentation of $F^{2q}.$
 
\end{definition}

{\bf Remark.} A  $q-$handlebody  has the homopoty type of a "bouquet" of $q-$spheres and its boundary  is always $(q-2)-$connected. When $q\geq 2$, a handle presentation is orientable because  the intersection of   each handle with $B^{2q}$ is connected (it is diffeomorphic to $(S^{q-1} \times B^q)$).  If an oriented   handle presentation $W^{2q}$ is embedded in $S^{2q+1}$ it is parallelisable.  We can identify  a tubular neighbourhood $N(W)$ with the   normal interval   bundle of  $W^{2q}$ in $S^{2q+1}$.  Let $i_+$ be a section of this bundle  in the positive direction.  "The" Seifert matrix of the handle presentation $W^{2q}$ is the matrix of the Seifert form   in a particular  basis  $ \lbrace c_j, \ 1\leq j \leq k \rbrace $ of  $H_q(W^{2q} , \bf Z)$ constructed as follows:

  Each  connected component  $L_j$ of the attaching  link of the handle presentation bounds  an  oriented $q-$ball  $D _j$ in $B^{2q}$. Let $c_j$ be the $q-$cycle obtained as the union of $D _j$ with  the core $C_j .$ The corresponding  Seifert matrix  $ A= (a_{ij})$ is given by  $a_{ij}= \mathcal {L}_{S^{2q+1}} (c_i ; i_+ (c_j))$.

{\bf Remark.}  Let $A$ be "The" Seifert matrix of a  handle presentation of $F^{2q}$.  Assume that  $I=(-1)^q(A+ \ (-1)^q A^T)$ is unimodular, then by Poincare duality:
\\
  If $q=1$,   the boundary of $F^{2}$ is $S^1$ and $bF^{2}\subset S^{3}$ is a classical knot.
  \\
   If  $q=2$ the boundary of $F^{2q}$ is a homology  sphere.
   \\
    If  $q\geq 3$ the boundary of $F^{2q}$ is a homotopy sphere and $bF^{2q}\subset S^{2q+1}$ is a knot.

\subsection{The realization theorem for Seifert matrices}

\begin{theorem}(Realization)
Suppose that $q \geq 2$. Let  $A$ be a square  $(k\times k)$ matrix with integral coefficients. Then there  exists  some  oriented  handle presentation  $W^{2q}$ embedded in $S^{2q+1}$ which has $A $ as Seifert matrix.
\end{theorem}

{\bf Comments:} M. Kervaire in \cite{kerv65} states this result as Th\'eor\`eme II.3 p.235 and proves it in p.255-257. In his proof, Kervaire constructs first an abstract handlebody, essentially by following the data     
given by the intersection form. He then embeds the handlebody without control in $S^{2q+1}$. Finally he modifies the embedding in order to realize the Seifert matrix. Levine refers to Kervaire. Here, we propose a different approach. We construct the handlebody directly in $S^{2q+1}$, following step by step the data given by the Seifert matrix. This construction works for any integral square matrix $A$. And, when  $q\geq 3$,  the boundary of the  obtained embedded $q-$handelbody is a knot if and only if  $I=(-1)^q(A+ \ (-1)^q A^T)$ is unimodular. 
\\
  Moreover,  if  $I=(-1)^q(A+ \ (-1)^q A^T)$ is unimodular  the construction  given in the following proof also works when $q=1.$

{\bf Proof} (Realization Theorem.)
We have a $(r\times r)$-matrix $A=(a_{ij}), \ a_{ij} \in {\bf Z} $ and  the matrix $I=(-1)^q(A+ \ (-1)^q A^T )$. 

1) {\bf First  step.}  We choose  $r$ disjoint   embeddings $\phi'_j : (S^{q-1})  \rightarrow S^{2q-1}$ such that:
  $$(-1)^q(a_{ij} +(-1)^qa_{ij})= \mathcal {L}_{S^{2q-1}} (L_i ; L_j)$$
   where $ 1\leq i<j\leq r, \ L_i= \phi'_i( S^{q-1}) ,\ L_j= \phi'_j( (S^{q-1}) $.
   
   We choose embeddings  $\phi'_j$ such that  $ L_j= \phi'_j( S^{q-1}  )  \subset  S^{2q-1}$  is a trivial knot (this is automatic if  $q\neq  2$). Each  $L_j$ bounds  a  differetial ball  $D _j$ in $B^{2q}$. We have:
    $$(-1)^q(a _{ij} +(-1)^qa_{ij})= \mathcal {L}_{S^{2q-1}} (L_i ; L_j) = I_{B^{2q}} (D_i  ;   D_j),  \ 1\leq i \leq  j \leq r. $$
   
2)  { \bf Second step.}

  Let $\mathcal S^{2q}$ be the boundary of  some  $(2q+1)$-ball embedded  in $S^{2q+1 }$. Let $i'_+$ be a  nornal vector field defined on  $\mathcal S$ and pointing to the exterior of the ball.  We identify the sphere  $ S^{2q-1}$ which contains  the  attaching  link $ L=\amalg_j \ ( L_j= \phi'_j( S^{q-1})  )$ with an  equator of the  sphere $\mathcal S$. This equator divides $\mathcal S$ in two hemispheres, the North  hemisphere denoted by $\bf {B}$ and the South hemisphere denoted by $\bf {B'}$. 
  \\
   We choose identifications of $B^{2q}$ with $\bf {B}$ and $\bf {B'}$.   We thus obtain $((\cup _j D _j) \subset   \bf{B})$ and  $((\cup _j  D' _j) \subset   \bf {B'})$ where  $D'_j$ are  the images, in $\bf {B'}$,  of $D_j$  by the chosen identifications.  Using  the normal vector field $i'_+$ we push the $D'_j$ slightly outside $\bf {B'}$ in $S^{2q+1}$ (keeping them fixed on $L_j$) getting disjoint balls $C'_j$. Each union $ c'_j=D_j \cup C'_j$ is a  $q$-sphere embedded in $S^{2q+1}$.  Let $N(c'_j)$,  be a  tubular neighbourhood  of $c'_j$   in $S^{2q+1}$ such that  the $N(C'_j)$ are disjoint. The normal vector field $i'_+$  restricted  on the cores $C'_j$, is a nowhere vanishing section of $N(C'_j)$. Let $B'_j$ be its  orthogonal complement  in $N(C'_j).$ The union 
$W'=  {\bf{B}} \cup_j B'_j$ is an oriented   $q$-handelbody  embedded in $S^{2q+1}$ which has a presentation with    cores $C'_j$.
Moreover $b(N(c'_j))$ is diffeomorphic to $c'_j\times S^q$. 

The handlebody $W'$ is essentially flat, hence it is not the one we are looking for to realize  the matrix $A$.  In particular,  we have   $\mathcal {L}_{S^{2q+1}} (D_j \cup D'_j; i'_+( D_j \cup D'_j) )=0$ by construction. This  implies that    $ \mathcal{A}(c'_j,c'_j)= \mathcal {L}_{S^{2q-1}} (c'_j; i'_+(c'_j))=0$.

 3) { \bf Third step: realization of the non diagonal coefficients.}
 
 From the embedded $W'\subset S^{2q+1}$ we  only keep the skeleton $\sigma_1 = {\bf{B}}  \cup  (\amalg _j C'_j)$.  For $1\leq   i<j\leq r$, let $a'_{ij}= \mathcal {L}_{S^{2q-1}} (c'_i ; i'_+(c'_j))$ and   let $b_{ij}=a_{1j}-a'_{ij}$. We denote by $C_{1j} $  a well chosen  connected sum of $C'_j$ with the boundary of  $b_{1j}$ distinct fibers of the normal bundle $N (C'_1)$ (which is isomorphic to $C'_1\times B^{q+1}$). These  fibers are oriented by the sign of $b_{1j}$.  We perform the  connected sum cautiously   in such a way that the new skeleton $\sigma_2=  {\bf{B}}  \cup C'_1 \cup _{2\leq j \leq r}C_{1j}$ has disjoint cores  and has  new cycles $c_{1j}=D _j \cup C_{1j}$ which satisfy:

$$ \mathcal {L}_{S^{2q-1}} (c'_1 ; i'_+(c_{1j}))=a_{1j}, \ 2\leq j \leq r.$$
  We also  denote  by $i'_+$ any extension  of the positively oriented  section $i'_+$,  of the normal bundle of ${ \bf{B}}$,  to the new cores $C_{1j}.$ The indeterminacy in  the choice of these extensions does not affect the value of the non diagonal coefficients of $A$.
 
 \newpage
 
  By the same construction, we obtain   for all $i, \ 2\leq i \leq (r-1),$ such a  skeleton $\sigma_{i+1}$.   The cores of $\sigma_{i+1}$ are constructed inductively  from the cores  $C_{(i-1)j}$ of  $\sigma _i$ such that  $C_{i1}= C_{(i-1)1}=...=C_{11}=C'_1,  \ C_{i2}=C_{12}, \ ... , \ C_{ii}=C_{(i-1)i} \ $   and,  for $i<j\leq r, \ C_{ij}$ are a  connected sum of $C_{(i-1)j}$ with the boundary of $b_{ij}=a_{ij}-a'_{ij}$ fibers of the normal bundle $N(C_{ii})$. At the end of this step,  we have a skeleton $\sigma _r$ with a basis  of $q$-cycles $c_{rj} = D_j \cup C_{rj} , \ 1\leq j\leq r.$ 
 
 We fix the skeleton $\sigma =\sigma _r$, and we simplify the notations as follow:  $C_j=C_{rj}$ and  $c_j=D_j \cup C_j$. The  construction of $\sigma $ implies:
  $$ (*) \   \  \mathcal {L}_{S^{2q-1}} (c_{i} ; i'_+(c_{j}))=a_{ij}, \ 1\leq  i < j \leq r.$$
  
  4) { \bf Fourth  step: realization of the  diagonal coefficients.}
 
 The normal bundle  of $c_j$  in $S^{2q+1} $ is  trivial. We identify a tubular neighbourhood $N(c_j)$ of $c_j$ in $S^{2q+1} $ with the associated normal disc bundle.  The restriction of $i'_+$  on $c_j$ is  a nowhere vanishing section of this bundle. Let $a'_{jj}$ be defined as follows: 
$$  a'_{jj} = \mathcal {L}_{S^{2q-1}} (c_{j} ; i'_+(c_{j}))  \ 1\leq  j \leq r.$$

We modify $i'_+$ (only  on $C_j$)  to obtain a new section $i_+$ such that:
$$a_{jj}=  \mathcal {L}_{S^{2q-1}} (c_{j} ; i_+(c_{j})) , \ 1\leq   j \leq r.$$

Let $B_j$ be the orthogonal complement to $i_+$ in $N(C_j)$. The handle presentation  

$$W^{2q}={ \bf{B}}\cup_j  B_j$$
has $A$ as Seifert matrix.

{\bf End of proof of the realization theorem for Seifert matrices}

{\bf Remarks.} In  the handlebody    $W^{2q}={ \bf{B}}\cup_j  B_j$,  each connected  component $L_j$ of the attaching  link $(L\subset  b{ \bf{B}})$ has a tubular neighbourhood $N(L_j)= B_j \cup \  b{ \bf{B}}$. One identifies   $ N(D _j)$  with the normal disc bundle of  $D _j$ in ${ \bf{B}}$ and $B_j$ with the normal disc bundle of  $C_j$ in $W^{2q}$.  Both of them have a unique  trivialization . The   trivialization of $B_j $ restricted on $N(L_j)$  gives  disjoint embeddings  $\phi_j : (S^{q-1} \times B^q)_j  \rightarrow b{ \bf{B}}$, $1\leq j \leq r$,  which are used to  define the handlebody presentation of $W^{2q}$ and his  $q$-handles $B_j.$

 Moreover $B_j\cup ({ \bf{B}}\cap N(D_j ))$ is isomorphic to the  normal disc  bundle  of  $c_j$ in $W^{2q}$ which is  is stably trivial,  but  can be non-trivial.  The  trivializations  of its restriction on $D_j$ and $C_j$ can be glued on $L_j$.  This gluing is well defined, up to isomorphism,  by an element of $\pi_{q-1}(SO_q)$ . When $q$ is even this element is given by the Euler number which is equal to 
$ I(c_j;c_j)= 2\times \mathcal{A}(c_j,c_j)=  2\times \mathcal {L}_{S^{2q-1}} (c_j; i_+(c_j))$.

\subsection{Levine's classification of embeddings of handlebodies in codimension one}

\begin{theorem}  (Levine) Let $F^{2q}$ be a $q$-handelbody embedded in $S^{2q+1}$ and let $\mathcal A $ be its Seifert form.   If $q\geq 3$, $\mathcal A$ classifies the embedding $F^{2q}\subset S^{2q+1}$ up to orientation preserving diffeomorphism of the pair $F^{2q}\subset S^{2q+1}$.
\end{theorem}

{\bf Comments: }  Of course an  isotopy between two embedded  handlebodies produces an isotopy between their boundaries. In \cite{levi70}, Levine proves the above theorem and uses it  to  prove the following classification theorem:

\begin{theorem} If $q\geq 2$, two simple knots are isotopic if and only if they have $S-$equivalent Seifert forms.
\end{theorem}

 The notion of $S-$equivalence is well established for knots  in   \cite{levi70}. Today there is no available generalisation of $S-$equivalence  to links. When $q=2$, Levine uses that $K^3$ is diffeomorphic to $S^3$ to apply an  important C.T.C. Wall result on the closed simply connected four manifolds. When $q\geq 3$, there is no condition on  the boundary of $F^{2q}$ ($bF^{2q}$ has not to be a homotopy sphere).

{\bf Sketch of proof.}  
 Let $F^{2q}$ and  $F'^{2q}$ be two handlebodies embedded in $S^{2q+1}$ with isomorphic Seifert forms. 
 
1) The first step of the proof  consists in showing that they have handle presentations   $W^{2q}$ and $W'^{2q}$ with the same Seifert matrix $A$. This is implied by  C.T.C. Wall \cite{wall62} (see also  corollary13.1.1 Section 13( Appendix IV):

\begin{proposition}
Let $q \geq 3$ and let $F^{2q}$ be a q-handlebody. Then there is a natural bijection between handle presentations of $F^{2q}$ and basis of $H_q(F^{2q} ; {\bf Z}) = \pi_q (F^{2q})$.
\end{proposition}

2) Let $L$ (resp.  $L'$) be the  attaching $(q-1)$-link of the handle presentation of $W^{2q}$ (resp. $W'^{2q}$ ). Let $\lbrace c_j, 1\leq j\leq r \rbrace$ and  $\lbrace c'_j, 1\leq j\leq r \rbrace$ be the $q$-cycles corresponding to the handle presentations (defined as in the proof of the realization theorem). 
\\ By hypothesis $ \mathcal {L}_{S^{2q-1}} (L_i ; L_j)=I(c_i,c_j)= I(c'_i,c'_j)=\mathcal {L}_{S^{2q-1}} (L'_i ; L'_j)$. The following theorem of Haefliger \cite{haef62-1} in the case of  $(q-1)$-dimensional  links implies that $L$ and $L'$ are isotopic.  Here we need $(q-1)\geq 2$ i.e. $q\geq 3.$

\begin{definition}
A q-link  with $r$ components  in $S^{2q+1}$ is a differential submanifold $L \subset S^{2q+1}$, where $L= \amalg _{1\leq j \leq r} L_j  $,  the components  $L_j$ being orientation preserving diffeomorphic to $S^q$. Note that the components are labelled. 
\end{definition}

Next theorem is due to Andr\'e Haefliger. See Appendix IV, Subsection 13.3.

\begin{theorem}(Haefliger \cite{haef62-1})
Let $q \geq 2$.  Let $L =\amalg _j L_j \subset S^{2q+1}$ and $L' =\amalg _j L'_j \subset S^{2q+1}$ be two q-link with $r$ components. There exists a diffeomorphism $\Phi : S^{2q+1} \rightarrow S^{2q+1}$ such that:
\\
i) $\Phi$ is isotopic to the identity;
\\
ii) $\Phi (L_j) = L'_j $ for $1 \leq j \leq r$ preserving both orientations;
\\
if and only if Condition $\heartsuit$ is satisfied:
\\
$\heartsuit ~~~\mathcal {L}_{S^{2q-1}}  (L_i ; L_j) = \mathcal {L}_{S^{2q-1}}  (L'_i  ;L'_j)$ for $1 \leq i < j \leq r$. 
\end{theorem}

3) The $q$-cycles $c_j$  and $c'_j$ are $q$-spheres embedded in $S^{2q+1}.$ By hypothesis
 $$(-1)^q a_{ij}= \mathcal {L}_{S^{2q-1}} (c_{i} ; i_+(c_{j}))=  \mathcal {L}_{S^{2q-1}} (c'_{i} ; i_+(c'_{j})), \ i\neq j.$$

 We can adapt the   theorem of Haefliger  to obtain   that the skeletons $\sigma (W)$ and $\sigma (W')$ are isotopic in $S^{2q+1}$. Here, we need  $q\geq 2.$

4)  Point 3) implies that we can assume that $ \ \sigma (W) =\sigma (W') \subset S^{2q+1}$.  Then the  $q$-cycles $c_j= c'_j$ are embeddings of $S^q$ in $W^{2q}$ and $W'^{2q}$. It is sufficient to prove that the normal disc bundles of   $c_j$ in $W^{2q}$ and  in $W'^{2q}$ are isotopic in $S^{2q+1}$. As the handle presentations are embedded in $S^{2q+1}$ the bundles are stably trivial. Let $N(c_j)$ be a tubular neighbourhood of $c_j$ in $S^{2q+1}$.  It is  is isomorphic to the trivial   normal  disc bundle of $c_j$ in $S^{2q+1}$. The boundary of $N(c_j)$ is diffeomorphic, via any trivialisation of $N(c_j)$, to $(S^q\times S^q)$.  But  $c_j= c'_j$ and $ \mathcal {L}_{S^{2q+1}} (c_{j} ; i_+(c_{j}))= \mathcal {L}_{S^{2q+1}} (c'_{j} ; i'_+(c'_{j}))$. This  implies that  $i_+(c_j)$ and $ i'_+(c_j)$ are isotopic on $bN(c_j)$. As the normal disc bundle of $c_j$ in $W^{2q}$ (resp. in $W'^{2q}$) is the orthogonal supplement to $i_+(c_j)$  (resp. to $i'_+(c_j)$) in $N(c_j)$, we have obtained the ambient isotopy between $W^{2q}$ and $W'^{2q}.$

{\bf End  of the sketch of proof of  Levine's classification of embeddings of handlebodies}.

{\bf Note.} The classification of simple $2q$-knots is difficult, much more than the classification of simple $(2q-1)$-knots. It was performed by Cherry Kearton (see \cite{kear83}) except for some difficulties with the 2-torsion. It was thoroughly treated by Mickael Farber; see \cite{farb84}.

The abstract structure of parallelisable handlebodies is described   in Section 13 (Appendix IV). The Seifert form provides an easy way to get it when the handlebody is embedded. In particuliar, when  $3\leq q$ the Seifert form determines  the diffeomorphic class of the boundary of a $q-$handle body. When $bW^{2q}$ is a homotopy sphere,  we give  a study of its  differential type  in Section 13.

Let $Q^{A}: H(W^{2q};{\bf Z}) \rightarrow  {\bf Z}$ be the quadratic form associated to the Seifert form i.e. $Q^{A}(x)=\mathcal A (x,x)$ and let $Q^{A}_2$ be its reduction modulo $2$.

\begin{theorem}   Let $F^{2q}$ be a $q$-handelbody embedded in $S^{2q+1}$ and let $\mathcal A $ be its Seifert form,  $I=(-1)^q(\mathcal A +(-1)^q \mathcal A^T) $ its intersection form and $Q^{A}_2$ the associated quadratic form.
\\
1) If $q$ is even and $q>3$, then $I$ classifies  $F^{2q}$ up to diffeomorphism.
\\
2) If $q=1,3,7$, then $I$ classifies $F^{2q}$ up to diffeomorphism. 
\\
3) If $q$ is odd and $q\neq 1,3,7,$ then $I$ and $Q^{A}_2 $ classify $F^{2q}$ up to diffeomorphism.
\end {theorem}

The proof of the theorem results from what is presented in the appendix on handlebodies up to one argument. We need to use the lemma on normal bundles p.512 of  Kervaire-Vasquez \cite{keva67} in order to know that  $Q^{A} $  and   $Q^{A}_2 $  coincide  with the forms $Q_W $ defined in Section 14, see 14.2.1 (Appendix V).

\vskip.5in

\section{Knot cobordism}

\subsection{Definitions}

Let $K^{n} \subset S^{n+2}$ be a $n$-dimension knot. We denote by  $-K^{n} \subset  -S^{n+2}$ the knot obtained after taking  the opposite orientation on $K^n$ and on $S^{n+2}$. Let  $I$ be the interval  $ \lbrack 0 , +1 \rbrack$.

\begin{definition}  Two knots  $K_0^{n} \subset  S^{n+2}$ and   $K_1^{n} \subset  S^{n+2}$ are {\bf cobordant} if there exists a smooth ($n+1$)-dimension oriented manifold $C^{n+1}$ imbedded in $I\times S^{n+2}$ such that:
\\
i) The boundary of  $C^{n+1}$ is equal to $-K_0^n\cap K_1^n$ and $C^{n+1}$ meets orthogonally the boundary of $I\times S^{n+2}$.
\\
ii) $C^{n+1}$ is diffeomorphic to $ I\times \Sigma^n  $, where $\Sigma^n$ is the homotopy sphere diffeomorphic to $K_0$ and $K_1$. 
\\
By definition, a  {\bf null-cobordant knot} is a knot cobordant to the trivial knot.
\end{definition} 

{\bf{Remark.}} To be cobordant is an equivalence relation on the set of $n-$knots.  Moreover, the connected sum operation induces an  abelian group structure on  the set of  cobordism classes of  $n-$knots.  The  trivial  knot represents  the zero element  and the class of   $-K^{n} \subset  -S^{n+2}$ is the inverse of the class of  $K^{n} \subset  S^{n+2}$.

\begin{definition}The standard sphere $S^{n+2}$ is the boundary of the standard ball $B^{n+3}$.
A knot  $K^{n} \subset  S^{n+2}$ is  {\bf slice} if  $K^n$ bounds,  in $B^{n+3}$,  an oriented  manifold $\Delta ^{n+1} $,   diffeomorphic to the standard ball $B^{n+1}$,  and $\Delta ^{n+1}$  meets  orthogonally $S^{n+2}$. 
\end{definition}

It is obvious that a knot is slice if and only if it is null-cobordant.

\subsection{The even dimensional case}

In his paper \cite{kerv65}  Michel Kervaire considers the group   $C_{2q}$ of the cobordism classes of  knots $K^{2q} \subset S^{2q+2}$, when $K^{2q}$ is diffeomorphic to the standard sphere $S^{2q}$. 

{\bf{Remark.}}  If $2q \neq 4$, an embedded homotopy sphere $K^{2q} \subset S^{2q+2}$  is diffeomeorphic to the standard $2q$-sphere. When $2q=4$, Michel  Kervaire needs  that $K^4$ is diffeomorphic to $S^4$. We will see why in the sketch of proof below.

The  main contribution  of  Kervaire's article, in  knot-cobordism theory, is  the following theorem (Theorem III, p.262 in \cite{kerv65})

\begin{theorem}
For all $q \geq 1$,  the group $C_{2q}$ is trivial. 
\end{theorem}

It seems that  Michel Kervaire judged the proof of   $C_{2q}=\lbrace 0 \rbrace $  easy.  Indeed,  it is the case if one is willing to admit the twenty pages or so of chapters 5, 6 and 7 of Kervaire-Milnor.  These pages are needed  to prove the following main statement:

{\bf{Main statement:}} The knot $K^{2q}$ bounds in the ball $B^{2q+3}$ 
 a contractible manifold $V^{2q+1}$  which meets orthogonally $S^{2q+2}$.

The main statement implies the theorem. There are two cases:

1)  If $2q \geq 5$,  $V^{2q+1}$ is a ball (Smale).

 2) When $2q=4$, Michel  Kervaire uses that $K^4$ is diffeomorphic to $S^4$. In this case   the manifold $V^5$ is a priori  "only" contractible. But, since  $K^4=bV^5$ is diffeomorphic to  the standard sphere, $V^5$ is in fact a ball. Indeed,  we can attach a 5-ball on $bV^5$. We get a homotopy 5-sphere $\Sigma^5$. Fortunately, Kervaire-Milnor prove that $\Theta^5 = 0$. Hence $\Sigma^5$ is diffeomorphic to $S^5$. By Cerf isotopy theorem the 5-ball attached to $bV^5$ is isotopic to a hemisphere. Then,   $V^5$ is isotopic to the other hemisphere and is  diffeomorphic to $B^5$. This  implies that  $C_{4}={0}$.

 Let us remark that,  if we admit homotopy spheres as 4-knots, in other words possibly exotic differential structures on $S^4$, the cobordism group of such $4-$knots could be non-trivial  (if such exotic structures exist).
 
 3) When $2q=2$, more surgery is needed  on the homotopy ball $V^3$ (Michel Kervaire did it in \cite{kerv65} , Lemme III.7, p.265),  to obtain a standard $3-$ball.  Kervaire's proof does not require a proof of the  Poincare Conjecture.
 
 The proof of the main statement has two steps.
 
 I) Let $F^{2q+1}$ be a Seifert hypersurface for $K^{2q}$. It is a framed differential manifold which has  the standard  sphere as boundary. Now, we can use the following theorem:

\begin{theorem}(Kervaire-Milnor \cite{kemi63})
There exists a framed manifold $W^{2q+2}$ obtained from $F^{2q+1} \times I$ by attaching handles  of index $j \leq q+1$ on $F \times \lbrace 1 \rbrace$ and such that
$$bW = (F \times \lbrace 0 \rbrace) \cup ( bF  \times I ) \cup V$$
with $V^{2q+1}$ contractible and $(F \times \lbrace 0 \rbrace )\cap (bF \times I )= bF\times \lbrace 0 \rbrace$ and also $ ( bF  \times I ) \cap V=bF\times \lbrace 1 \rbrace = bV$.
\end{theorem}

II) The second step of the proof  consists in embedding  $W^{2q+2}$ in the  ball $B^{2q+3}$ in such a way that the image of $F^{2q+1} \times \lbrace 0 \rbrace$ is equal to the already  chosen  Seifert hypersurface of  $K^{2q}$. The image of $( bF  \times I ) \cup V^{2q+1}$  provides the contractible manifold which has  $K^{2q}$ as boundary . 
 
Michel  Kervaire provided two proofs of the existence of the embedding of $W^{2q+2}$. One is given in  \cite{kerv65},  the other is  in his Amsterdam paper \cite{kerv70}.

{\bf The Amsterdam paper proof.} We begin by invoking Hirsch's theory of immersions (in fact a relative version of it). Since $W$ is framed with non-empty boundary it is parallelisable. Hence there exists an immersion of $W$ in the ball which coincides with the identification of $F^{2q+1} \times \lbrace 0 \rbrace$ with a Seifert hypersurface. Then we use a trick (much used by Moe Hirsch and others) which consists in jiggling the images of the cores of the handles to get them embedded and disjoint. The jiggling is possible since the highest dimension of the cores is $(q+1)$ and since $2(q+1) < 2q+3$ (general position). Since immersions constitute an open set among the differentiable maps, we get a new  immersion which is an embedding on a neighbourhood of a spine of $W$ (rel. $F^{2q+1} \times \lbrace 0 \rbrace)$ and  which is  therefore regularly homotopic to an embedding. 

{\bf Kervaire's proof from his Paris paper.} People who prefer algebra to topology hate the jiggling arguments. Here is a sketch of a more reasonable proof. Kervaire-Milnor's proof of the existence of the manifold $W$ results from a careful succession of framed surgeries. Kervaire shows that these surgeries can be embedded in the ball $B^{2q+3}$. Two ingredients are important. 
\\
1) In the Kervaire-Milnor approach obstructions to perform the surgeries take place in homotopy groups which are ``stable". It happens that in order to perform embedded surgeries the obstruction we meet are ``already" stable and hence vanish, since they vanish in the non-embedded case. 
\\
2) To embed the cores of the handles we need  a general position argument, which works as above since again $2(q+1) < 2q+3$. 

{\bf Remark.} The general position argument fails in the case of odd dimension knots since we have cores of handles  of dimension $q+1$ in the  ball $B^{2q+2}$.

\subsection{The odd dimensional case}

  Here we denote by  $C_{2q-1}$, the group   of the cobordism classes of  knots $K^{2q-1} \subset S^{2q+1}$, when $K^{2q-1}$  a homotopy sphere. 
  
 In contrast to the even dimensional  case, the odd dimensional  cobordism groups are non-trivial. Kervaire, who considers only the cobordism classes of knots diffeomorphic to the standard $S^{2q-1}$,  proved that the cobordism groups $C_{2q-1}$  are not finitely generated. His  ``easy proof",  uses the  Fox-Milnor condition on  the Alexander polynomial. Next step was performed by Levine, who reduced the computation of $C_{2q-1}$ to an algebraic problem.

We  summarize Levine's  theory. From the Seifert form, an obstruction is constructed for a knot to be null-cobordant (which is also valid  when $q=1$). An  equivalent of  the Witt-relation can be defined on the (non-symmetric) Seifert forms. If $q > 1$ the ``Witt group"  of Seifert forms produces  a complete classification of knots in higher dimensions up to  cobordism. Let  us be more explicit.

\begin{definition} An  $\epsilon -$form is a bilinear form $\mathcal{A}: M\times M \rightarrow  {\bf{Z}}$, where $M$ is a free ${\bf Z}-$module of finite rank $r$, such that $\mathcal{I}=\mathcal{A}+\epsilon \mathcal{A}^T$ is an  $\epsilon-$symmetric unimodular form (i.e. the  determinant of $\mathcal{I}$ is  equal to $\pm 1$).  
\\
An  $\epsilon -$form is null-cobordant if $r=2r'$ is even and if there exists  a   submodule $H$  of $M$, of rank $r'$, which is pure (i.e. the quotient $M/H$ is {\bf Z-}torsion free),  and   such that $\mathcal{A} (x,y)=0 \ $ for all  $ (x,y)\in H\times H$.  Such a $H$ is called a {\bf  metaboliser} of $\mathcal{A}$.
\\
 The unimodularity of $\mathcal{I}$ implies that a metaboliser $H$ of $\mathcal{A}$ is equal to its orthogonal for $\mathcal{I}.$
\end{definition}

\begin{definition} 
Two $\epsilon-$forms $\mathcal{A}_1$ and $\mathcal{A}_2$ are cobordant if the orthogonal sum $\mathcal{A}_1  \boxplus  -\mathcal{A}_2 $ is null-cobordant.
\end{definition}

{\bf Remark} To be cobordant is a equivalence relation on the set of the $\epsilon-$forms. The transitivity needs a proof (see Levine \cite{levi69} or Kervaire \cite{kerv70}). The orthogonal sum provides a structure of abelian group on the set of equivalence classes. We denote this abelian group by ${\bf C}_{\epsilon}$.

\begin{theorem}  (Levine) 
If $1\leq q$  the Seifert forms  of    cobordant $(2q-1)-$knots  are   cobordant $(-1)^q- $forms.
\end{theorem}

{\bf Remark} A knot has several Seifert forms since  it has several  Seifert hypersurfaces. In particular, the theorem implies that two Seifert forms of the same knot are cobordant.

\begin{corollary} Let $\epsilon  =(-1)^q.$  The correspondence  which associates to each knot one of its Seifert forms induces a group homomorphism:
$$ \gamma _q  : C_{2q-1} \rightarrow {\bf C}_{\epsilon}.$$

\end{corollary}

One can find many proofs of this theorem. The  first proof was given  by Levine (see also Kervaire's Amsterdam paper  \cite{kerv70}).  It is sufficient to prove the following statement:  null-cobordant knots have null-cobordant Seifert forms. We now give an idea of the proof of the statement.

 Let $F^{2q}$ be a Seifert hypersurface of  a $(2q-1)-$knot  $K^{2q-1} \subset S^{2q+1}$ which bounds in $B^{2q+2}$ an  embedded ball $\Delta ^{2q}$. Let us consider the closed oriented  submanifold  $M^{2q}=F^{2q} \cup \Delta ^{2q}$  in $B^{2q+2}$. As $M^{2q}$ is a  codimension two embedded submanifold of  $B^{2q+2}$, by  an argument similar to the existence of the Seifert hypersurfaces for links Section 11 (Appendix II),  $M^{2q}$ bounds in $B^{2q+2}$  an  oriented sumanifold $W^{2q+1}$. Let $j_q: \ H_q(M^{2q},{\bf Z}) \rightarrow H_q(W^{2q+1},{\bf Z})$ be  the homomorphism induced by the inclusion. Using the homology exact sequence  of the pair $ (W^{2q+1},M^{2q})$ and  Poincare duality, one can  show that  $Ker j_q$ is a metaboliser of the Seifert form  $\mathcal{A}$ associated to $F^{2q}.$
 
 In any odd dimensions, the theorem gives  also a necessary condition for a knot to be null-cobordant which is readable on the Alexander polynomial. It is this argument which enables  Michel Kervaire to obtain the last main theorem of is paper \cite{kerv65}: `` Pour $1\leq q$, le groupe $C_{2q-1}$ n'est pas de type fini". We  explain now  this result with the help of the theorem above.
 
 \begin{definition}  Let   $A$  be  a Seifert matrix of  $K^{2q-1} \subset S^{2q+1}$ and let $ \ P(t )=det ( \ t  A  +(-1)^q A^T )$. The Alexander polynomial  of  $K^{2q-1}$  is the class of $P(t)$,  modulo multiplication by a unit of ${\bf Z} \lbrack t , t^{-1} \rbrack$.
 
 \end{definition}
 
 By the previous theorem, if $K^{2q-1}$ is null-cobordant any of his Seifert matrix $A$ is null-cobordant. Then,  a direct computation of $det ( \ t A  +(-1)^q  A^T )$ implies  a Fox-Milnor condition on  the Alexander polynomial of null-cobordant  knots:

 \begin{corollary}
 Let $1\leq q$.  If $P(t)$ is the Alexander polynomial of a null-cobordant $(2q-1)$-knot, there exist  a polynomial $Q(t)$ in  $  {\bf Z} \lbrack t  \rbrack$
such that $P(t)=Q(t)Q(t^{-1})$ modulo multiplication by a unit of  ${\bf Z} \lbrack t , t^{-1} \rbrack$.

 \end{corollary}
 
 In his paper \cite{levi69},  Levine  considers the homomorphism $\gamma _q$ and shows that:
 \begin{theorem}
 
   1) If $q\neq 2$, $\gamma _q$ is surjective.
   \\ 2) In ${\bf C_{+1}} $, let  $C_0$ be the subgroup of  the classes of $+1$-forms  $\mathcal{A}$ such that $16$ divides the signature of $\mathcal{A} + \mathcal{A}^T$. Then,  $\gamma _2$ is an isomorphism on $C_0$.
   \\  3) If $2<q$, $\gamma _q$ is an isomorphism.
    \end{theorem}
 
 {\bf Remark.} To obtain the  surjectivity in statement 3),  we need to  consider  knots which are homotopy spheres.

 Now, we  give some indications on the proof of the above statements.
 
 When $q\geq 3 $, the realization theorem  of Seifert forms  \cite{kerv65} implies  the surjectivity of $\gamma_q$.  The first step to prove  the   injectivity  consist in the following   theorem  of Levine  \cite{levi69} (also proved by Kervaire in \cite{kerv70} when $q\geq 3$):

\begin{theorem}
Let $q\geq 2$. Every $(2q-1)$-knot is cobordant to a simple knot. 
\end{theorem}

Let $F^{2q}$ be a Seifert hypersurface of a knot $K^{2q-1}$. Kervaire and  Milnor provide (Theorem 6.6 in \cite{kemi63}) the following result:

\begin{theorem}
  Let $q\geq 2$ and let   $F^{2q}$ be a framed differential manifold whose boundary is   a homotopy sphere. Then   there exists a framed manifold $W^{2q+1}$ obtained from $F^{2q} \times I$ by attaching handles  of index $j \leq q$ on $F^{2q}\times \lbrace 1 \rbrace$ and such that
$$bW^{2q+1} = (F^{2q} \times \lbrace 0 \rbrace) \cup (bF^{2q}  \times I) \cup X^{2q}$$
with $X^{2q}$ $(q-1)$-connected. Moreover,  we have:  $ (F^{2q} \times \lbrace 0 \rbrace) \cap (bF^{2q}  \times I) = bF^{2q}\times \lbrace 0 \rbrace$ and  $ bX^{2q}=bF^{2q}\times \lbrace 1 \rbrace .$
\end{theorem}

The manifold $W^{2q}$ can be embedded in the ball $B^{2q+2}$ in a way similar to the even dimensional  case (see Theorem 8.2.2). An ``engulfing argument" produces an embedding of $X^{2q}$ in a $(2q+1)-$sphere $S^{2q+1} \subset B^{2q+2}$. For the use of  engulfing see Levine's paper  (proof of Lemma 4). A simplified use of engulfing is given by Kervaire in his Amsterdam paper middle p.91. Then, $K^{2q-1}$ is cobordant to the boundary of  $X^{2q}$ which is a simple knot.

Next step consists in taking  a knot   $K^{2q-1}$ in the kernel of $\gamma _q$ . From the last result we can assume that $K^{2q-1}$ bounds a $(q-1)$-connected  Seifert surface $F^{2q}$ and that  $F^{2q}$ has a null-cobordant  Seifert form  $\mathcal A$. This  implies that  we can find a basis
 $\lbrace c_j,\ 1\leq j \leq r \rbrace $  of a metaboliser of  $\mathcal A$ which can be completed in a basis $ \lbrace c_j,c^*_j, \ 1\leq j \leq r \rbrace $ of $H_q(F^{2q},{\bf Z}) $  which is a hyperbolic basis for the intersection form $I$ when $(-1)^q=+1$ and symplectic when $(-1)^q=-1.$ We push $F^{2q}$  inside  the ball $B^{2q+2}$ with a small collar around $K^{2q-1}=bF^{2q}$.  Let $F$ be the result of this push. If $q\geq 2$ there is no obstruction to perform an embedded surgery on the $q$-cycles  $\lbrace c_j,\ 1\leq j \leq r \rbrace $ because $\mathcal A$ vanishes on the metaboliser.  As the basis $ \lbrace c_j,c^*_j, \ 1\leq j \leq r \rbrace $ is either hyperbolic or symplectic for $I$ the result of this surgery is diffeomorphic to a $2q$-ball. The knot $K^{2q-1}$ is null-cobordant. This ends the sketch of proof.
 
 The hard work is  now to describe  the  group  ${\bf C}_{\epsilon}$.   In  \cite{levi69}, Levine  gives a complete list of  cobordism invariants of an  $\epsilon$-form. In   \cite{kerv70}  Kervaire obtains partial results on the group ${\bf C}_{\epsilon}$ and in  \cite{stol76}  Neal Stoltzfus describes   it   in details.
 
 Andrew Ranicki has constructed a thorough algebraisation of higher dimensional knot theory, including Seifert surfaces and Blanchfield duality. He applies his concepts to various situations such as simple knots, fibered knots and knot cobordism. See \cite{rani98}, \cite{rani03}.

 \vskip.5in

\section{Singularities of complex hypersurfaces}
In this section we  summarize some results  about the topological type of singularities of complex hypersurfaces. We  only present results or problems  related to Kervaire's work on homotopy spheres and knots in higher dimensions (Seifert forms, knot cobordism,..).

\subsection{The theory of Milnor}

The connection between higher dimensional  homotopy spheres and isolated singularities of complex hypersurfaces was established in Spring 1966. The story is beautifully (and movingly) told by Egbert Brieskorn in \cite{brie00} in pages 30-52. Several mathematicians took part in the events: Egbert Brieskorn, Klaus J\"{a}nich, Friedrich Hirzebruch, John Milnor and John Nash. In June 1966, Egbert Brieskorn proved the following theorem, which is a corollary of his thorough study of the now called Pham-Brieskorn singularities. The proof rests on the work of Fr\'ed\'eric Pham.

\begin{theorem}  Let   $ \  \Sigma ^{2q-1}$ be  a $(2q-1)-$homotopy sphere which bounds a parallelisable manifold. Then there exists $(q+1)$ integers $a_i, \ 0\leq i \leq q,$ such that the link  $ L_{f}\subset S^{2q+1}$ 
associated to  $f(z_0,...,z_q)= \Sigma  _{i=0}^{i=q}  \ z_i^{a_i}$  is a knot diffeomorphic to $ \  \Sigma ^{2q-1}.$ 
\end{theorem}

 We now recall the main features of Milnor's theory contained in  Milnor's famous book \cite{miln68}.

Let $f : ({\bf C}^{q+1} , 0 )\rightarrow ({\bf C} , 0)$ be a germ of  holomorphic  function with an isolated critical  point at the origin in  ${\bf C}^{q+1}$.

 We denote by $B^{2q+2}_r$  the $(2q+2)$-ball, with radius $r>0$ centered at the origin of ${\bf C}^{q+1}$ and by $S^{2q+1}_r$  the boundary of $B^{2q+2}_r$.  Let  us take the following notations:   
               $F_0 = B^{2q+2}_{\epsilon} \cap f^{-1}(0)$  and $L_{f} = S^{2q+1}_{\epsilon} \cap f^{-1}(0).$

 I)  For a  sufficiently small $\epsilon$,  $ f^{-1}(0)$ meets  $S^{2q+1}_{\epsilon '}$  transversally for all $\epsilon ',\ 0<\epsilon ' \leq \epsilon  ,$ and the homeomorphism  class of the pair $(F_0 \subset  B^{2q+2}_{\epsilon} )$ does not depend on  such  a sufficiently small $\epsilon.$ By definition,  it is the {\bf topological type } of $f.$ 
 
 II) The orientation preserving  diffeomorphism  class of the pair  $(L_{f} \subset S^{2q+1}_{\epsilon} )$ does not depend on  such a sufficiently small $\epsilon.$  By definition,  it is the {\bf link  } of $f.$ By Milnor's conic structure theorem  (\cite{miln68}, p.18), it determines the topological type of $f.$
 
 III) There exist a sufficiently small $\eta, \ 0<\eta <<\epsilon ,$ such that $ f^{-1}(t)$ meets  $S^{2q+1}_{\epsilon '}$  transversally for all $t\in B_{\eta}^2.$ For   such  a  sufficiently small $\eta ,\ $ $N(L_f)= f^{-1}(B_{\eta}^2) \cap S^{2q+1}_{\epsilon} $  is a tubular neighbourhood   of $L_f$ and the restriction  $f_{\vert  N(L_f)}$ of $f$ on $N(L_F)$ is a  (trivial) fibration over $B_{\eta}^2$. 
 
 IV) The manifold $f^{-1}(t) \cap B_{\epsilon}^{2q+2}$ for some $t \neq 0$ is the {\bf Milnor fiber} $F_f$ of $f$. Its boundary is diffeomorphic to the link $L_f$ since the critical point is isolated. The Milnor fiber is oriented by its complex structure and $L_f$ is oriented as its boundary.

  Open books were introduced (without the name) by John Milnor in \cite{miln68}.  The reader will find more details on open books in Section 12 (Appendix III). We use the following definition:

\begin{definition}

 Let $ L$  be  a closed oriented $(2q-1)-$dimension submanifold of $S^{2q+1}.$   An {\bf open book decomposition} of $S^{2q+1}$ with {\bf binding}  $L$ is given by a differentiable fibration $\psi$ over the circle: $$\psi :\ E(L) \rightarrow {\bf S}^1,$$  and a (trivial)  fibration $\phi$:
  $$\phi :\ N(L) \rightarrow {\bf B}^2,$$ 
  such that the restrictions of  $\psi$  and $\phi$ on $bN(L)=bE(L)$  coincide.
  For  $z\in S^1$, a  fiber $F=\psi ^{-1}(z)$ of $\psi $ is a {\bf page } of the open book.

   Let  $h: F\rightarrow F $ be a diffeomorphism  of $F$. The  {\bf mapping torus} $T(F; h)$ of $F$ by $h$ is the quotient of $F \times \lbrack 0 , 1 \rbrack$ by the glueing $(x , 1) \sim (h(x) , 0)$ for $x \in F$. 
 A diffeomorphism   $h$ of $F$  is a   {\bf monodromy} of $\psi$  if  the   projection on the second factor,  $(x,t)\mapsto t,$ induces a fibration  $\pi : T(F; h) \rightarrow S^1$  isomorphic to $\psi.$

\end{definition}

{\bf Remark 1} 
The open book decomposition implies that the boundary $bF$ of the page  defined above is isotopic to $L$ in $N(L),$  and that there exists  a monodromy of $\psi $ which is the identity on  $bF$.

\begin{theorem} (Milnor fibration theorem) Let $E(L_f) = S^{2q+1}_{\epsilon} \setminus \mathring {N}(L_f)$ be the exterior of $L_f.$ The restriction of $f\over \Vert  f \Vert$ on $E(L_f)$ provides an open book decomposition of $S^{2q+1}$ with binding $L_f$ and page $F_f$.
 \end {theorem}
  
{\bf Proof.} In \cite{miln68}, J. Milnor shows that the restriction $\psi $ of  $f\over \Vert  f \Vert$ on  $E(L_f) $ is  a differentiable fibration.  We know, (see  Point III above),  that $ \phi={1\over \eta }f_{\vert  N(L_f)}$ is a (trivial) fibration  over $B^2.$ As $bN(L_f)=f^{-1}(S_{\eta}^1) \cap S^{2q+1}_{\epsilon} $,   $\psi $ restricted on $bE(L_f)=bN(L_f)$ is equal to $ \phi .$
\\
  {\bf End of proof.}

{\bf Remark 2}    When $f$ has a non-isolated singular locus at the origin,    Milnor's proof   that  $f\over \Vert  f \Vert$  restricted on $E(L_f)$ is a fibration  still works. But, in general,  this  fibration does not provide  an open book decomposition.  In many cases $bF_f$ is not  homeomorphic to $L_F.$

\begin{theorem}  
If $f$ has an isolated critical  point at the origin its Milnor fiber $F_f$ is a $q-$handelbody.
 \end {theorem} 
 
 {\bf Proof.}
 Without  assuming  that $f$ has an isolated critical point, Milnor proves in Chapter 5, in \cite{miln68}, that:
\\
 i)  $F_f$ has the homotopy type of a finite C.W. complex of dimension $q.$
 \\
  ii)  $L_f$ is $(q-2)-$connected.
  \\
   If $f $ has an isolated critical point, Milnor proves in Chapter 6, using the existence of an open book decomposition, that:
   \\
    iii) $F_f$ is $(q-1)-$connected and hence has the homotopy type of a bouquet of $q-$spheres.
    \\
     iv) As  $L_f=bF_f$, the boundary $bF_f$ is $(q-2)-$connected.
 \\
 Therefore, if $q\geq 3$, from Smale's Recognition Theorem for handlebodies, Milnor deduces that the fiber $F_f$ is a $q-$handlebody and he  conjectures that it is still true when $q=2.$ Notice that Theorem 1.2 is trivially true when $q=1.$

In \cite{lepe79}, D.T. L\^e and B. Perron prove that $F_f$ is  a $q-$handelbody, for all $q\geq 1.$ Their   proof  based on  L\^e's carrousel construction is independent of  Milnor's   technique  and they don't use Smale's results. For $q=2$, it  is a positive answer to  this  Milnor's conjecture.
 \\
{\bf End of proof.}

\begin{definition}
A $(2q-1)$ {\bf simple fibered   link}  $ L\subset S^{2q+1}$ is a closed oriented $(2q-1)-$dimensional submanifold of $S^{2q+1}$  which is the binding of an open book decomposition,   the     page being   a  $q-$handelbody. A simple fibered link such that $L$ is a homotopy sphere is a {\bf simple fibered knot} (also called  {\bf simple fibered spherical link}).
\end{definition}

\begin{corollary}
If $L_f$ is the link associated to a germ $f$ which has an isolated critical  point,  $L_f$ is a simple fibered link.
 \end {corollary}

\subsection{Algebraic links and Seifert forms}

Let  $L\subset S^{2q+1}$ be a simple fibered link.  By definition $L$ is the binding of an open book decomposition which has  a  $q-$handlebody   $F$ as page.  Then $H_q (F ; {\bf Z})$ is  a free ${\bf Z}$-module of finite rank. We follow the notations  and definitions of Subsection 6.5. for the  page $F.$ We have   the maps: 

  $$i_+:F\rightarrow E(F) ,   (resp. \ i_-:F\rightarrow  E(F)),$$

 Let $(x',y')$ be a pair of $q-$cycles representatives of  $ (x,y)\in H_q (F ; {\bf Z})\times H_q (F ; {\bf Z}).$ The Seifert form $\mathcal{A}_{F}$ associated to $F$ is defined by:
 $$\mathcal {A}_{F} (x , y) = \mathcal {L}_{S^{2q+1}} (x' , i_+ (y'))$$ 
Let $\mathcal {I} :   H_q (F ; {\bf Z})\times H_q (F ; {\bf Z}) \rightarrow {\bf Z}$ be the intersection form on $H_q (F ; {\bf Z})$. Let $\mathcal {A}_{F}^T$ be the transpose of $\mathcal {A}_{F}.$ Let us recall that, with our conventions of signs, we have:

$$(-1)^q \mathcal {I}  = \mathcal {A} + (-1)^q \mathcal {A}^T$$

\newpage

\begin{proposition}
1) The Seifert form $\mathcal{A}_{F}$, associated to a page $F$ of a simple fibered link,  is unimodular.
\\
2)  The fibration $\psi :\ E(L) \rightarrow {\bf S}^1,$ given by the open book decomposition,   admits   $ h= i_{-1}^{-1}\circ i_+$  as monodromy.
\\
3) Let  $A$ be the matrix of $\mathcal{A}_{F}$ in   a ${\bf Z }$-basis  $\mathcal{B}$ of  $H_q (F ; {\bf Z}).$ Then, $((-1)^{(q+1)}(A^T)^{-1}A)$ is a matrix of $h_q$ in the basis $\mathcal{B},$ where $h_q$ denotes the homomorphism induced by $h$ on $H_q (F ; {\bf Z}).$
\end{proposition}

{\bf Proof.}
By definition a bilinear form, defined on a free module of finite rank, is unimodular if it has an invertible matrix.
 As $F$ is a fiber of a fibration $\psi :\ E(L) \rightarrow {\bf S}^1,$ the map $i_+: F\rightarrow E(F)$ will induce an isomorphism, also written $i_+$, from   $H_q (F ; {\bf Z})$ to $H_q (E(F) ; {\bf Z}).$ So, the determinant of a matrix of $i_+$ is $\pm 1.$ Proposition 6.5.1 implies that the determinant of a matrix of $\mathcal{A}_{F}$ is equal, up to sign, to the determinant of a matrix  of $i_+.$
 Point 2) is obvious if we define $i_+$ as half a turn in the positive direction and $i_{-}$ in the negative direction. By Proposition 6.5.1 Point 2) implies  Point 3).
{\bf End of proof.}

Let $U(B)$ be the set of  equivalence classes,   modulo isomorphism,  of  unimodular bilinear forms defined  on free ${\bf Z}-$modules of finite rank. One can find the following theorem in \cite{durf74}.

\begin{theorem}
If $q\geq 3$,  to associate  the Seifert  form  of the  page of an  open book decomposition of a simple fibered link induces a bijective map $\sigma $ between the isotopy classes of  simple fibered links in $S^{2q+1}$ and $U(B).$
\end{theorem}

{\bf Sketch of  proof.} 
If $q\neq 2$,  Browder's Lemma 2  (see Appendix  III)  implies that the Seifert forms associated to isotopic  simple fibered links are isomorphic and hence $\sigma$ is well defined. If $q \geq 3$  Levine's embedding classification  theorem (see 7.3.1)  implies the injectivity of  $\sigma$. Kervaire's realization theorem (see 7.2.1) and the h-cobordism theorem  imply  the surjectivity of $\sigma .$ 
\\
{\bf End of  sketch of  proof}

{\bf Remark 3} If $q=1$ the isomorphism class of the Seifert form associated to a page of the open book decomposition is only an invariant of the link. If $q=2,$ it is an invariant of the open book decomposition.

\begin{definition}
A $(2q-1)-$simple fibered   link,  $ L\subset S^{2q+1},$  is  an  {\bf algebraic link} if   $ L\subset S^{2q+1}$   is isotopic to   a link $ L_{f}  \subset S_{\epsilon}^{2q+1}$  associated to a germ  $f : ({\bf C}^{q+1} , 0 )\rightarrow ({\bf C} , 0)$  of  holomorphic  function with an isolated critical  point at the origin in  ${\bf C}^{q+1}$.

\end{definition}

{\bf Remark  4}  Algebraic links  are a very special  kind  of simple fibered links. The question of the characterization  of  algebraic links among  simple fibered links is an open problem.  This problem is presented in details  in \cite{mich83}. An important  necessary condition to be an algebraic link is  given by the monodromy theorem: ``The monodromy $h_q$ is a quasi-unipotent endomorphism". This necessary condition is far from being sufficient.

{\bf Remark  5} For $q\geq 3, $ Theorem 1.3 implies that  $(2q-1)-$algebraic links are classified, up to isotopy,  by the isomorphism class of the Seifert form associated to their Milnor fiber.  But this result is not true if $q \leq 2$.  In \cite{bomi94}, the authors construct families of pairs  of germs  $f$ and $g$  defined on  $({\bf C}^{2} , 0),$ with isomorphic Seifert forms and non-isotopic links in $S^3$.  By Sakamoto  \cite{saka74} the germs $f_1=f(x,y)+z^2$ and $g_1=g(x,y)+z^2$ obtained by suspension  have also isomorphic Seifert forms. In \cite{arta91}, Enrique Artal shows that  the 3-manifolds  $L_{f_1}$ and $L_{g_1}$  are not homeomorphic. Hence $f_1$ and $g_1$ do not have the same topological type.

\subsection{Cobordism of  Algebraic links}
 
 In \cite{blmi97},  a equivalence relation, called algebraic cobordism, is defined between unimodular bilinear forms with integral coefficients. As a generalization of Kervaire-Levine's theory of knot cobordism, the following theorem is obtained:

\begin{theorem} For all $q>0,$ the Seifert forms associated to a page of  two cobordant $(2q-1)-$dimensional  algebraic links are algebraically cobordant. If $q \geq 3$ this necessary condition to be cobordant is sufficient.
\end{theorem}

 Of course this theorem can be applied to algebraic links. But, as  algebraic links are "rare", their classification up to cobordism is particular. Here are some typical results:

  I)  In \cite {ledu72}, D.T. L\^e shows that ``Cobordant algebraic links in $S^3$ are isotopic."

 II) ``A null-cobordant  $(2q-1)-$algebraic knot,  ($q>0$), is always trivial" (i.e. a null-cobordant algebraic knot  is never  associated to a non-singular germ). This result is proved in  \cite{mich83} with the help of the following  result of A`Campo \cite{acam73} (see also \cite{ledu75}): an algebraic link  has always a monodromy diffeomorphism  without  fixed point. 
 
 III) In  higher dimensions, cobordism of algebraic links does not anymore imply isotopy.
  Explicit examples are constructed in \cite{bomi93}.

\subsection{Examples}

There are methods of computation of Seifert forms associated to algebraic links in $S^3$ (for example in \cite{bomi94}). We  state  Sakamoto's formula which computes,  by induction,   Seifert forms associated to some higher dimension algebraic links. As the determination of signs  is delicate, we  give  explicit matrices of the Seifert form, the monodromy and the intersection  form associated to a ``cusp" singular point. 

We use the following notations:

Let $u = (u_0 , \dots , u_{n}) \in {\bf C}^{n+1}$ ~, ~$v = (v_0 , \dots , v_{m}) \in {\bf C}^{m+1}$, and 
\\ 
$w = (u_0 , \dots , u_{n} , v_0 , \dots , v_{m}) \in {\bf C}^{n+1} \times {\bf C}^{m+1} = {\bf C}^{n+m+2}$.
\\
Let $f : ({\bf C}^{n+1} , 0)  \rightarrow  ({\bf C} , 0) $ and $g : ({\bf C}^{m+1} , 0)  \rightarrow ({\bf C} , 0)$ be germs of holomorphic functions with an isolated singularity at the origin.
Define $h : ({\bf C}^{n+m+2} , 0) \rightarrow ({\bf C} , 0),$ by $h(w) = f(u) + g(v)$. Clearly $h$ has also an isolated singularity at the origin.
\\
Let us denote by $\mathcal {A}_f$ ~, ~$\mathcal {A}_g$ ~,~$\mathcal {A}_h$ the respective Seifert forms, associated to the Milnor fibres $F_f$ ~,~$F_g$ ~,~$F_h$. Remember that the complex structure provides orientations (we have no choice here). 

 Sakamoto's formula  \cite{saka74} is:

$$\mathcal {A}_h = (-1)^{(n+1)(m+1)} \mathcal {A}_f \otimes \mathcal {A}_g$$

i) Let $f_n(u_0,u_1,...u_n)=u_0^2+u_1^2+...+u_n^2$ be the germ  for  the ordinary quadratic singularity.
Then  the Seifert matrix for $f_n$ is equal to $(1)$ when $n\equiv 0\ or \ 3 \ (mod4)$ and is equal to $(-1)$ when $n\equiv 1\ or \ 2 \ (mod4).$ By  Point 3) of Proposition 1.1, its  monodromy  is equal to $(-1)^{n+1}.$

ii)Let $g_0 : ({\bf C} , 0) \rightarrow ({\bf C} , 0) $ be the germ define by:
$g_0(u_0)=v_0^3$.  Let $F_0$ be the Milnor fiber of $g_0$. There exist a basis of $H_0 (F_0 ; {\bf Z})$
such that the Seifert matrix $A_0$ and the monodromy matrix  $h_0$ of $g_0$ in this basis are:

\def\R{{\mathbb R}}
\def\N{{\mathbb N}}
\def\Z{{\mathbb Z}}
\def\Q{{\mathbb Q}}
\def\C{{\mathbb C}}
\def\K{{\mathbb K}}

$$A_0= \left( \begin{array}{cc}
1 & 0\\
-1& 1
\end{array}\right),
  \qquad
h_0=\left( \begin{array}{cc}
0& -1\\
1& -1
\end{array}\right). $$
 \qquad

 iii) Let $n\geq 1.$ We consider the  singularity given by  $g_n(u_0,u_1,...,u_n)=u_0^2+u_1^2+...+u_{n-1}^2+u_n^3$. Let $F_n$ be its Milnor fiber. Using  Sakamoto's formula and the given above formulas,   the Seifert matrix $A_n$, the matrix  $h_n$ of the  monodromy and $I_n$ of the intersection form in a chosen basis  $H_n (F_n ; {\bf Z})$ can be easily computed. For example:

\def\R{{\mathbb R}}
\def\N{{\mathbb N}}
\def\Z{{\mathbb Z}}
\def\Q{{\mathbb Q}}
\def\C{{\mathbb C}}
\def\K{{\mathbb K}}

$$A_1= \left( \begin{array}{cc}
-1 & 0\\
1& -1
\end{array}\right),
  \qquad
h_1=\left( \begin{array}{cc}
0 & 1\\
-1& 1
\end{array}\right), 
  \qquad
I_1=\left( \begin{array}{cc}
0 & 1\\
-1& 0
\end{array}\right),$$

\def\R{{\mathbb R}}
\def\N{{\mathbb N}}
\def\Z{{\mathbb Z}}
\def\Q{{\mathbb Q}}
\def\C{{\mathbb C}}
\def\K{{\mathbb K}}

$$A_2= \left( \begin{array}{cc}
-1 & 0\\
1& -1
\end{array}\right),
  \qquad
h_2=\left( \begin{array}{cc}
0 & -1\\
1& -1
\end{array}\right), 
  \qquad
I_2=\left( \begin{array}{cc}
-2 & 1\\
1& -2
\end{array}\right) .$$

\vskip.3in

{\bf Comments.} 1) Egbert Brieskorn recalls in \cite{brie00} p.46 that Friedrich Hirzebruch revealed already in 1965 the connection between higher dimension homotopy spheres and singularities. He proved that the link of the singularity $f : {\bf C}^6 , 0 \rightarrow {\bf C} , 0$ given by $f(z_1 , ... , z_6) = z_1^2 + \cdots + z_5^2 + z_6^3$ is Kervaire's homotopy sphere of dimension 9.
\\
2) If we follow Egbert Brieskorn, the germ $g_n(u_0,u_1,...,u_n)=u_0^2+u_1^2+...+u_{n-1}^2+u_n^3$ should be called the higher dimension $A_2$ singularity.

\vskip.5in

\section{Appendix I: Linking numbers and signs }

In the literature there is an abundance of definitions and formulae about Seifert forms and matrices. Apparently they look more or less the same, but often they differ in details, in particular in signs. However, in applications signs may be important and not arbitrary. This is especially the case in complex geometry, in particular in singularity theory. Therefore, in the hope to avoid ambiguity, we make explicit the rules we follow. Basically there are two of them: the boundary of an oriented manifold and linking numbers.

\subsection {The boundary of an oriented manifold}

{\bf RULE 1: The oriented boundary of an oriented manifold.} Let $M^m$ be a compact connected manifold of dimension $m$, with non-empty boundary $\partial M$. An orientation of $M$ produces (in fact is) a well defined generator $\lbrack M , \partial M \rbrack \in H_m(M , \partial M ; {\bf Z})$. By definition, the orientation induced by $\lbrack M , \partial M \rbrack$ on $\partial M$ is the image $\lbrack \partial M \rbrack$ of this generator by the boundary homomorphism

$$\partial : H_m (M , \partial M ; {\bf Z}) \rightarrow H_{m-1} (\partial M ; {\bf Z})$$

In the simplicial case (which is adequate for the structures we consider) the boundary homomorphism is induced by the (classical) boundary homomorphism in simplicial chain complexes defined by 

$$\partial (A_0 , A_1 , \cdots , A_j , \cdots , A_q) = \sum_{j=0}^{j=q} (-1)^j (A_0 , A_1 , \cdots , \widehat {A}_j , \cdots , A_q)$$

See for instance Eilenberg-Steenrod \cite{eist52}  p. 88.

If $M^m$ is a differential manifold, an orientation of $M$ may also be defined as a coherent choice of equivalence classes of $m$-frames in each fibre of the tangent bundle. Equivalently and more conveniently, it may also be defined as the choice of an atlas such that the coordinate changes have positive determinant everywhere. In this setting, the orientation induced by $\lbrack M , \partial M \rbrack$ on $\partial M$ is as follows: for $x \in \partial M$ choose a frame $(e_1 , \cdots , e_m)$ such that $e_1$ is orthogonal to $\partial M$ at $x$ and is pointing to the exterior of $M$. Then $(e_2 , \cdots , e_m)$ is a frame tangent to $\partial M$ which produces the desired orientation at $x$. It is important that this rule involving frames is compatible with the rule which involves the boundary operator in homology (in fact it can be deduced from it). These rules are also compatible with Stokes formula. See  Conlon \cite{conl93} p. 230.

{\bf Consequence 1.} Let $U^u$ and $V^v$ be two compact, connected and oriented manifolds, such that $\partial U \neq \emptyset$ and $\partial V = \emptyset$. If we wish to take their product in the oriented category, we should use the product $U \times V$ rather than $V \times U$. In this way, the oriented boundary $\partial (U \times V)$ is naturally equal to $(\partial U) \times V$. For instance if $U = I = \lbrack -1 , +1 \rbrack$ it is recommendable to use $I \times V$ in that order, although many topologists consider the product in the reverse order.

\subsection {Linking numbers} 

{\bf RULE 2: Linking numbers.} (Lefschetz). Let $\gamma_1$ be a $(k-1)$-cycle and $\gamma_2$ be a $(l-1)$-cycle in  $S^{m-1}$. Suppose that $\gamma_1 \cap \gamma_2 = \emptyset$ and that $m = k + l$. Suppose that $B^m$ is oriented; hence $S^{m-1} = \partial B^m$ is oriented too.  Then the linking number $\mathcal {L}_{S^{m-1}} (\gamma_1 ; \gamma_2)$ of $\gamma_1$ and $\gamma_2$ in $S^{m-1}$ is by definition the intersection number $\mathcal {I}_{B^m} (c_1 , c_2)$ where $c_i$ is a chain in $B^m$ such that $\partial c_i = \gamma_i$ for $i = 1 , 2$.

We recall the usual definition of the intersection number. Without loss of generality we may assume that the intersection points of $c_1$ and $c_2$ are finite in number and that the intersection is transversal at such a point $P$. At $P$ we choose a frame representing the orientation of $c_1$ and a frame for $c_2$. We place side by side these two  frames and compare the orientation of $B^m$ at $P$ thus obtained with the given one. The sum of $\pm 1$ we obtain is the intersection number we are looking for. It follows immediately that 

$$\mathcal {I}_{B^m} (c_1 , c_2) = (-1)^{kl}\mathcal {I}_{B^m} (c_2 , c_1)$$

Hence the same symmetry rule is valid for the linking number:

$$\mathcal {L}_{S^{m-1}} (\gamma_1 ; \gamma_2) = (-1)^{kl}\mathcal {L}_{S^{m-1}} (\gamma_2 ; \gamma_1)$$

A consequence of Rules 1 and 2 is the following.

{\bf Consequence 2.}  Let $\gamma_1$ and $\gamma_2$ be cycles in $S^{m-1}$ as in the statement of Rule 2. Let $\beta_i$ be a chain in $S^{m-1}$ such that $\partial \beta_i = \gamma_i$ for $i = 1 , 2$. Then

$$\mathcal {L}_{S^{m-1}} (\gamma_1 ; \gamma_2) = \mathcal {I}_{S^{m-1}} (\gamma_1 , \beta_2) = (-1)^k \mathcal {I}_{S^{m-1}} (\beta_1 , \gamma_2)$$

\vskip.5in

\section{ Appendix II:  Existence of Seifert  hypersurfaces}

Let $M^{n+2}$ be a closed, oriented, 2-connected, $(n+2)$-dimension differential manifold, and let $L^n$ be a closed, oriented, $n$-dimension submanifold of $M^{n+2}$. (We do not assume that $L^n$ is connected).

\begin{definition}
{\bf A Seifert hypersurface of $L^n$} is a smooth, oriented and connected,   $(n+1)$-dimensional  submanifold $F^{n+1}$ in $M^{n+2}$  which has  $L^n$  as boundary.
\end{definition}

The aim of this section is to prove the following theorem:

\begin{theorem}
A submanifold $L^n \subset M^{n+2}$ has Seifert hypersurfaces.
\end{theorem}

\begin{lemma}
The normal bundle of $L^n$ is trivial.
\end{lemma}

{\bf Proof of the lemma.} The normal bundle of $L^n$ is orientable. Since it has rank $2$, it is enough to prove that it has a nowhere vanishing section. The only obstruction to construct such a section is the Euler class $e \in H^2 (L^n , {\bf Z})$. But  $ H^2(M^{n+2} , {\bf Z})=0$ since $M^{n+2}$ is 2-connected. Hence  Thom's formula for the Euler class implies that $e=0$.

{\bf End of proof of the lemma.}

Let $N(L)$ be a  compact tubular neighbourhood of $L^n$. Let   $\nu : \ N(L)  \rightarrow L^n$ be the projection of the normal disc bundle associated to $N(L)$. Let $bN(L)$ be the boundary of $N(L)$ and $E(L)=M^{n+2} \setminus \mathring{N(L)}$ be the exterior of $L^n.$

\begin{definition} For each  connected component $L_i$ of $L^n$,$1\leq i \leq r$, we choose $x_i\in L_i$.  { \bf A meridian $m_i$} of $L_i$ is a  circle $m_i=bN(L)\cap \nu^{-1}(x_i)$ oriented as the boundary of $\nu^{-1}(x_i).$
Let $\pi : (L^n \times B^2) \rightarrow L^n$ be the projection on the first factor. {\bf  A trivialisation of $\nu $} is a diffeomorphism $\phi : N(L) \rightarrow (L^n\times B^2)$ such that $\nu =\pi \circ \phi .$

\end{definition}

{\bf Remark.} The homology group  $H_n (N(L) , {\bf Z})$ is free with a basis given by the fundamental classes $l_i$ of $L_i$ in $N(L)$. By Alexander duality $ H_1(E(L) , {\bf Z})$ is free with a basis given by the homology classes of the meridians that we will also denote  by $m_i$ .

{\bf Proof of the theorem.}

 We will show the existence of a differentiable  map $\psi :\ E(L) \rightarrow {\bf S}^1$ such that the restriction $\psi'$ of $\psi$ on $bN(L)$ extends to a fibration $\tilde \psi : N (L)\rightarrow B^2$  where  $ L^n =\tilde \psi ^{-1}(0)$. Then, we choose a regular value $z$ of $\psi$ and we denote by $I$ the radius in $B^2$ of extremity $z.$ By construction $ F_1^{n+1}=( \psi ^ {-1} (z) )\cup (\tilde \psi ^{-1}(I))$ has  $L^n$ as boundary.  If $F_1^{n+1}$ is not connected, first we remove its closed components, then  we connect the connected components of $ F_1^{n+1}$ with the help of carefully chosen thin  tubes and we obtain  a connected $ F^{n+1}$.

Differentiable maps $\psi : E(L) \rightarrow  \bf {S^1}$  are classified up to homotopy by $ H^1(E(L) , {\bf Z})$. The universal coefficient theorem  and the above remark imply that $ H^1(E(L) , {\bf Z})$ is free and has a basis given by  the duals $m^*_i , 1\leq i \leq r,$ of the meridians.  We   choose  a  differential  map $\psi : E(L) \rightarrow  \bf {S^1}$,    which has degree $+1$ on each $m_i,1\leq i \leq r$.

On the other hand let $\phi$ be a  trivialisation,   $\phi : N(L) \rightarrow  L^n\times B^2$,  of the normal bundle of $L^n$. 
 Let $\pi ' :L^n \times B^2  \rightarrow B^2$ be the projection on the second factor. 
The restriction $\phi'$ of  $\pi '\circ \phi$ on $bN(L)$ has degree $+1$ on each meridian $m_i,1\leq i \leq r$. We will say that such a map $\phi' : bN(L) \rightarrow  {\bf S^1}$ is {\bf given by the parametrisation} $\phi.$

As $ H^1(M^{n+2} , {\bf Z})=0$ and  $ H^2(M^{n+2} , {\bf Z})=0$, we have:
  $$ H^1(bN(L) , {\bf Z})= i_N( H^1(N(L) , {\bf Z})) \oplus i_E( H^1(E(L) , {\bf Z}))$$
where  $i_N$ and $i_E$  are  respectively  induced by   the inclusions of $bN(L)$ in $N(L)$ and $E(L)$.

So, a differentiable map $\alpha : bN(L) \rightarrow \bf {S^1}$   is  homotopic  to a map given by a parametrisation if and only if $\alpha$ has degree $+1$ on each meridian $m_i,1\leq i \leq r$. Let $\psi'$ be the restriction on $bN(L)$ of the already chosen map $\psi : E(L) \rightarrow  \bf {S^1}$. Hence  $\psi'$ is homotopic to a map $\phi'$ given by a parametrisation $\phi$. If necessary, we perform this homotopy in a small collar around  $bN(L)$ and we can extend $\psi'$ to the  fibration $\pi '\circ \phi : N(L) \rightarrow B^2$. 

{\bf End of the proof of the theorem.}

\begin{proposition} If   $F^{n+1}$ be a Seifert hypersurface  of $L^{n}$,  there exists  a differentiable  map $\psi :\ E(L) \rightarrow {\bf S}^1$ such that $\psi$ is of degree $+1$ on  a meridian $m$ and  $ F=F^{n+1}\cap E(L)=( \psi ^ {-1} (z) ) $ where $z$ is a regular value of $\psi$.
  \end{proposition}

 {\bf Proof}.
   As  $F=F^{n+1}\cap $ is oriented,   $F$  has a trivial normal bundle. Hence, we can choose an embedding $\alpha :I \times F \rightarrow E(K)$ where
   $I = \lbrack -1 , +1 \rbrack$, such that  the image $\alpha (I\times F)=N(F)$  is  an oriented compact tubular  neighbourhood of $F$ in $S^{n+2}$ and  $F= \alpha (0\times F)$. We use the following notation: $F_+=\alpha (+1\times F)$ and $ F_-=\alpha (-1\times F)$. 
   
    Let $E(F)=(E(K)\setminus  \alpha (F\times \rbrack -1 , +1 \lbrack))$  be the exterior of $F$. We can   define $\psi$ as follows:
   \\
   when $ y=\alpha (x,t)\in N(F)$  let  $\psi (y)=e^{i\pi t}$,  when $y\in E(F)$  let $\psi(y)=-1$.
   \\
   {\bf end of proof}

\vskip.5in

\section{ Appendix III: Open book decompositions }

\subsection{Open books }

 Open books were introduced (without the name) by John Milnor in \cite{miln68}. We  use the following definition:

\begin{definition}
Let $M^{m}$ be a closed, oriented, $m$-dimensional  differential manifold, and let $L^{m-2}$ be a closed, oriented, $(m-2)$-dimensional differentiable map submanifold of $M^{m}$. Let $N(L)$ be a compact regular tubular neighbourhood of $L^{m-2}$ and $E(L)=M^{m}\setminus \mathring N(L)$ be the exterior of $L^{m-2}.$
\\
An {\bf open book decomposition} of $M^m$ with {\bf binding}  $L^{m-2}$ is given by a differentiable fibration $\psi$ over the circle: $$\psi :\ E(L) \rightarrow {\bf S}^1,$$ such that the restriction $\psi'$ of $\psi$ on $bN(L)$ extends to a fibration $\tilde \psi : N (L)\rightarrow B^2$  and  $ L^{m-2} =\tilde \psi ^{-1}(0)$.

\end{definition}

A fiber $F$ of $\psi$ is a {\bf page} of the open book decomposition. The fibration $\tilde \psi$ produces an isotopy, in $N(L)$, between $L^{m-2}$ and the boundary of $F$.

{\bf Remarks}

1) When  we have an  {\bf open book decomposition} of $M^m$ with binding $L^{m-2}$ the fibration $\tilde \psi$ defined on $N(L)$ is trivial and  the normal bundle of $L^{m-2}$ is trivial.
\\
Hence there exists a trivialization  of $N(L)$ such that the  restriction of the fibration $\psi :E(L) \rightarrow S^1$ to the boundary $ bE(L)$ coincides with the projection of $L \times S^1$ on the second factor.

2) Let $F$ be a fiber of $\psi $ and let   $h: F \rightarrow F$ be a monodromy of $\psi $. The exterior $E(L)$ is diffeomorphic to the {\bf mapping torus} $T(F ; h)$ of the monodromy $h$. Recall that it is the quotient of $F \times \lbrack 0 , 1 \rbrack$ by the glueing $(x , 1) \sim (h(x) , 0)$ for $x \in F$. Point 1) implies that we can choose $h$ such that its restriction on $bF$ is the identity.

\begin{theorem} (Uniqueness theorem for open book decompositions of $S^{n+2}$)
\\
Let $L^n \subset S^{n+2}$ be a compact, simply connected submanifold in $S^{n+2}$ such that $\pi_1( E(L) )= {\bf Z}$. Suppose that $n \geq 4$. Then two open book decompositions of $S^{n+2}$ with binding $L^n$ are isomorphic.
\end{theorem}

It is not difficult to replace $S^{n+2}$ by more general manifolds.

Since $L^n$ is simply connected, there is  a unique trivialization of its normal bundle. Hence the main point  to prove the uniqueness is to prove the uniqueness of the fibration over the circle.

Explicitly, we are given two fibrations $f~ \mathrm {and}~\widetilde {f} : E(L) \rightarrow S^1$. We wish to prove that they are isomorphic.

The proof proceeds in two steps. 

{\bf First step.} By Browder's Lemma 2 (see \cite{brow67}), there exists a diffeomorphism $\Phi : E(L) \rightarrow E(L)$ such that $\Phi (f^{-1} (1)) = \widetilde {f}^{-1} (1)$. In other words, $\Phi $ sends one fiber of $f$ onto one fiber of $\widetilde {f}$.

{\bf Second step.} By Cerf pseudo-isotopy theorem (see \cite{cerf70}), one can modify $\Phi$ to a diffeomorphism $\Psi$ which is an isomorphism between the fibrations. For details and explicit formulae, see Mitsuyoshi Kato \cite{kato74} p.461.

\subsection{Browder's Lemma 2}

\begin{lemma} (Browder's Lemma 2 in \cite{brow67})
\\
Let $W$ be a closed, connected differential manifold of dimension $\geq 6$. Let $f$ and $\widetilde {f}$ be two fibrations $W \rightarrow S^1$. We suppose that:
\\
1) $f$ and $\widetilde {f}$ are homotopic;
\\
2) the fibres $f^{-1} (1)$ and $\widetilde {f}^{-1} (1)$ of $f$ and resp. $\widetilde {f}$ are simply connected.
\\
Then  there exists a diffeomorphism $h: W \rightarrow W$, pseudo-isotopic to the identity, such that 
\\
$h(f^{-1}(1)) = \widetilde {f}{-1}(1)$.
\end{lemma}

{\bf Remark.} If we consider disjoint  lifts of the fibers $f^{-1} (1)$ and $\widetilde {f}^{-1} (1)$ in the infinite cyclic covering of $W$, the region between them is an h-cobordism. Hence the two fibers are diffeomorphic. But we need a diffeomorphism of $W$ which sends one to the other.

The main ingredient in Browder's  proof of the lemma is Browder-Levine's existence theorem for fibrations over the circle in a relative form and Smale's h-cobordism theorem.

There exists a version of the lemma if $\partial W \neq \emptyset$.  The proof is essentially the same as the one for the empty boundary case. It is this version that we use. Additional conditions are:
\\
3) The boundary $\partial W$ of $W$ is a product $Z \times S^1$;
\\
4) the projections $f$ and $\widetilde {f}$ coincide on the  boundary $\partial W$ and are equal to the projection $Z \times S^1 \rightarrow S^1$ on the second factor.

In our use of the lemma, the assumption 1) is satisfied since both projections represent the same  generator of $H^1 (E(L) ; {\bf Z})$, taking orientations into account. The assumption 2) is satisfied since we require that $\pi_1 E(L) = {\bf Z}$. The assumptions 3) and 4) are satisfied by the open book conditions.

\vskip.5in

\section{ Appendix IV: Handlebodies}

\subsection{Bouquets of spheres ans handlebodies}

Handlebodies are the simplest of manifolds with boundary. They were studied by Smale and Wall in the early sixties. See \cite{smal62} and \cite{wall63-1}, \cite{wall63-2},  \cite{wall65}. We  present here the basic facts  on  parallelisable  handlebodies of type $(q,r)$.

\begin{definition}
 Let  $\varphi_j : (bB_j^q) \times B_j^{q} \rightarrow bB^{2q},\ 1\leq j\leq r,$ be $r$ disjoint embeddings. A handle presentation of type $(q,r)$ is  a manifold $W^{2q}$ obtained from the disjoint union $B^{2q} \amalg_j B_j^q \times B_j^q$   by identifying $x \in (bB_j^q) \times B_j^q $ with $\varphi_j(x) \in bB^{2q}$. 
 \\ A $q$-handlebody  is a manifold diffeomorphic to a handle presentation of type $(q,r).$
\end{definition}

The image of $B_j^q \times B_j^q$ in $W^{2q}$ is called a {\bf handle of index $q$} and the image $C_j$ of  $B_j^q \times \lbrace 0 \rbrace $ is its {\bf core}. The collection of embeddings $\Phi = \lbrace \varphi_j \rbrace_j$ is the {\bf attaching map of the presentation} and $L=\amalg _j  \varphi _j(S_j^{q-1} \times \lbrace 0 \rbrace)$ is its  {\bf attaching  link}. If $q\neq 2 $ each  connected component  $L_j$ of the link $L$ of the handle presentation bounds  an  oriented $q-$ ball  $D _j$ in $B^{2q}$. Let $c_j$ be the $q-$cycle obtained as the union of $D _j$ with  the core $C_j .$  We say that    $ \lbrace c_j, \ 1\leq j \leq k \rbrace $  is an {\bf  adapted } (to the handle presentation ) basis of  $H_q(W^{2q} , \bf Z)$.

Inspired by Smale, we denote by  $\mathcal {H} (q,r)$ the set of diffeomorphism classes of manifolds which admit a presentation as above. One can say equivalently that such manifolds admit a Morse function with one minimum and $r$ critical points of index $q$. This is true without any restriction on the integers $(q,r)$. 

The case $q = 1$ is easily understood, often with arguments different from those presented here. It is omitted from our discussion. When $q \geq 2$ a q-handlebody is orientable. We assume that it is oriented.

Clearly a $q-$handlebody has the homotopy type of a bouquet of $k$ spheres of dimension $q$. The converse is trivially true if $q = 1$ but wrong if $q \geq 2$. Counterexamples are for instance provided by homology spheres of dimension $(2q-1)$, with non-trivial fundamental group, which bound a contractible $2q-$manifold. This contractible manifold cannot be a q-handlebody since it would be a $2q-$ball and hence its boundary would be $S^{2q-1}$. 

Differential homology spheres which bound differential contractible manifolds are plentiful. Indeed, Michel Kervaire proved in \cite{kerv69} the following theorem.

\begin{theorem}
Let $M^n$ be a differential, oriented homology sphere with $n \neq 3$. Then there exists a homotopy sphere $\Sigma^n$ such that the connected sum $M^n \sharp \Sigma^n$ bounds a differential contractible manifold.
\end {theorem}

{\bf Remarks.} 1) Kervaire proves also that $\Sigma^n$ is unique.
\\
2) To take the connected sum with a homotopy sphere amounts to change the differential structure of $M^n$ in the neighborhood of a point. 
\\
3) Rohlin's theorem prevents the theorem to be true in the differential category if $n = 3$. However Michael Freedman proved that every homology $3-$sphere bounds a topological contractible 4-manifold. See \cite{free82}.

\begin{theorem}(Recognition theorem)
Suppose that $q \geq 3$. Let $V$ be a $2q-$dimensional  manifold which has the homotopy type of a bouquet of spheres of dimension $q$. Then $V$ is a $q-$handlebody if and only if $bV$ is simply connected.
\end{theorem}

{\bf Remark.} An easy computation reveals that the reduced homology groups $H_i^{red} (bV ; {\bf Z})$ vanish for $i \leq q-2$. Hence, if $q \geq 3$, the boundary $bV$ is $(q-2)$-connected if and only if it is simply connected.

{\bf Proof of the theorem.} By general position the boundary of a $q-$handlebody is simply connected if $q \geq 3$. Conversely, assume that $bV$ is simply connected and  let us prove that $V$ is a handlebody. The original proof is due to Smale. Here is a slightly different one which makes the simple connectivity of the boundary quite visible.  This proof rests on  a classical engulfing argument and Smale's h-cobordism theorem.

Choose a basis $\lbrace e_j \rbrace$ for $j = 1 , \dots , r $,  of $H_q(V ; {\bf Z}) = \pi_q (V)$. Since $q \geq 3$ and since $\pi_1 (V) = 1$ we can represent each $e_j$ by an embedded sphere $S_j^q \hookrightarrow V$ (by Whitney). By general position we can assume that $S_i^q \cap S_j^q$ is a finite set of points. In each sphere $S_j^q$ we construct an embedded arc  $\gamma_j$ which contains all the intersection points of $S_j^q$ with the other spheres $S_i^q$ for $i \neq j$. The union of the arcs $\gamma_j$ for $j = 1 , \dots , r,$ is a compact graph $\Gamma$ embedded in $V$. Let $C\Gamma$ be the abstract cone with base $\Gamma$. Since $V$ is simply connected the inclusion $\Gamma \subset V$ extends to a map $\eta : C\Gamma \rightarrow V$. Since dim$(V) \geq 6$ we can find a $\eta$ which is a piecewise regular embedding (in the sense of Morris Hirsch) such that $\eta (C\Gamma) \cap S_j^q = \gamma_j$ for $j = 1 , \dots , r$.

Let $N$ be regular neighbourhood of $\eta (C\Gamma)$ in $V$. By the Whitehead-Hirsch theory of regular neighbourhoods in the differential category, $N$ is $2q-$ball $B^{2q}$ in $V$ which meets each sphere $S_j^q$ in a regular neighbourhood of $\gamma_j$ hence in a $q-$ball $D_j^q$. Let $B_j^q = S_j^q \setminus \mathring {D}_j^q$. By construction we have $B_i^q \cap B_j^q = \emptyset$ for $i \neq j$. Hence the $B_j^q$'s are the cores of handles of index $q$ attached to the ball $B^{2q}$. Let $h_j$ be the handle with core $B_j^q$ obtained by thickening the core a little bit. The union $B^{2q} \cup_j h_j$ is a $q-$handlebody $W$ embedded in the interior of $V$.

Consider $X = V\setminus \mathring {W}$. We claim that $X$ is a simply connected h-cobordism. By construction the inclusion $ W \subset V $ induces an isomorphism on homology. Hence by excision $H_* (X$ , $bW ; {\bf Z}) = 0$. By Poincar\'e duality $H_* (X$ , $bV) = 0$. The boundary $bW$ is simply connected, since $W$ is a q-handlebody with $q \geq 3$. Now $V = W \cup X$ with $W \cap X = bW$. The manifolds $V , W$ and $bW$ are simply connected. Hence by van Kampen $X$ is also simply connected. To conclude that $X$ is a simply connected h-cobordism we still need that $bV$ is simply connected. It is here that  the hypothesis is used. By Smale $X$ is a product since dim$X \geq 6$. 

{\bf End of proof of the theorem.}

\begin{corollary}
(Of the proof of the theorem)
If $q\geq 3,$ any basis of $H_q(V ; {\bf Z})$ can be realized by the cores of a handle decomposition.
\end{corollary}

We wish now to characterize the attaching map   $\Phi = \lbrace \varphi_j \rbrace_j$ up to isotopy. We assume that $q \geq 3$. 
\\
Each $\varphi_j ((bB_j^q) \times \lbrace 0 \rbrace)$ is an oriented sphere $K_j^{q-1}$ of dimension $(q-1)$ embedded in the sphere $bB^{2q}$ of dimension $(2q-1)$. This knot is trivial and the isotopy class of the oriented link $\mathcal {K} = \lbrace K_j \rbrace_j$ is completely determined by the linking numbers $\mathcal {L}_{bB^{2q}} (K_i , K_j) \in {\bf Z}$ for $i \neq j$. See a subsection below for more details.
\\
The sphere $K_j^{q-1}$ bounds a differential ball $D_j^q$ in $B^{2q}$. The normal bundle of this ball has a unique trivialisation, which by restriction provides a canonical trivialisation of the normal bundle of $K_j^{2q-1} \subset bB^{2q}$. The attaching map $\varphi_j$ provides another trivialisation of this normal bundle. The comparison between the two trivialisations produces an element $Q_j \in \pi_{q-1} (SO_q)$. Next lemma summarises what we have obtained so far.

\begin{proposition}
Up to isotopy the attaching map $\Phi$ is characterised by the linking coefficients $\mathcal {L}_{bB^{2q}} (K_i , K_j)$ for $i \neq j$ and the set of elements $ \lbrace Q_j \rbrace_j$ with $Q_j \in \pi_{q-1} (SO_q)$.
\end{proposition}

The proof of the proposition follows easily from Haefliger's classification of spherical links presented below in Subsection 13.3.

{\bf Comment.} The homotopy groups $\pi_{q-1} (SO_q)$ are not stable.  Their value is also periodic of period eight as for the stable ones. See  \cite{wall65}. We shall not need this value here.

\subsection{Parallelisable handlebodies}

We now consider parallelisable (equivalently stably parallelisable) q-handlebodies. This is what we really need from the theory of handlebodies since Seifert hypersurfaces are parallelisable. 

Consider a part of the exact homotopy sequence of the locally trivial fibration $SO_{q+1} \rightarrow S^q$ with fibre $SO_q$:
$$\cdot \cdot \cdot \longrightarrow \pi_q(S^q) \stackrel {\partial}{\longrightarrow} \pi_{q-1}(SO_q) \stackrel{\sigma}{\longrightarrow} \pi_{q-1}(SO_{q+1}) = \pi_{q-1}(SO) \longrightarrow \cdot \cdot \cdot$$

{\bf Comments.}
\\
1) the group $\pi_{j-1}(SO_k)$ classifies oriented real vector bundles (with fibre ${\bf R}^k$) over the sphere $S^j$;
\\
2) the generator of $\pi_q(S^q)$ is the identity map $i : S^q \rightarrow S^q$; its image $\partial (i) \in \pi_{q-1}(SO_q)$ represents the tangent bundle of $S^q$ (this is due to Steenrod,  Section 23 of his book \cite{stee51});
\\
3) $\sigma$ is the stabilisation homomorphism.

Therefore $Ker(\sigma) = Im(\partial) \subset \pi_{q-1}(SO_q)$ classifies  stably trivial q-vector bundles over $S^q$. This subgoup takes the following values (most of this is also due to Steenrod, together with Kervaire's and Bott-Milnor's non-parallelisability of spheres of dimension $\neq 1 , 3 , 7$; see Kosinski's book \cite{kosi93} p.231).

\begin{theorem}
i) When $q$ is even, $Ker(\sigma) = Im(\partial)$ is isomorphic to the integers ${\bf Z}$; its elements are distinguished by the Euler number of the vector bundle, which can take any even value (since the Euler characteristic of an even dimension sphere is equal to 2).
\\
ii) When $q$ is odd, there are two possibilities:
\\
if $q = 1 , 3 , 7$ ~~ $Ker(\sigma) = Im(\partial)$ vanishes;
\\
if $q \neq 1 , 3 ,7$ ~~ $Ker(\sigma) = Im(\partial)$ is isomorphic to ${\bf Z} / 2$. The non-trivial element is represented by the tangent bundle to $S^q$. 
\end{theorem}

\newpage

Let $W$ be a parallelisable  $q$-handlebody  with a basis     $ \lbrace c_j, \ 1\leq j \leq k \rbrace $  of $H_q(W , \bf Z)$ which is adapted to a handle presentation. As explained above each $q$-cycle $c_j$ provides an element  $Q_j \in \pi_{q-1} (SO_q)$. Since  $W$ is parallelisable, we have   $ Q_j\in Ker(\sigma) = Im(\partial) \subset \pi_{q-1}(SO_q)$. Hence  $Q_j\in {\bf Z}$ if $q$ is even, $Q_j=0$ if $q=1,3,7$ and $Q_j\in {\bf Z} / 2$ if $q$ is odd and $q\neq 1,3,7.$
\\
 One can prove that  $\mathcal {Q}_W(c_j)=Q_j$ extends to a quadratic form 
\\
$\mathcal {Q}_W: H_q(W ; {\bf Z}) \rightarrow {\bf Z}$ if $q$ is even,
\\
    $\mathcal {Q}_W: H_q(W ; {\bf Z}/2) \rightarrow {\bf Z} / 2$  if  $q$ is  odd and $\neq 1 , 3 , 7$. 

 See \cite{kosi93} p.206-208. Next proposition is an immediate consequence of what has just been said.

\begin{proposition}
1) Suppose that $W$ is parallelisable and let $q$ be even. Then the intersection form $\mathcal {I}_W$ is even and one has the equality $2\mathcal {Q}(x) = \mathcal {I}_W(x,x)$ for all $x \in H_q(W , {\bf Z})$. 
\\
2) Suppose that $W$ is parallelisable and let $q$ be odd and $\neq 1 , 3 , 7$. Then the map $\mathcal {Q}_W: H_q(W ; {\bf Z}/2) \rightarrow {\bf Z} / 2$ is a quadratic form with associated bilinear form the reduction $mod ~2$ of the antisymmetric intersection form $\mathcal {I}_W$. In other words we have in ${\bf Z} / 2$ the equality $\mathcal {Q}_W(x + y) = \mathcal {Q}_W(x) + \mathcal {Q}_W(y) + \mathcal {I}_W(x,y)$. If this equality is satisfied one says that $\mathcal {I}$ is associated to $\mathcal {Q}$.
\end{proposition}

\begin{theorem}
Let $q \geq 3$. There is a natural bijection between oriented diffeomorphism classes of oriented parallelisable q-handlebodies and:
\\
(i) if $q$ is even, isometry classes of integer valued, symmetric, even bilinear forms over free ${\bf Z}$-modules of finite rank;
\\
(ii) if $q$ is equal to $3$  or  $7$,  isometry classes of antisymmetric bilinear forms over ${\bf Z}$-modules of finite rank;
\\
(iii) if $q$ is odd and $\neq 3 ,7$,  free ${\bf Z}$-modules of finite rank equipped with an integer valued antisymmetric bilinear form  together with a  quadratic form  with values in ${\bf Z} / 2$ (associated as  in the proposition above).
\end{theorem}

The theorem is a consequence of what we have done so far and from the fact that the intersection determines the linking coefficients (use the corollary of the recognition theorem). It is a special case of Wall's vast study of handlebodies.

\begin{definition}
{\bf A q-handlebody $W$ is unimodular} if the intersection form $\mathcal {I}_W$ is unimodular; in other words if the discriminant of $\mathcal {I}_W$ is equal to $\pm 1$.
\end{definition}

{\bf Comments.} 1) If $q \geq 3$ unimodularity is equivalent to require that the boundary $bW$ is a homotopy sphere.
\\
2) If $q = 2$, unimodularity is equivalent to $bW$ a homology sphere.

\subsection{$m$-dimensional spherical links in $S^{2m+1}$}

\begin{definition}
A m-dimension spherical r-link in $S^{2m+1}$ is a differential submanifold $L^m \subset S^{2m+1}$, with r connected components $\lbrace L_1^m , \cdots , L_r^m \rbrace$. Each $L_i^m$ is diffeomorphic to the standard sphere $S^m$ and oriented. Note that the components are labelled. 
\end{definition}

\begin{theorem} (Haefliger \cite{haef62-1})
Let $m \geq 2$.  Let $L^m \subset S^{2m+1}$ and $\widehat {L}^m \subset S^{2m+1}$ be two spherical r-links. There exists a diffeomorphism $\Phi : S^{2m+1} \rightarrow S^{2m+1}$ such that:
\\
i) $\Phi$ is isotopic to the identity;
\\
ii) $\Phi (L_i^m) = \widehat {L}_i^m$ for $1 \leq i \leq r$ preserving both orientations;
\\
if and only if Condition $\heartsuit$ is satisfied:
\\
$\heartsuit ~~~Lk (L_i^m ; L_j^m) = Lk (\widehat {L}_i^m ; \widehat {L}_j^m)$ for $1 \leq i < j \leq r$. 
\end{theorem}

{\bf Comments.}
\\
1) Note that the diffeomorphism $\Phi$ preserves all orientations and the labels. Compare with Haefliger's definitions in \cite{haef62-1}.
\\
2) Condition $\heartsuit$ is clearly necessary for the existence of $\Phi$ for all $m \geq 1$. The point is that the condition is sufficient if $m \geq 2$.
\\
3) Let $\lbrace a_{ij} \in {\bf Z} \rbrace$ be a set of integers with $1 \leq i < j \leq r$.  Then it is easy to prove that there exists a spherical r-link $L^m \subset S^{2m+1}$ such that $a_{ij} = Lk (L_i^m ; L_j^m)$. The $(-1)^{m+1}$-symmetry of linking coefficients enables us  to state the following corollary.

\begin{corollary}(Classification)
\\
Fix $m \geq 2$. There exists a natural bijection between
\\
i) the set of isotopy classes of  spherical r-links in $S^{2m+1}$
\\
and 
\\
ii) the set of  $(-1)^{m+1}$-symmetric $(r \times r)$ matrices with integer coefficients and zeroes in the diagonal.
\end{corollary}

Implicit in the statement of the theorem is the fact that the individual components of the link are trivial knots.

\vskip.5in

\section{ Appendix V: 
Homotopy spheres embedded in codimension two 
 and  the Kervaire-Arf-Robertello-Levine invariant}

\subsection{Which homotopy spheres can be embedded in codimension two?}

Let $\Sigma^n$ be a homotopy sphere differentiably embedded in $S^{n+2}$.

To answer the question ``which homotopy spheres can be embedded in codimension two", let us first discuss the case $n \leq 4$. 
\\
If $n = 1$ or $2$, there is no problem since $\Sigma^n = S^n$. 
\\
If $n = 3$ we could invoke Perelman. But there is no need to do that, since Wall proved the lovely result  that any 3-manifold embeds differentiably in  $S^5$.  
\\
If $n = 4$ we have to argue a little bit. By Kervaire-Milnor's Theorem 6.6 in \cite{kemi63} we know that every homotopy 4-sphere $\Sigma^4$ bounds a contractible 5-manifold $\Delta^5$. Let us take the double of $\Delta^5$, i.e. two copies of $\Delta^5$ glued along their boundary. We thus get a homotopy 5-sphere $\Sigma^5$. But Kervaire-Milnor prove that any homotopy 5-sphere is diffeomorphic to $S^5$ (see p.504 in \cite{kemi63}). Therefore $\Sigma^4$ embeds in $S^5$ and a fortiori in $S^6$. 
\\
{\bf Note.} The arguments for $n = 3$ and $4$ are those of Kervaire in his paper about the group of a knot.

We consider now the set $G^n$ of orientation preserving diffeomorphism classes of oriented homotopy spheres $\Sigma^n$ which can be embedded in $S^{n+2}$. 

\begin{theorem}
Suppose that $n \geq 5$. Then:
\\
$1$) $G^n$ is a finite cyclic group, canonically isomorphic to ${\bf Z} / d_n$.
\\
$2$) If $n$ is even, $d_n = 1$. In other words $\Sigma^n = S^n$.
\\
$3$) Suppose that $n$ is odd, equal to $4k-1$, with $k \geq 2$. Then 
$$d_n = 2^{2k-2}(2^{2k-1} - 1) \mathrm{numerator} \left( \frac{4B_k}{k} \right)$$
The canonical generator $1 \in {\bf Z} / d_n$ is represented by the boundary of the $E_8$ plumbing.
\\
$4$) If $n = 4k+1$ with $k \geq 1$ the group $G^n$ is either trivial or  isomorphic to ${\bf Z} / 2$. More precisely:
\\
$4_1$) $G^n$ is trivial if $n = ~5~ ,~13~ , ~29~ , 61$ and possibly $125$ (these are the {\bf exceptional integers}). 
\\
$4_2$) $G^n$ is isomorphic to ${\bf Z} / 2$ if $n$ is not exceptional. The non-trivial element is represented by the boundary of the Kervaire plumbing. 
\end{theorem}

The integer $d_n$ when $n = 4k-1$ was called by Kervaire in a letter to Haefliger ``l'entier bord\'elique bien connu". The numerator of Bernoulli numbers is indeed a tough object, contrary to the denominator.

The items $1$, $2$ and $3$ are  an immediate consequence of Kervaire-Milnor, up to a factor 2, which was settled by the solution of the Adams Conjecture. The item $4$ is intimately related to the Kervaire invariant problem  solved in 2009.

Next proposition is the first step towards the proof of the theorem above.
 
\begin{proposition}
A homotopy n-sphere can be embedded differentiably in $S^{n+2}$ if and only if it bounds a parallelisable (n+1)-manifold.
\end{proposition}

{\bf Sketch of proof.} If $\Sigma^n$ is embedded in $S^{n+2}$ it bounds a Seifert hypersurface, which is parallelisable since it is orientable and embedded in codimension 1. (Of course this argument works with no restriction on $n$.)
\\
Conversely if $\Sigma^n$ bounds abstractly a parallelisable manifold, it also bounds a parallelisable handlebody if $n \geq 5$, by easy surgery below the middle dimension and by the recognition theorem (it is here that $(n+1) \geq 6$ is required). It is then easy to embed such handlebodies in codimension 1. This implies the following proposition.

\begin{proposition}
Let $n \geq 5$.  The group $G^n$ is isomorphic to Kervaire-Milnor's group  $bP^{n+1}$.
\end{proposition}

{\bf Question.}  If $\Sigma^{2q-1} \subset S^{2q+1}$ is an odd dimensional  knot,   can we detect by knot invariants which homotopy sphere is represented by $\Sigma^{2q-1}$ ?

The answer is yes, and the Seifert form does the work.
 \\
 If $n=4k-1=2q-1$ ( i.e. $q$ is even)  the answer is straightforward. Take the Seifert form of any Seifert hypersurface for $\Sigma^n \subset S^{n+2}$, symmetrize it, take the signature of this bilinear symmetric form (which is the intersection form of the chosen Seifert hypersurface) and divide it by $8$. The reduction modulo $d_n$ of the  integer such  obtained is the element of   ${\bf Z} / d_n$ represented by $\Sigma^n$.
 \\
 If $n = 4k+1=2q-1$ (i.e. $q$ is odd) the answer is provided by the Kervaire-Arf-Robertello-Levine invariant as we shall see now.
 
 \newpage

\subsection{The Kervaire-Arf-Robertello-Levine invariant }

KARL stands for Kervaire, Arf, Robertello and Levine. They are the main contributors.

The origin of the construction of the KARL invariant goes back to the thesis presented by Kervaire in Paris in June 1964, precisely with the purpose to answer the question raised just above. See \cite{kerv65} p.236. The study of the invariant was then developed by Robertello in his NYU thesis (Kervaire was the director) also in 1964. See \cite{robe65}.  Soon afterwards, Levine simplified in \cite{levi66}  Robertello's presentation of the invariant defined via Seifert forms. Actually,  Robertello gave in the case of classical knots two equivalent definitions. One of them can easily apply to higher dimensional  knots. The data is the following.

We have $n = 4k+1$ and   $K^{4k+1} \subset S^{4k+3}$ is a  $(4k+1)$-dimensional knot, where $K^{4k+1}$ is diffeomorphic to the homotopy sphere $\Sigma ^{4k+1}$ .  The case $n = 1$ (i.e. $k = 0$) is admitted. Let $F^{4k+2}$ be  a Seifert hypersurface for $K^{4k+1}$. 
\\
We denote by $H$ the quotient of $H_{2k+1} (F ;  {\bf Z})$ by its ${\bf Z}$-torsion subgroup. Let $\mathcal {A} : H  \times H \rightarrow {\bf Z}$ be the Seifert form  defined by $\mathcal {A} (x , y) = Lk(x ; i_+y)$.

\begin{definition}
We define the quadratic form $Q_F : H  \otimes {\bf Z} / 2 \rightarrow {\bf Z} / 2$ by $Q_F(x) = Lk_2(x ; i_+(x))$, where $Lk_2$ denotes the linking number reduced mod 2.
\end{definition}

\begin{lemma}
We have the equality $Q_F(x+y) = Q_F(x) + Q_F(y) + x \cdot y$ where $x \cdot y$ denotes the intersection number in $F$ reduced mod 2.
\end{lemma}

We denote by Arf$(Q_F)$ the Arf invariant of the quadratic form $Q_F$. Note that the bilinear form $(x , y) \mapsto x \cdot y$ is unimodular, since $\partial F$ is a homotopy sphere.

\begin{theorem}(Robertello; Levine \cite{levi66})
Arf$(Q_F)$ is a knot invariant, i.e. it does not depend on the choice of a Seifert hypersurface. 
\end{theorem}

{\bf Sketch of Levine's proof.} Let $p : \widehat {E(K)} \rightarrow E(K)$ be the infinite cyclic covering of the knot complement $E(K) = S^{4k+3} \setminus N(K)$.  Let $\Delta_K (t) \in {\bf Z} \lbrack t , t^{-1} \rbrack$ be the Alexander polynomial of the knot module $H_{2k+1} (\widehat {E(K)} ; {\bf Q})$ normalised \`a la Conway by $\Delta (t^{-1}) = \Delta (t)$ and $\Delta (1) = 1$. Then Levine proves that $\Delta (-1) \equiv 1 + 4$Arf$(Q_F)$ mod 8 (Levine's proof is a beauty). QED.

\begin{definition}
We define the KARL invariant KARL$(K)$ of the knot $K$ to be Arf$(Q_F)$ for any Seifert hypersurface $F$.
\end{definition}

\begin{theorem} (Robertello; Levine)
KARL(K) is a knot cobordism invariant.
\end{theorem}

{\bf Sketch of  proof, following Levine.} By Levine \cite{levi69}, the Fox-Milnor result is true in higher dimensions. More precisely,  if $K$ is cobordant to zero,  there exists an element $P(t) \in {\bf Z} \lbrack t , t^{-1} \rbrack$ such that $\Delta (t) = P(t) P(1 / t)$. Since $\Delta (1) = 1$ we deduce that $P(-1)$ is odd and hence $\Delta (-1)$ is an odd square and therefore congruent to $1$ mod 8. Hence KARL$(K) = 0$ if $K$ is cobordant to zero.

\subsection{The Hill-Hopkins-Ravenel result and its influence on the KARL invariant}

\begin{theorem}(Browder \cite{brow69} and Hill-Hopkins-Ravenel \cite{hihora09})
The ``boundary" homomorphism
 $$b : P^{4k+2} = {\bf Z} / 2 \rightarrow \Theta^{4k+1}$$
 is trivial if $4k+2 = 2 , 6 , 14 , 30 , 62$ and possibly $126$. (These are the exceptional integers +1).
\\
In all other cases it is injective and hence $bP^{4k+2} = {\bf Z} / 2$.
\end{theorem}

\begin{theorem}
If (4k+1) is not exceptional,  the Arf invariant of the quadratic form defined via the Seifert form (i.e. the KARL invariant of the knot) detects which homotopy sphere  of dimension $4k+1$ is embedded.
\end{theorem}

{\bf Proof of the theorem.}
Let us first suppose that the knot $K^{4k+1}$ is simple. Let $F^{4k+2}$ be a $2k$-connected Seifert hypersurface. Hence the quadratic form $H_{2k+1} (F ; {\bf Z} / 2) \rightarrow {\bf Z} / 2$  \`a la Kervaire-Milnor (see \cite{kosi93} p.202-207) is defined. 
\\
By \cite{keva67} Lemma on normal bundles p.512, this last quadratic form coincides with the quadratic form which gives rise to the KARL invariant. Hence KARL$(K)$  detects which element of $P^{4k+2}$ is represented by $F^{4k+2}$ and hence which element of $bP^{4k+2}$ is represented by $K^{4k+1}$.
\\
Let us now withdraw the assumption that $K$ is simple. We know that its KARL invariant is defined. By \cite{levi69} there exists a simple knot $K^*$ which is knot-cobordant to $K$.
By \cite{levi66} the knots $K$ and $K^*$ have the same KARL invariant. Since $K$ and $K^*$ are diffeomorphic, KARL$(K^*)$ = KARL$(K)$ detects which element of $bP^{4k+2}$ is represented by $K$. 
{\bf End of proof of the theorem.}

{\bf Recapitulation.} In terms of $bP^{4k+2}$ the interpretation of the KARL invariant splits in two cases. 

i) If $bP^{4k+2} = {\bf Z} / 2$ there are two homotopy spheres  of dimension $(4k+1)$ which can be embedded in $S^{4k+3}$ to represent  knots. The KARL invariant detects which homotopy sphere it is. Thus the differential structure on the homotopy sphere is detected by the Alexander module (in dimension (2k+1)) of the knot. In other words, different differential structures give rise to different knot modules in dimension $(2k+1)$.
\\
{\bf Caution.} The KARL invariant detects the Arf  invariant of the handlebody (i.e. of the Seifert hypersurface). If you know beforehand that the Kervaire invariant is non-trivial, then it detects which homotopy sphere  is in the boundary of the handlebody.

ii) If $bP^{4k+2} = 0$ all  $(4k+1)$-knots are represented by an embedding of the standard sphere.  From a knot theory point of view, one could argue that this situation is more interesting, since the KARL invariant is a genuine invariant of the embedding. The case of classical knots is typical.

An (easy) consequence of our discussion of the KARL invariant produces a proof of the following result, when $n = 4k+1$. If $n = 4k+3$ it can be obtained by using  the signature and Blanchfield duality. See Dieter Erle \cite{erle72}.
 
\begin{theorem}
Let $n \geq 5$.  Let $\Sigma^n$ and $\widehat{\Sigma}^n$ be homotopy spheres embedded in $S^{n+2}$. Suppose that the two homotopy spheres are not diffeomorphic. Then the two  knots are not topologically isotopic.
\end{theorem}

\newpage

%Voici comment faire la bibliographie%

%Majuscules dans les titres en anglais. Est-ce une bonne idee?

\end{document}